\documentclass[12pt]{article}
\title{Henstock Lectures on Integration Theory}
\author{Edited by P.~Muldowney}
\date{}

\newcommand{\bS}{\mathbf{S}}
\newcommand{\A}{\mathcal{A}}
\newcommand{\Q}{\mathcal{Q}}
\newcommand{\G}{\mathcal{G}}
\newcommand{\mcC}{\mathcal{C}}
\newcommand{\cS}{\mathcal{S}}
\newcommand{\vs}{\mathcal{VS}}
\newcommand{\VS}{\overline{\mathcal{VS}}}
\newcommand{\mbX}{\mathbf{X}}
\newcommand{\mcB}{\mathcal{B}}

\newtheorem{theorem}{Theorem}

\newtheorem{example}{Example}

\newcommand{\vt}{\vspace{5pt}\\}

\newcommand{\D}{\mathcal{D}}

\newcommand{\Pa}{\mathcal{P}}

\newcommand{\ve}{\varepsilon}

\newcommand{\T}{\mathcal{T}}

\newcommand{\proof}{\noindent\textbf{Proof. }\noindent}
\newcommand{\nproof}{\hfill$\mathbf{\bigcirc}\vspace{12pt}$ }

\date{}
\begin{document}
\maketitle

\section{Introduction}\label{Introduction}
These are the notes of lectures given by Ralph Henstock at the New University of Ulster in 1970--71. 

Sections 2 to 19 (pages 1--70) deal with the Riemann-complete (or generalized Riemann) integral, also known as the gauge integral, or Henstock-Kurzweil integral. 
These sections cover, essentially, the same ground as (\cite{H1}, 1963), the first book on the subject; but perhaps in a clearer and simpler style.

Robert Bartle's paper, \emph{Return to the Riemann integral} (\cite{RB}, 1980), is a good introduction
to the Riemann-complete integral.

The rest of these notes---Section 20 onwards---deal with Henstock's abstract or general theory of integration, which in \cite{MTRV} is called the {Henstock integral}. Originally mooted in \cite{H2} (1968), this general theory was still in a formative stage in 1970--71, and received fuller expression in  \cite{H4}, (1991). Also in MTRV (\cite{MTRV}, 2012).

$\;\;\;\;\;\;\;\;\;\;\;\;\;\;\;\;\;\;\;\;\;\;\;\;\;\;\;\;\;\;\;\;\;\;\;\;\;\;\;\;\;\;\;\;\;\;\;\;\;\;\;\;\;\;
$
--  \emph{P.~Muldowney, November 16 2015}

\section{Riemann Integration}\label{Riemann Integration}
We consider this integration in Euclidean space $E^n$ of $n$ dimensions, assuming that we are given $n$ co-ordinate\footnote{The symbol $P$ usually appears as an alternate for $x$, $y$, $\ldots$, generally denoting associated point or tag-point of interval. \emph{--P.~Muldowney}} axes at right angles, the axes $x_1, \ldots , x_n$. Then for each collection of $n$ pairs of numbers $a_1<b_1$, $a_2<b_2$,$\ldots$,$a_n<b_n$, we can define the \emph{brick} $I$ formed of all points $P$ with co-ordinates $x_1, \ldots , x_n$ satisfying
\[
a_1 \leq x_1 \leq b_1, \ldots , a_n \leq x_n \leq b_n.
\]
The volume of $I$ is defined by
\[
\mu(I) = (b_1-a_1) \cdots (b_n-a_n).
\]
The integration process depends upon cutting up $I$ into a finite number of smaller bricks $J$ and the result of this cutting up is called a \emph{division} $\D$.
We also call $I$ alone a division of $I$.


We are integrating a function $f(P) = f(x_1, \ldots , x_n)$, and to integrate by Riemann's original method we sum the terms $f(P)\mu(J)$,
for an arbitrary choice of points $P$ in $J$, one $P$ for each such $J$, and we write
\[
R(f;\D) = (\D) \sum f(P)\mu(J).
\]
We would expect that for ``smooth'' functions, as the bricks $J$ shrink in size,
\[
R(f;\D) \rightarrow \mbox{limit}.
\]
We define the \emph{diameter} of $I$ to be
\[
\mbox{diam} (I) = \sqrt{(b_1-a_1)^2 + \cdots + (b_n-a_n^2)}.
\]
The \emph{norm} of the division $\D$, norm$(\D)$, is defined to be the greatest diam$(J)$ for all $J$ in $\D$. Riemann's original definition is that
\[
R(f;\D) \rightarrow (\mbox{Riemann})\int_If(P)d\mu
\]
of $f(P)$ over the brick $I$ as norm$(\D) \rightarrow 0$.

More strictly (and holding even if $f$ is complex-valued), a number $r$ is the Riemann integral $(R)\int_I f(P)d\mu$ of $f(P)$ over the brick $I$ if, given $\ve>0$, there is a $\delta>0$ such that
\[
\left| R(f;\D) -r \right| < \ve
\]
for all divisions $\D$ of $I$ with norm$(\D) <\delta$ and for all choices of points $P$ in the bricks $J$ of $\D$. In this case we could have $f$ complex-valued or real.

But Darboux's modification has to assume $f$ to be real. Instead of using $f(P)$, Darboux uses
\[
M(f;J) = \mbox{l.u.b.} \left\{f(P):{P \in J}\right\},\;\;\;\;\;\;m(f;J) = \mbox{g.l.b.}\left\{f(P):{P \in J}\right\},
\]
and then considers
\[
S(f;\D) = (\D) \sum M(f;J)\mu(J),\;\;\;\;\;\;s(f;\D) = (\D) \sum m(f;J)\mu(J).
\]
The Riemann-Darboux integral is defined in the following way. It exists if and only if, given $\ve>0$, there is a $\delta>0$ such that for all divisions $\D$ with norm$(\D)< \delta$, we have
\[
\left|S(f;\D) - r \right| < \ve,\;\;\;\;\;\; \left|s(f;\D) - r \right| < \ve.
\]
\begin{theorem} Riemann's and Darboux's methods are equivalent for real-valued functions.
\end{theorem}
\proof
For each $P \in J$,
\[
m(f;J) \leq f(P) \leq M(f;J).
\]
Therefore
\[
s(f;\D) \leq R(f;\D) \leq S(f.\D).
\]
So if $s(f;\D)$ and $S(f;\D)$ tend to the same number $r$, then $R(f;\D)$ must also tend to $r$.
Conversely, let $\ve>0$. Then as $M(f;\D)$ is the lub of $f(P)$ in $J$ there is a point $Q$ in that range for which
\[
M(f;\D) - \ve < f(Q) \leq M(f;J).
\]
Taking such $Q$ for each $J$ in $\D$ we have
\[
S(f;\D) - \ve (\D) \sum \mu{(J)} < R(f;\D) \leq S(f;\D).
\]
Hence if $R(f;\D) \rightarrow$ limit, for all choices of $P$, we take $P=Q$ and this shows $S(f;\D) \rightarrow$ limit $r$. Similarly $s(f;\D) \rightarrow$ limit $r$.
\nproof

A corollary of this is that if the Riemann integral exists the function is bounded; otherwise either $m(f;J)=-\infty$ and then $s(f;\D)$ might not be defined, or $M(f;J) = +\infty$ and $S(f;\D)$ might not be defined, for some intervals $I$, and the Riemann-Darboux integral cannot exist.

\section{The Calculus Indefinite Integral}\label{The Calculus Indefinite Integral}
An earlier definition of an integral is due to Newton. If $f(x)$ has finite values in $b < x <c$ in one dimension, and if by some means we can find a function $F(x)$ continuous in $b \leq x\leq c$ such that $F'(x) = f(x)$ in $b <x<c$ then we say that the \emph{calculus integral} of $f(x)$ in $b \leq x \leq c$ is
\[
(\mbox{Cal})\int_b^c f(x)dx = F(c) - F(b) = F]_b^c.
\]
In practice Newton only integrated continuous functions and then it is easily proved that the calculus integral is the same as the Riemann integral. But about the year 1900 mathematicians decided to throw away Newton's limitation to continuous $f$ and take the definition as it stands. In about 1910 Denjoy showed that the newly-defined Lebesgue integral was not strong enough to integrate every function that was calculus-integrable and had to define his own integral to deal with the problem; e.g. if $F(x) = 2x^{\frac 12}$ in $x \geq 0$, then $F'(x) = x^{-\frac 12}$ in $x \geq 0$, and that is not integrable by Riemann's method in $0 \leq x \leq 1$ as it is unbounded. But (Cal)$\int_0^1 F'(x) dx = 2-0 =2$.
If
\[
F(x) = x^p \sin x^{-q}\;\;\;\mbox{ for }\;\;\;x\neq 0,\;\;\;\;\;\; F(0)=0,
\]
then
\[
F'(x) = px^{p-1} \sin x^{-q} -qx^{p-q-1} \cos x^{-q}\;\;\;(x \neq 0),
\]
while
\[
\frac{F(h) - F(0)}{h} = h^{p-1} \sin h^{-q} \rightarrow 0\;\;\;\mbox{ as }\;\;\;h \rightarrow 0 \;\;\;\mbox{ if }\;\;\;p>1.
\]
For instance, if $p=2$, then, for $x \neq 0$,
\[
F'(x) = 2x \sin x^{-q} - q x^{1-q} \cos x^{-q}.
\]
If, say, $q=3$ then the second term in $F'(x)$ is $-3 x^{-2} \cos x^{-3}$ and this is not integrable by Lebesgue's method, as $|3x^{-2} \cos x^{-3}|$ has an infinite integral over $0\leq x \leq 1$.

The exact definition of the calculus integral is as follows. Given $\ve>0$, a function $\delta(x)>0$ is defined in $b \leq x \leq c$:
\begin{eqnarray}
 \left|\frac{F(y)-F(x)}{y-x}\right| < \ve\;\;\;\mbox{ or }\;\;\; |F(y)-F(x)| < \ve,\;\;\;|f(x)(y-x) < \ve,\;\;\;b<x<c,
\label{(1)}
\vt
 |F(y)-F(x) |<\ve,\;\;\;|f(x)(y-x)| < \ve,\;\;\; x=b,c,\;\;\;F \mbox{ continuous at } b,c.\label{(2)}
\end{eqnarray}
If $ b<x<c$, $b \leq u < v \leq c$, $x-\delta(x) <u \leq x \leq v <x+\delta (x)$, (\ref{(1)}) gives
\begin{equation}\label{(3)}
\begin{array} {rll}
|F(v)-F(u) -f(x)(v-u) & \leq & |F(v)-F(x) - f(x) (v-x)|\;\;+  \vt
&& \;\;\;\;\;\;+\;\; |F(x)-F(u) - f(x) (x-u)| \vt
&\leq & \ve (v-x) + \ve (x-u) \;\;\;\;=\;\;\;\; \ve(v-u).
\end{array}
\end{equation}
\begin{theorem} \label{theorem 2}
If a finite number of closed intervals $[u,v]$ of the type (\ref{(3)}) form a division $\D$ of $b \leq x \leq c$, then for the various $x$ and $u,v$,
\[
\left|(\D) \sum f(x)(v-u) -(F(c) - F(b))\right| < \ve (c-b+4).
\]
\end{theorem}
\proof
If the division $\D$ is given by $b=t_0<t_1< \cdots < f_m =c$ where each $[t_{j-1}, t_j]$ is a $[u,v]$ and if the corresponding $x$ is $\xi_j$, then
\[
\begin{array}{rll}
&&\left|(\D)\sum f(x)(v-u) - (F(c) - F(b))\right| \vt
&=&
\left| \sum_{j=1}^m \left(F(\xi_j)(t_j-t_{j-1}) -F(t_j) + F(t_{j-1})\right)\right| \vt
&\leq & \sum_{j=1}^m \ve (t_j - t_{j-1} + 2 \ve + 2 \ve \;\;\; =\;\;\;\ve(c-b+4).
\end{array}
\]
\nproof

\noindent
This shows that if we can arrange our divisions to satisfy the conditions involving the $\delta(x)>0$, for all such $\delta(x)>0$, we can consider a limit process in which $\ve \rightarrow 0$ and so
\[
(\D) \sum f(x)(v-u) \rightarrow F(c) - F(b).
\]
\begin {theorem} \label{theorem 3}
Let $b<c$ and let $\delta(x)>0$ be a positive function defined in the closed interval $b \leq x \leq c$. Then there is a finite number of points $b \leq x_1 < x_2 < \cdots < x_n \leq c$ such that each $x$ in $b \leq x \leq c$ lies in at least one of the open intervals $x_j - \delta(x_j) < u < x_j + \delta (x_j)$. A corollary of this is that we can arrange that each point of $[b,c]$ lies in at most two of the finite number of open intervals.
\end{theorem}
\proof
Let $\Delta(x)$ be the interval $x-\delta(x) <t<x+\delta(x)$. Then for $x=u,v,w$ with $u<v<w$, let the $\Delta(x)$ have a common point. If $v-\delta(v) \leq u - \delta(u)$ then
\[
\delta(v) \geq \delta(u) + v-u > \delta(u),\;\;\;\;\;\; v+\delta(v) > v+\delta(u) > u+\delta(u),\;\;\;\;\;\; \Delta(v) \supset \Delta(u),
\]
and we can omit $\Delta(u)$. Similarly, if $v+\delta (v) \geq w + \delta (w)$ then $\Delta (v) \supset \Delta (w)$, and we can omit $\Delta (w)$. Thus we can assume
\[
u-\delta (u) < v- \delta(v) < v + \delta (v) < w + \delta(w).
\]
As all these intervals have a common point, $\Delta(v) \subseteq \Delta(u) \cup \Delta (w)$
and we can omit $\Delta(v)$. Thus if three intervals have a common point we can omit at least one. Eliminating the finite number of intersections of intervals in this way, we obtain the corollary.
Finally, to obtain the required divisions $\D$ we use the corollary and have $\Delta(\xi_j)$ intervals, $1 \leq j \leq m$, where $b \leq \xi_1 < \cdots < \xi_m \leq c$, and where we can have $\xi_1=b,$, $\xi_m=c$. If
\[
t_1 \in \Delta(\xi_1) \cap \Delta(\xi_2),\;\;\;\;t_2 \in \Delta(\xi_2) \cap \Delta(\xi_3), \ldots ,
\]
we put $t_0 = b$, $t_m =c$, and then
\[
[t_{j-1},t_j] \subset \Delta(\xi_j),\;\;\;\; t_{j-1} \leq \xi_j \leq t_j.
\]
Hence in one dimension we can find divisions $\D$ of $[b,c]$ that are \emph{compatible with} $\delta(x)>0$.
\nproof

\section{The Generalised Riemann Integral}\label{The Generalised Riemann Integral}
Theorem 2 shows that if $n=1$ and if we restrict the divisions of $[b,c]$ forming them from intervals as in (\ref{(2)}), (\ref{(3)}), and if we have a limiting process that uses smaller and smaller $\delta(x)>0$ so that we can let $\ve \rightarrow 0$, then we can integrate the derivative $f(x) = F'(x)$.

Let $\delta(P)>0$ be a positive function of points $P$ in $n$-dimensional space, and let $S(P, \delta(P))$ denote the open sphere (ball) with centre $P$ and radius $\delta(P)$. So if $P=(p_1, \ldots ,p_n)$ the ball (\emph{open ball}) is the region of all points $Q=(q_1, \ldots , q_n)$ with
\[
\sum_{j=1}^n \left(q_j - p_j\right)^2 < \delta(P)^2.
\]
If $I$ is a \emph{brick} $a_j \leq x_j \leq b_j$, $1 \leq j \leq n$, with $b_j>a_j$, and if $\D$ is a division of the brick $I$ that has a point $P$ associated with each brick $J$ of $\D$ such that $P \in J$, and $J \subseteq S(P, \delta(P))$, we say that $\D$ is \emph{compatible with} $\delta(P)$. We call $P$ an \emph{associated point} of $J$.

\begin{theorem}
Given $\delta(P)>0$ defined in a brick $I$, there is a division $\D$ of $I$ that is compatible with $\delta(P)$.
\end{theorem}
\proof
Suppose false, i.e.~there does not exist such a division. We bisect $I$ in the direction of each of the co-ordinate hyperplanes to obtain $2^n$ smaller bricks. (If the brick is given by $a_j \leq x_j \leq b_j$, we use $\frac 12 (a_j +b_j)$, $1 \leq j \leq n$.)
If each smaller brick has a division compatible with $\delta(P)$ we could put the divisions together to form a division of $I$. Thus at least one of the smaller bricks has no division compatible with $\delta(P)$. Take one such brick and call it $I_1$, and repeat the construction to obtain $I_2$, etc. We thus obtain an infinite sequence of bricks in $I$. Taking the point $P_k$ whose coordinates are $p_j = \min \{x_j\}$ for all $(x_1, \ldots , x_n) \in I_k$, the sequence $\{P_k\}$ is monotone increasing in each co-ordinate. Also each $P_k$ lies in $I$, so that each co-ordinate $p_{kj}$ is bounded above. Therefore each $p_{kj} \rightarrow $ limit $\pi_j$ as $k \rightarrow \infty$. Call the limit $P$. Then $\delta(P>0$ so that, for some $k$,
\[
\mbox{norm}(I_k)  = 2^{-k} \mbox{norm}(I) < \delta(P).
\]
therefore $I_k \subset S(P, \delta(P))$ and $P \in I_k$. Therefore $I_k$ forms a division of itself that is compatible with $\delta(P)$, giving a contradiction.
\nproof

Having shown the possibility of the construction of suitable divisions we can now define the (generalised Riemann) integral of a function $f(P)$ of points $P$ over a brick $I$, with real or complex values. We say that $f$ is \emph{integrable} over $I$, with respect to the volume function $\mu$, and that $F$ is the \emph{integral} of $f$ over $I$, $F = \int_I fd\mu$, if
\begin{equation}\label{def of integral}
\left|(\D) \sum f(P) \mu(J) - F\right| < \ve
\end{equation}
for all divisions $\D$ of $I$ that are compatible with $\delta(P)$ (i.e.~ $P \in J \in \D$, $J \subseteq S(P, \delta(P))$).
In passing, it is just as easy to define the (generalied Riemann) integral $H=\int_I dh$ of a function $h(P,J)$ of a brick $J$ and its associated point $P$ over the brick $I$, simply by replacing (\ref{def of integral}) by
\begin{equation}\label{def of integral2}
\left|(\D) \sum h(P,J) - H\right| < \ve
\end{equation}
\textbf{Uniqueness:} For fixed $f, \mu, I$, let $F$ and $G$ have the property (\ref{def of integral}) for all $\ve>0$, i.e.~we can replace $F, \delta, \D$ by $G, \delta_1, \D_1$, respectively. We take
\[
\delta_2(P)=\min\{\delta(P), \delta_1(P)\}>0.
\]
Then each division $\D_2$ compatible with $\delta_2(P)$ is also compatible with $\delta(P)$ and $\delta_1(P)$. Hence $\D_2$ is a $\D$ and a $\D_1$, so
\[
\left|(\D_2) \sum f(P) \mu(J) - F\right| < \ve,\;\;\;\;\;\;\left|(\D_2) \sum f(P) \mu(J) - G\right| < \ve,
\]
and
\[
|F-G| = \left| F - (\D_2) \sum f(P)\mu(J) + (\D_2) \sum f(P)\mu(J) -G \right| < \ve + \ve.
\]
As $\ve$ is arbitrary we get $F=G$.
Similarly, if $H$ exists it is uniquely defined.

\section{Elementary Properties of the Integral}\label{Elementary Properties of the Integral}
We begin with a theorem on the geometry of the system.
\begin{theorem}\label{geometry}
Let bricks $I_j$ form a division of a brick $I$, and let $\delta_j(P)>0$ be defined in $I_j$ ($1 \leq j \leq m$). Then there is a $\delta(P)>0$ defined in $I$ such that
$\delta(P) \leq \delta_j(P)$ ($P \in I_j$) for $1 \leq j \leq m$, and such that if $J \subseteq I$, with associated point $P$, compatible with $\delta(P)$, and does not lie in any $I_j$ but overlaps with at least two, then $P$ is on the boundary of two or more of the $I_j$, and $J$ can be divided up into bricks $J_1, \ldots , J_k$ such that each $J_k$ lies in an $I_j$, has associated point $P \in I_j$, and os compatible with $\delta_j(P)$.
\end{theorem}
\proof
When $P \in I_j$, $P$ not on the boundary $B_j$, let $\rho(B_j,P)$ be the least distance from $P$ to $B_j$ and put
\[
\delta(P) = \min\left\{ \delta_j(P), \frac 12 \rho (B_j, P)\right\}.
\]
If $P \in B_j$ then $P$ might be on several boundaries. Let $d>0$ be the distance from $P$ to the nearest vertex (or nearest other vertex) of the $B_k$ that contain $P$. Then we define
\[
\delta(P) = \min \left\{\frac 1 2 d, \delta_k(P) \right\}
\]
for all $k$ with $P \in B_k$. By constrction, if a brick $J \subseteq I$, with associated point $P$, is compatible wth $\delta(P)$ and overlaps with several $I_j$, then $P$ lies  on a $B_j$ and that part of $B_j$ on which $P$ lies in  part of an $(n-1)$-dimensional hyperplane that cuts across $J$. Hence the result.
\nproof

The function $h(J)$ of bricks $J$ is additive in $I$ if, for each $I_1 \subset I$ and each division of $\D$ of $I$, $(\D)\sum h(J) = h(I_1)$. For example, the volume has this property.
\begin{theorem} \label{additivity of integral}
Let bricks $I_j$ form a division of a brick $I$, such that $f(P)$ is integrable over $I_j$, $1 \leq j \leq m$. Then $f(P)$ is integrable over $I$, and the integral is additive over divisions of $I$. More generally, a similar result holds for a function $h(P,J)$ of bricks $J$ and associated points $P$ if $h$ is additive with respect to $J$ or $h$ is integrable in $I$.
\end{theorem}
\proof
Given $\delta_j(P)>0$ in $I_j$, $1 \leq j \leq n$, we construct the $\delta(P)$ of Theorem \ref{geometry} so that a division $\D$ over $I$ and compatible with $\delta(P)$ has the properties given in Theorem \ref{geometry}. In particular if a $J$ of $\D$ overlaps with two or more $I_j$, we can cut up the $J$ into two or more bricks $J_1, \ldots , J_r$, each with associated point $P$, each lying in an $I_j$, each compatible with $\delta(P)$, and
\[
f(P)\mu(J)= f(P)\mu(J_1) + \cdots +f(P)\mu(J_r).
\]
Thus we can reduce the sum $(\D)\sum f(P)\mu(J)$ to a sum in which there is no overlapping, and then we can separate the sum into $m$ parts, each part being the sum for a division of $I_j$ compatible with $\delta_j(P)$, $1 \leq j \leq m$. By choice of $\delta_j(P)$ we can ensure that each such part is within $\ve m^{-1}$ of the corresponding integral. Therefore
\[
|(\D)\sum f(P)\mu(J) -A| < \ve\;\;\;\mbox{ where }\;\;\;A= \int_{I_1}fd\mu + \cdots + \int_{I_m}fd\mu,
\]
and where $A$ is independent of $\D$. As $\ve$ is arbitrary, the integral of $f$ over $I$ exists and is equal to $A$. For the more general case just replace $f(P)\mu(J)$ by $h(P,J)$. If we already know that $h(P,J)$ is integrable over $I$, we need only choose the original $\D$ so that there is no overlapping.
\nproof

\begin{theorem} \label{integrability on sub-interval}
If $f(P)$ is integrable in $I$, and if the brick $J\subset I$, then $f(P)$ is integrable in $J$. A similar result holds for $h(P,J)$.
\end{theorem}
\proof
If $\delta(P)>0$ is such that the sums over every division of $I$ compatible with $\delta(P)$ are within $\ve>0$ of the integral of $f$ over $I$, then two such sums differ by at most $2\ve$. Now let $s_1,s_2$ be the values of two sums for divisions of $J$ compatible with $\delta(P)$. By extending the sides of $J$ to cut the sides of $I$, the part of $I$ lying outside of $J$ can be cut up into one or more bricks $I_1, \ldots , I_r$. There is a sum $s_J^*$ over a division of $I_j$ compatible with $\delta(P)$ for $1 \leq j \leq r$ so that
\[
s_1+ s_1^*+ \cdots + s_r^*\;\;\;\;\;\mbox{ and }\;\;\;\;\;s_2+ s_1^*+ \cdots + s_r^*
\]
are two sums over $I$ for divisions compatible with $\delta(P)$. Hence $s_1 - s_2$ satisfies
\begin{equation}\label{cauchy cond}
|s_1-s_2| < 2 \ve.
\end{equation}
Using a Cauchy sequence argument, take $\ve = m^{-1}$ and then write the $\delta(P)$ as $\delta_m(P)$. The function
\[
\delta_m^*(P) = \min\{ \delta_1(P), \ldots , \delta_m(P)\} >0
\]
has the same properties as $\delta_m(P)$, and is monotone decreasing in $m$. Thus there is no loss of generality in assuming that
\begin{equation} \label{cauchy cond 2}
\delta_m(P) \leq \delta_{m-1}(P) \;\;\;\;\;(m \geq 2).
\end{equation}
Let $s_m$ be a fixed sum and $s'_m$ an arbitrary sum for divisions of $J$ compatible with $\delta_m(P)$. By (\ref{cauchy cond 2}), $s_k$ is an $s_m$ when $k>m$, so that, by (\ref{cauchy cond})
\begin{equation}\label{cauchy cond 3}
\left|s_m - s'_m\right| < \frac 2 m,\;\;\;\;\;\;\;\;\left|s_m - s_k\right| < \frac 2 m\;\;\;\;\;(k>m).
\end{equation}
Thus the sequence $\{s_m\}$ of special sums, one for each $m$, is a Cauchy sequence of real or complex numbers, and so converges to the value $s$, say. Therefore by (\ref{cauchy cond 3}), for all sums $s'_m$ over divisions of $J$ compatible with $\delta_m(P)$,
\[
\left|s_m'-s\right| \leq \left|s_m'-s_m\right| + \left|s_m-s\right|< \frac 4 m,
\]
and the integral of $f$ over $J$ exists and is equal to $s$.
\nproof

\begin{theorem} \label{Henstock Lemma}
Let $F(J)$ be the integral of f(P) over the brick $J \subseteq I$. Given $\ve>0$, let $\delta(P)>0$ be defined in $I$ so that $|s-F(I)|<\ve$ for all sums $s$ over $I$ compatible with $\delta(P)$. Suppose $p$ is a sum of terms $\{f(P_j)\mu(J_j) - F(J_j)\}$ for any number of distinct bricks $J_j$ of a particular division $\D$ of $I$ compatible with $\delta(P)$, together with the associated points $P_j$. Then
\begin{equation}\label{HL1}
|p| \leq \ve,
\end{equation}
and
\begin{equation}\label{HL2}
(\D) \sum |f(P)\mu(J) -F(J)| \leq 4\ve.
\end{equation}
A similar result holds for $h(P,J)$ when it is integrable in $I$.
\end{theorem}
\proof
We number the bricks of $\D$ so that $J_1, \ldots , J_k$ are the bricks used for $p$, and $J_{k+1}, \ldots , J_m$ are the rest. For each $j$ in $k < j \leq m$, there is a sum $s_j$ for a division of $J_j$, compatible with $\delta(P)$, that is as near as we please to $F(J_j)$. If $q=\sum_{j=1}^k f(P_j) \mu(J_j)$, then $q+s_{k+1} + \cdots + s_m$ is a sum over a division of $I$ compatible with $\delta(P)$, so
\[
|q+s_{k+1} + \cdots + s_m -F(I)|< \ve,
\;\;\;\;\;\;
\left| q + F(I_{k+1}) + \cdots + F(I_m) -F(I) \right| \leq \ve.
\]
By the additivity of the integral,
\[
|p| = \left| \sum_{j=1}^k \left(f(P_j)\mu(J_j) - F(I_j)\right) \right| = \left| q + F(J_{k+1}) + \cdots + F(J_m) -F(I)\right| \leq \ve.
\]
If $p$ is complex then the real and imaginary parts have modulus $\leq \ve$. Taking in turn the set of bricks of $\D$ for which the real part of $f(P_j)\mu(J_j) - F(J_j)$ is positive, then those for which it is negative, we have
\[
(\D) \sum \left| \Re \left(f(P_j) \mu(J_j) -F(J_j) \right)\right| \leq 2\ve.
\]
Similarly for the imaginary part, and hence (\ref{HL2}).
\nproof
\begin{theorem}\label{HL converse}
Let $F(J)$ be additive for bricks $J \subseteq I$ and satisfy (\ref{HL2}) for every division $\D$ of $I$ compatible with $\delta(P>0$, where $\delta(P)$ depends on an arbitrarily small $\ve>0$. Then $f(P)$ is integrable in $I$, the integral over each brick $K \subseteq I$ being equal to $F(K)$. Similarly for $h(P,J)$.
\end{theorem}
\proof
Let $\D$ be a division of $K$ compatible with $\delta(P)$. Extending the sides of $K$ to meet the sides of $I$, we can cut up the closure of $I \setminus K$ into a finite number of bricks over which we can take divisions  compatible with $\delta(P)$, giving a division $\D_1$ of $I$ compatible with $\delta(P)$. (Thus $\D$ is a \emph{partial division} of $\D_1$.) Then
\[
\begin{array}{rll}
|(\D)\sum f(P)\mu(J) - F(K)| &=& |(\D)\sum \left(f(P)\mu(J) - F(J)\right)| \vt
&\leq & (\D) \sum |f(P)\mu(J) - F(J)| \vt
&\leq & (\D_1)\sum |f(P)\mu(J) - F(J) | \;\;\;\leq \;\;\;4\ve.
\end{array}
\]
\nproof

\begin{theorem}\label{Levi}(\textbf{Levi monotone convergence theorem}) Let $\{f_r(P)\}$ be a sequence of real functions, each integrable in $I$.
For each $P$ let $\{f_r(P)\}$ be monotone increasing and convergent to the value $f(P)$. If $F_r(I) = \int_I f_r(P)d\mu$  has a finite least upper bound $F$ as $r$ varies, then $f(P)$ is integrable in $I$ with integral equal to $F$. A similar result holds if $\mu$ is replaced by $h(J) \geq 0$ where $h(J)$ is additive.
\end{theorem}
\proof
For each integer $r$ there exists $\delta_r(P)>0$ such that if $\D_r$ is an arbitrary division of $I$ compatible with $\delta_r(P)$,
\begin{equation}\label{LMCT1}
\left|(\D_r)\sum f_r(P) \mu(J) - F_r(I) \right| < \frac \ve{2^r}.
\end{equation}
Further, as $\{f_r(P)\}$ is monotone increasing in $r$ for each $P$, $F_r$ is monotone increasing, and so tends to $F$, its l.u.b.~(\emph{least upper bound, supremum, $\sup$}). Thus there exists integer $r_0$ such that
\begin{equation}\label{LMCT2}
|F-r(I) - F| < \ve\;\;\;\;\;(r \geq r_0).
\end{equation}
Also, for each $P \in I$, there exists integer $r=r(P, \ve) \geq r_0$ such that
\begin{eqnarray}
\left| f_r(P) - f(P) \right| < \ve, \label{LMCT3} \vt
\delta (P) \equiv \delta_r(P, \ve) >0. \label{LMCT4}
\end{eqnarray}
By (\ref{LMCT3}), if $\D$ is a division of $I$ compatible with $\delta(P)$,
\begin{eqnarray} \label{LMCT5}
&&\left|(\D) \sum f(P) \mu(J) - (\D)\sum f_{r(P, \ve)}(P) \mu(J) \right| \nonumber \vt
 & \leq &(\D) \sum \left|f(P) \mu(J) - (\D)\sum f_{r(P, \ve)}(P) \mu(J) \right| \nonumber \vt
&<& (\D)\sum \ve \mu(J) \;\;\;=\;\;\; \ve \mu(I).
\end{eqnarray}
Grouping the $f_{r(P, \ve)}(P) \mu(J)$ into brackets with equal values of $r(P, \ve)$, and using (\ref{LMCT1}) and Theorem \ref{Henstock Lemma},
\begin{equation}\label{LMCT6}
\left|(\D) \sum f_{r(P, \ve)}(P) \mu(J) - (\D)\sum F_{r(P, \ve)}(J)  \right| \leq \sum_{r=1}^m \frac \ve{2^r} < \ve.
\end{equation}
Now $F_r(J) \leq F_{r+1}(J)$ since $f_r(P) \leq f_{r+1}(P)$, each $P$. Let $g$ and $h$ be the least and largest integers in the finite collection of values of $r(P, \ve)$. Then by (\ref{LMCT2}),
\begin{eqnarray*}
F- \ve &<& F_g(I) \;\;\;=\;\;\; (\D) \sum F_g(J) \;\;\;\leq \;\;\; (\D) \sum F_{r(P, \ve)} (J) \vt
&\leq & (\D) \sum F_h(J) \;\;\;=\;\;\;F_h(I) \;\;\;\leq \;\;\;F.
\end{eqnarray*}
Therefore
\[
F-\ve\left(\mu(I)+2\right) \;\;\;< \;\;\;(\D) \sum f(P) \mu(J) \;\;\;\leq\;\;\; F + \ve\left(\mu(I) +1\right).
\]
As $\ve>0$ is arbitrary, we have proved the result.
\nproof

If $X$ is a set of points in $n$-dimensional space let $\chi(X,P)$ be its characteristic function, i.e.,
\[
\chi(X,P) = \left\{
\begin{array}{lll}
1 & \mbox{if} & P \in X,\vt
0 & \mbox{if} & P \notin X.
\end{array}
\right.
\]
Let $ \delta(P)>0$ in a brick $I$ and let
\begin{equation}\label{def var 1}
V(\chi(X,\cdot)\mu;I;\delta):= \sup\{(\D) \sum \chi(X,P) \mu(J)\},
\end{equation}
the  supremum being taken over all divisions $\D$  of $I$ { compatible with } $\delta(P)$.
Then
\begin{equation}\label{def var 2}
0 \leq V(\chi(X,\cdot)\mu;I;\delta) \leq \mu(I).
\end{equation}
If $0 < \delta_1(P) \leq \delta(P)$ at each $P \in I$, then every division of $I$ compatible with $\delta_1(P)$ is also compatible with $\delta(P)$, so
\begin{equation}\label{def var 3}
V(\chi(X, \cdot)\mu;I;\delta_1) \;\;\leq \;\;V(\chi(X, \cdot)\mu;I;\delta).
\end{equation}
The g.l.b.\footnote{{Greatest lower bound, infimum $\inf$}}~of $V(\chi(X,\cdot)\mu;I;\delta)$, for all $\delta(P)>0$ is written as
\[
V(\chi(X,\cdot)\mu;I) = V(\mu;I;X) ,\;\;\;\;\;\;(=\mu^*(X\cap I)\;\mbox{ in Lebesgue case,})
\]
and is called the \emph{outer measure} of $X$ in $I$. (\ref{def var 3}) shows that the g.l.b. is the limit as $\delta(P)$ shrinks. If $\chi(X,P)$ is integrable in $I$ to the value $K$ then, for suitable $\delta(P)>0$,
\begin{eqnarray}
K-\ve &<& (\D)\sum \chi(X,P) \mu(J)\;\;<\;\;K+\ve, \nonumber \vt
K-\ve &<& V( \chi(X,P) \mu(J);I;\delta)\;\;<\;\;K+\ve, \nonumber \vt
V(\mu;I;X) &=& \mu^*(X \cap I) \;\;=\;\;\int_I \chi(X,P)d\mu. \label{def var 4}
\end{eqnarray}
In this case we say that $X$ is \emph{measurable} in $I$, or that $X\cap I$ is measurable, and we call $\mu^*(X \cap I)$ the \emph{measure} $\mu(X \cap I)$ of $X$ in $I$. In particular this happens when $\mu^*(X \cap I)=0$, for then there exists $\delta(P)$ so that
\[
0 \leq (\D) \sum \chi(X,P)\mu(J) \leq V(\chi(X, \cdot)\mu;I;\delta) < \ve,
\]
i.e.~$\int_I\chi(X,P)d\mu$ exists and equals $0$.

More generally, let $h(P,J)$ be a function of bricks $J$ and their associated points $P$.
Let
\begin{equation}\label{def var a}
V(h(P,J):I:\delta)\;\;:=\;\;\sup\{(\D)\sum |h(P,J)|,
\end{equation}
the supremum being taken over all divisions $\D$ of $I$ compatible with $\delta(P)$. Sometimes the sums are unbounded, in which case the supremum is written as $+\infty$. However if $V$ is finite for some $\delta(P)>0$ then $V$ is monotone decreasing as $\delta(P)$ decreases at each $P$. Thus we use the g.l.b.~of $V(h;I;\delta)$ over all $\delta(P)>0$ and we call it the \emph{variation} $V(h(P,J);I)$ of $h(P,J)$ in $I$, and we say that $h(P,J)$ is of \emph{bounded variation} in $I$. On the other hand, if $V(h;I;\delta) = +\infty$ for all $\delta(P)>0$, we write symbolically $V(h(P,J);I)=+\infty$. Replacing $h(P,J)$ by $\chi(X,P)h(P,J)$ we can define
\begin{eqnarray*}
V(h(P,J);I;X;\delta) &=& V(\chi(X,P)h(P,J);I;\delta),\vt
V(h(P,J);I;X)&=& V(\chi(X,P)h(P,J);I).
\end{eqnarray*}
In the definition of the integral of $f$ over $I$, using $\mu$, we can replace $f(P)\mu(J)$ by $|h(P,J)|$ and then we have defined the integral of $|h(P,J)|$. If it exists we have a result like (\ref{def var 4}). As a special case, if $f(P)\geq 0$ is integrable then
\begin{equation}\label{def var b}
\int_If(P)d\mu = V(f\mu;I).
\end{equation}

\begin{theorem}\label{properties of variation}

\begin{enumerate}
\item
If $X \subseteq X_1$ then $V(h(P,J);I;X) \leq V(h(P,J);I;X_1)$.
\item
If $I \subseteq I_1$ then $V(h(P,J);I) \leq V(h(P,J);I_1)$.
\item
If $X = \bigcup_{j=1}^\infty X_j$ then $V(h(P,J);I;X) \leq \sum_{j=1}^\infty V(h(P,J);I;X_j)$. \newline\emph{Corresponding results for $\mu^*$ are:}
\item
If $X \subseteq X_1$ then $\mu^*(X \cap I) \leq \mu^*(X_1 \cap I)$.
\item
If $X = \bigcup_{j=1}^\infty X_j$ then $\mu^*(X \cap I) \leq \sum_{j=1}^\infty \mu^*(X_j \cap I)$.
\end{enumerate}
\end{theorem}
\proof
For 1 we can use $\chi(X,P) \leq \chi(X_1,P)$. For 2 the bricks used to cover $I$ can be used to cover $I_1$ if we use extra bricks to cover the closure of $I_1 \setminus I$. Using 1 in 3, we can assume that the $X_j$ are mutually disjoint (i.e.~ no pair has a point in common). If the sum in 3 is $+\infty$ there is nothing to prove. If the sum is finite then, given $\ve>0$ and an integer $j$, there is a $\delta(P)>0$ such that
\begin{equation} \label{properties of variation inequality}
V(h(P,J);I; X_j; \delta_j) < V(h(P,J);I; X_j)  + \ve 2^{-j}.
\end{equation}
Let $\delta(P) =1$ outside of $X$ and $=\delta_j(P)$ in $X_j$, $j=1,2,3, \ldots$, the $X_j$ being mutually disjoint. Then if $\D$ is a division compatible with $\delta(P)$ we group together the non-zero terms of $(\D)\sum \chi(X,P)|h(P,J)|$ according to the $j$ for which $P \in X_j$ and we find
\begin{eqnarray*}
(\D)\sum \chi(X,P)|h(P,J)| & \leq & \sum_{j=1}^\infty V(h(P,J);I;X_j;\delta) \vt
&=& \sum_{j=1}^\infty V(h(P,J);I;X_j;\delta_j) \vt
&<& \sum_{j=1}^\infty V(h:I;X_j) + \ve,
\end{eqnarray*}
giving 3.
\nproof

We say that a result, which depends on points $P$, is true \emph{almost everywhere} in $I$ if it is true except for the $P$ in a set $X$ for which $\mu(X \cap I)=0$.

\begin{theorem}\label{variation in null sets}
If $V(h(P,j);I)=0$ and if $f(P)$ is a function of points $P$. Then $V(f(P) h(P,J);I) =0$.
Conversely, if $V(f(P) h(P,J);I) =0$ and if $f \neq 0$ in $X$, then $V(f(P) h(P,J);I;X) =0$.
\end{theorem}
\proof
Let $X_j$ be the set where $|f(P)| \leq j$, $j = 1,2,3, \ldots$. Then $I = \bigcup_{j=1}^\infty X_j \cap I$.
By Theorem \ref{properties of variation}, 3,
\begin{eqnarray*}
V(f(P)h(P,J);I) & \leq & \sum_{j=1}^\infty V(f(P)h(P,J);I; X_j) \vt
& \leq & \sum_{j=1}^\infty j V(h(P,J);I; X_j)\vt
&\leq & \sum_{j=1}^\infty jV(h(P,J);I) \;\;\;=\;\;\;0.
\end{eqnarray*}
For the second result we replace $f(P)$ in the first result by $\chi(X,P) (f(P))^{-1}$, and $h(P,J)$ by $f(P) h(P,J)$.
\nproof

\begin{theorem}\label{13: measure integral null set}
If $\mu^*(X \cap I) =0$ and if $f(P)$ is a function of points $P$, then
\[
\int_I f(P) \chi(X,P) d\mu =0 = \int_I |f(P)| \chi(X,P) d\mu =0.
\]
Conversely, if $\int_I |f(P)|d\mu =0$, then $f(P)=0$ except in a set $X$ with $\mu^*(X \cap I) =0$.
\end{theorem}
\proof
$|(\D) \sum f(P)\mu(J) -0| \leq V(f\mu;I;\delta)$ for all divisions $\D$ of $I$ compatible with $\delta(P)$. By the first part of Theorem \ref{variation in null sets}, $V(f\mu;I)=0$, so $V(f\mu;I;\delta)$ can be made arbitrarily small by choice of $\delta$. Therefore $f\chi(X,P)$ and $|f|\chi(X,P)$ are both integrable to $0$ in $I$. Conversely
\[
(\D)\sum |f(P)| \mu(J) = |(\D) \sum |f(P)| \mu(J) -0| < \ve
\]
for suitable $\delta(P)>0$ and all divisions $\D$ of $I$ compatible with $\delta(P)$. Therefore
\[
V(f\mu;I) \leq V(f\mu;I;\delta) \leq \ve\;\;\;\mbox{ and }\;\;\; V(f \mu;I)=0.
\]
By the second part of Theorem \ref{variation in null sets}, if $f(P) \neq 0$ in $X$ then $\mu^*(X \cap I) =0$.
\nproof

\begin{theorem} \label{14: additivity of V}
Let $I_1, \ldots , I_m$ form a division of $I$. Then
\begin{enumerate}
\item
$\sum_{j=1}^m V(h(P,J);I_j) \leq V(h(P,J);I)$,
\item
If, for each fixed $P$, $h(P,J)$ is additive in $J$, then there  is equality in 1. (E.g., $h(P,J) = f(P)\mu(J)$ is additive when $P$ is fixed; also $h(P,J) = f(P) \mu(J) - F(J)$ where $F$ is the integral of $f$.)
\end{enumerate}
\end{theorem}
\proof
For 1, let $\delta(P)>0$ be defined in $I$ so that
\begin{equation}\label{c}
V(h(P,J);I;\delta) < V(h(P,J);I) + \ve.
\end{equation}
Then there is a division $\D_j$ of $I_j$ compatible with $\delta(P)$ for which
\begin{equation} \label{d}
(\D_j)\sum |h(P,J)| < V(h(P,J);I_j ;\delta) - \frac{\ve}{ m},\;\;\;(j = 1, \ldots , m).
\end{equation}
The bricks of $\D_1, \ldots , \D_m$ form a division $\D$ of $I$ which is compatible with $\delta(P)$, so that, by (\ref{c}) and (\ref{d}),
\begin{eqnarray*}
\sum_{j=1}^m V(h(P,J);I_j) - \ve & \leq & \sum_{j=1}^m \left(V(h(P,J);I_j;\delta) -\frac \ve m\right) \vt
&<& \sum_{j=1}^m (\D_j)\sum |h(P,J)| \;\;\;=\;\;\;(\D) \sum |h(P,J)| \vt
& \leq & V(h(P,J);I;\delta) \;\;\;<\;\;\;V(h(P,J);I) + \ve.
\end{eqnarray*}
For 2, there is a $\delta_j(P)>0$ with
\begin{equation} \label{e}
V(h(P,J);I;\delta_j) < V(h(P,J);I_j) + \frac \ve m.
\end{equation}
By Theorem \ref{geometry} we can define a suitable $\delta(P)>0$ from the separate $\delta_j(P)>0$. If $\D$ is a division of $I$ compatible with $\delta(P)$, and if a brick $J$ of $\D$ overlaps with two or more $I_j$ we can cut up $J$ into $J_1, \ldots , J_r$ as before, and this time we have
\[
|h(P,J)| = |\sum_{j=1}^r h(P,J_j)| \leq \sum_{j=1}^r |h(P, J_j)|
\]
by the additivity of of $h(P,J)$ in $J$. Thus the sum over $\D$ when there is overlapping is $\leq$ the corresponding sum when the overlapping is removed. Thus, as we are taking the supremum, we need only assume the non-overlapping case, and the sum for $\D$ can be separated into $m$ parts, each part being a sum for a division  $\D_j$ of $I_j$ compatible with $\delta_j(P)$, $j=1, \dots , m$. Hence, by (\ref{e}),
\begin{eqnarray*}
(\D) \sum |h(P,J)| &=& \sum_{j=1}^m (\D_j) \sum |h(P,J)| \;\;\;<\;\;\; \sum_{j=1}^m V(h(P,j);I_j) + \ve, \vt
V(h(P,j);I) & \leq & V(h(P,J);I;\delta) \;\;\;\leq \;\;\;\sum_{j=1}^m V(h(P,J);I_j) + \ve.
\end{eqnarray*}
As $\ve>0$ is arbitrary, this gives the inequality opposite to 1, and so gives 2.
\nproof

\begin{theorem}\label{15}
Let $h(P,J)$ be additive in $J$ for fixed $P$, and let $\delta(P)>0$ give
\[
V(h(P,J);I;\delta) < V(h(P,J);I) + \ve.
\]
If $\D$ is a division of $I$ compatible with $\delta(P)$, and if $I_1, \ldots , I_r$ are a partial division of $I$ from $\D$ with associated points $P_1, \ldots , P_r$, then
\[
\sum_{j=1}^r |h(P_j,J_j)| < \sum_{j=1}^r V(h(P,J);I_j) + 2 \ve.
\]
\end{theorem}
\proof
If $I_{r+1}, \ldots , I_m$ are the other bricks in $\D$ we can choose divisions $\Delta_j$ of $I_j$ compatible with $\delta(P)$ so that
\[
|(\D_j) \sum h(P,J)| > V(h(P,J);I_j;\delta) - \frac \ve m \geq V(h(P,J);I_j) - \frac \ve m.
\]
Then $I_1, \ldots ,I_m, \D_{r+1}, \ldots , \D_m$ form a division of $I$ compatible with $\delta(P)$ and
\begin{eqnarray*}
&&\sum_{j=1}^r |h(P_j,I_j)| + \sum_{j=r+1}^m V(h(P,J);I_j)\vt
&<&
\sum_{j=1}^r |h(P_j,I_j) + \sum_{j=r+1}^m (\D_j) \sum|h(P,J)| + \ve
\vt
&\leq &
V(h(P,J);I;\delta) + \ve \;\;\;<\;\;\;V(h(P,J);I) + 2 \ve.
\end{eqnarray*}
We now use the additivity of $V$ (Theorem \ref{14: additivity of V}).

\begin{theorem}\label{16}
Let $X_1 \subseteq X_2 \subseteq \cdots \subseteq X_j \subseteq \cdots$ with $X= \bigcup_{j=1}^\infty X_j$. Then
\[
\mu^*(X \cap I) = \lim_{j \rightarrow \infty} \mu^*(X_j \cap I).
\]
More generally, if $h(P,J)$ is additive in $J$ for each fixed $P$ then
\[
V(h(P,J);I;X) = \lim_{j \rightarrow \infty} V(h(P,J);I;X_j).
\]
\end{theorem}
\proof
We need only prove the second result. By item 3 of Theorem \ref{properties of variation},
\[
V(h(P,J);I;X) \geq V(h(P,J);I;X_j)
\]
for all $j$, so
\[
V(h(P,J);I;X) \geq \lim_{j\rightarrow \infty}V(h(P,J);I;X_j).
\]
To prove the reverse inequality, let $\delta_j(P)>0$ in $I$ be such that
\[
V(h(P,J);I;X_j;\delta_j)    V(h(P,J);I;X_j)  + \frac \ve{2^j}.
\]
Let
\[
\delta(P) = \left\{
\begin{array}{lll}
1 & \mbox{if} & P \notin X,\vt
\delta_1(P)  & \mbox{if} &  P \in X_1, \vt
\delta_j(P)  & \mbox{if} & P \in X_j \setminus X_{j-1},\;\;j=2,3, \ldots.
\end{array}
\right.
\]
If $\D$ is a division of $I$ compatible with $\delta(P)$, then $\D$ splits up into partial divisions $\Pa_j$ of $I$ compatible with $\delta_j(P)$ with associated points in $X_j$.
Also there are bricks with associated points outside $X$. Thus by Theorems \ref{properties of variation}, \ref{14: additivity of V}, and \ref{15}, with $h(P,J)\chi(X_j,P)$ replacing $h(P,J)$,
\begin{eqnarray*}
(\D)\sum \chi(X,P)|h(P,J)| & = & \sum_{j=1}^m (\Pa_j) \sum |h(P,J)| \vt
&<& \sum_{j=1}^m \left((\Pa_j)\sum V(h(P,J);J;X_j) + \frac{2\ve}{2^j}\right) \vt
&<& \sum_{j=1}^m (\Pa_j)\sum V(h(P,J);J;X_m) + 2 \ve \vt
& \leq & V(h(P,J) ;I;X_m) + 2 \ve \vt
& \leq & \lim_{j \rightarrow \infty} V(h(P,J) ;I;X_j) + 2 \ve, \vt
V(h(P,J);I;X) & \leq & V(h(P,J);I;X;\delta) \vt
 & \leq & \lim_{j \rightarrow \infty}V(h(P,J);I;X_j) + 2 \ve,
 \end{eqnarray*}
 hence the result.
 \nproof

 We now return to Lebesgue's monotone convergence theorem.
\begin{theorem}\label{17}
\begin{enumerate}
\item
Let $f(y,P)$ depend on a parameter $y$ that takes all values in $y \geq 0$, or all positive integer values, such that $f(y,P)$
is monotone increasing in $y$ for each fixed $P$. If, for each fixed $y$, $f(y,P)$ is integrable in the brick $I$ with integral $F(y)$, and if $\lim_{y \rightarrow \infty} F(y)$  exists, then $f(P) = \lim_{y \rightarrow \infty} f(y,P)$ exists almost everywhere in $I$. If we put $0$ for $f(P)$ where it is otherwise undefined then $f(P)$ is integrable in $I$ with integral $F$, i.e.,
\[
\int_I \lim_{j \rightarrow \infty} f(y,P) d\mu = \lim_{j \rightarrow \infty} \int_I f(y,P) d\mu.
\]
\item
If in 1 we omit the hypothesis that $\lim_{y \rightarrow \infty} F(y)$ exists but add that $f(P)$ exists as the limit of $f(y,P)$ without being integrable, then $F(y) \rightarrow \infty$ as $y \rightarrow \infty$.
\end{enumerate}
\end{theorem}
\proof
We can assume $f(y,P) \geq 0$ or else replace it by $f(y,P)-f(0,P)$. For each fixed $P$, $f(y,P)$ is monotone increasing in $y$, and so tends to a limit or to $+\infty$. Also,
\[
f(n,P) \leq f(y,P) \leq f(n+1,P),\;\;\;\;\;n \leq y \leq n+1,
\]
and similarly for their integrals, so we can assume $y$ is an integer.
For integers $j,N$, let $X_j$ be the set of $P$ with $f(j,P)>N$, and let $Y$ be the set with \textbf{either} $\lim_{j \rightarrow \infty} f(j,P) >N$ \textbf{or} $f(j,P) \rightarrow \infty$ as $j \rightarrow \infty$. Then $X_{j+1} \supseteq X_j$, $Y = \bigcup_{j=1}^\infty X_j$, so that by Theorems \ref{properties of variation}, \ref{16}, and the definition of variation,
\begin{eqnarray*}
F & \geq & \int_I f(j,P)d\mu \;\;\;\geq \;\;\; V(f(j,P) \mu(J);I) \vt
&=& V(f\mu;I;I)\;\;\;\geq \;\;\; V(f(j,P) \mu; I;X_j) \;\;\;\geq\;\;\;NV(\mu;I;X_j).
\end{eqnarray*}
Then $\mu^*(X_j \cap I)=V(\mu;I;X_j) \leq FN^{-1}$, and, letting $j \rightarrow \infty$,
\[
\mu^*(Y \cap I) \leq \frac F N.
\]
If $X$ is the set of points $P$ where $F9j,P) \rightarrow \infty$ as $j \rightarrow \infty$ we have $X \subset Y$ for all values of $N$, so
\[
\mu^*(X \cap I) \leq \mu^*(Y \cap I) \leq \frac F N,
\]
and $\mu^*(X \cap I) =0$. Thus, by Theorem \ref{13: measure integral null set}, $\int_If(j,P) \chi(X,P)d\mu =0$ and
\[
\int_I f(j,P)d\mu = \int_I f(j,P) \chi(\setminus X,P)d\mu
\]
by the additivity of the integral with respect to the integand. Theorem \ref{Levi} now gives 1 since $f(P) = \lim_{j\rightarrow \infty}f(j,P)\chi(\setminus x,P)$.
For 2, if the final result is false so that $\int_If(j,P) d\mu$ does not tend to $+\infty$, the sequence of integrals is bounded and part 1 shows that $f(P)$ is integrable.

\section{The Integrability of Functions of Functions}
This section enables us to deal with another Lebesgue convergence theorem. We use a real- or complex-valued function of two real numbers. We assume that
\begin{enumerate}
\item
$r$ is homogeneous for non-negative numbers, i.e., $r(ax_1,ax_2) = ar(x_1,x_2)$ for all $a \geq 0$;
and  $r$ satisfies a Lipschitz condition:
\item
 $|r(y_1,y_2) - r(x_1,x_2) | \leq A|y_1-x_1| +B|y_2-x_2|$ for constants $A<0$, $B>0$.
 \end{enumerate}
For example, 2 holds if the partial derivatives of $r$ exist and are bounded.

 \begin{theorem} \label{18}
 Let $S$ be the real line or the subset $x \geq 0$, let $Z$ be the complex plane, and let $f_j(P):I \mapsto S$ be integrable to $F_j(J)$ in each $J$ contained in the brick $I$ ($j=1,2$). If $r(x_1,x_2): S \times S \mapsto Z$ satisfies 1 and 2, then
 \[
 V\left(r(f_1(P),f_2(P))\mu(J)-r(F_1(J),F_2(J));I\right)=0.
 \]
 \end{theorem}
\proof
By Theorem \ref{Henstock Lemma} (\ref{HL2}), given $\ve>0$ there exists $\delta(P)>0$ with
\[
V(f_j(P)\mu(J) - F_j(J);I;\delta) \leq 4\ve,\;\;\;\;j=1,2.
\]
This remains true for $j=1,2$ if we replace $\delta_j(P)$ by $\delta(P)=\min\{\delta_1(P),\delta_2(P)\}>0$. For each division $D$ of $I$ compatible with $z\delta(P)$, 1 and 2 give
\begin{eqnarray*}
&& (\D)\sum |r(f_1(P), f_2(P))\mu(J) -r(F_1(J),F_2(J))| \vt
&=& (\D)\sum |r(f_1(P)\mu(J), f_2(P)\mu(J)) -r(F_1(J),F_2(J))| \vt
&\leq & A(\D)\sum |f_1(P) \mu(J) -F_1(J)| + B (\D)\sum |f_2(P) \mu(J) -F_2(J)| \vt
&\leq & 4(A+B) \ve.
\end{eqnarray*}
\nproof

\begin{theorem}\label{19}
With the conditions of Theorem \ref{18} let $r$ also be real-valued and satisfy
\begin{equation}\label{19a}
r(x_1+y_1, x_2+y_2) \leq r(x_1,x_2) + r(y_1,y_2)\;\;\mbox{ for all }\;\; x_1,x_2,y_1,y_2 \in S.
\end{equation}
Then $r(f_1(P),f_2(P))$ is integrable if and only if, for some $\delta(P)>0$ in $I$, and for all divisions $\D$ of $I$ which are compatible with $\delta(P)$,
\begin{equation}\label{19b}
(\D)\sum r(f_(P),f_2(P))\mu(J) \;\;\mbox{ is bounded above.}
\end{equation}
\end{theorem}
\proof
Given $\ve>0$, by Theorem \ref{1} there exists $\delta_1(P)>0$ such that for all divisions $\D$ of $I$ compatible with $\delta_1(P)$,
\begin{equation}\label{19c}
(\D)\sum |r(f_1,f_2)\mu - r(F_1,F_2)|< \ve.
\end{equation}
Thus the boundedness above of (\ref{19b}) is equivalent to the boundedness above, for divisions $\D$ compatible with some $\delta_2(P)>0$ of
\begin{equation}\label{19d}
(\D)\sum r(F_1,F_2).
\end{equation}
Let (\ref{19d}) be bounded above with supremum $s$. Then there is a division $\D_1$ of $I$ compatible with $\delta_2(P)$, and formed of bricks $I_1, \ldots ,I_m$, for which
\begin{equation}\label{19e}
s-\ve < (\D)\sum r(F_1,F_2)\leq s.
\end{equation}
Replacing the $\delta_j(P)$ of the geometric theorem (Theorem \ref{geometry}) by $\delta_2(P)$, we can construct a new $\delta(P)>0$ which we denote by $\delta_3(P)>0$. Let $\D$ be a division of $I$ compatible with $\delta_3(P)>0$. Then if a brick $J$ of $\D$ overlaps with more than one $I_k$, we can cut it up into two or more bricks $J_1, \ldots , J_r$ with the same associated point $P$ and each $J_j$ lying in one $I_k$. The sum for $\D$ can thus be changed into a sum for a division $\D_2$ of $I$, each brick of which lies in one $I_k$. By the additivity of $F_1, F_2$, and by induction on (\ref{19a}), if $D_{2k}$ is the division of $I_k$ from $\D_2$, we have
\begin{equation}\label{19f}
r(F_1(I_k) F_2(I_k)) = r\left((\D_{2k})\sum F_1, (\D_{2k})\sum F_2\right) \leq (\D_{2k})\sum r(F_1,F_2).
\end{equation}
By (\ref{19c}), (\ref{19e}), (\ref{19f}), and by the definitions of $s$ and $\D_1$,
\[
\begin{array}{rllllll}
s- \ve &<& (\D_2)\sum r(F_1,F_2) & \leq & s,&& \vt
s-2 \ve & < &(\D_2) \sum r(f_1,f_2) \mu &=& (\D) \sum r(f_1,f_2) \mu &<& s+ \ve.
\end{array}
\]
Hence, if (\ref{19b}) is true, $r(f_1,f_2)$ is integrable in $I$ to $s$---the supremum of sums (\ref{19d}). Conversely, if $r(f_1,f_2)$ is integrable in $I$ to the value $R$, the sum in (\ref{19b}) lies between $R-\ve$ and $R+\ve$ for suitable $\delta(P)>0$, and so is bounded. \nproof

\begin{theorem}\label{20}
If $f_1,f_2,f_3$ are real and integrable in $I$ with $f_1 \leq f_3$, $f_2 \leq f_3$, then $\max\left\{f_1(P),f_2(P)\right\}$ is integrable in $I$.
\end{theorem}
\proof
(\ref{19a}) of Theorem \ref{19} is satisfied with $r(x_1,x_2)= \max\{x_1,x_2\}$, since, for $j=1,2$,
$
x_j+y_j \leq \max\{x_1,x_2\} + \max\{y_1,y_2\}$,
so
\[
\max\{x_1+x_2,y_1+y_2\} + \max\{x_1,x_2\} + \max\{y_1,y_2\}.
\]
For Lipschitz condition (2) with $A=B=1$, we have, for $j=1,2$
\[
\begin{array}{rll}
x_j & = & (x_j-y_j) + j_j \;\;\; \leq \;\;\; |x_j - y_j| + y_j \vt
&\leq & |x_1-y_1| + |x_2-y_2| + \max\{y_1,y_2\},\vt
\max\{x_1,x_2\} &\leq & |x_1-y_1| + |x_2-y_2| + \max\{y_1,y_2\}.
\end{array}
\]
Interchanging $(x_1,x_2)$, $(y_1,y_2)$, we get
\[
\max\{y_1,y_2\} - \max\{x_1,x_2\} \leq |x_1-y_1| + |x_2-y_2|.
\]
Hence
\[
\left|\max\{y_1,y_2\} - \max\{x_1,x_2\} \right| \leq |x_1-y_1| + |x_2-y_2|.
\]
Finally, (\ref{19b}) of Theorem \ref{19} is bounded above since $\max\{f_1,f_2\} \leq f_3$, and, as $f_3$ is integrable to $F_3$,
\[
(\D) \sum \max \{f_1,f_2\} \mu \;\;\leq \;\;(\D) \sum f_3 \mu \;\; < \;\; F_3 + \ve
\]
for $\D$ compatible with some $\delta(P)>0$.  \nproof

\noindent
\textbf{Corollary:} If $f_1,f_2,f_3$ are integrable in $I$ with $f_1 \geq f_2$, $f_2 \geq f_3$, then $\min\{f_1,f_2\}$ is integrable in $I$. This follows from Theorem \ref{20} if $f_j$ is replaced by $-f_j$.

\section{Majorised or Dominated Convergence}
Theorem \ref{20} and its Corollary, with Lebesgue's monotone convergence theorem, lead to Lebesgue's majorised convergence theorem and Fatou's lemma.

\begin{theorem} \label{t21}
Let $g_1(P), g_2(P), f_j(P)$ be real and integrable in $I$ with
\begin{equation}\label{21a}
f_j(P) \geq g_1(P),\;\;\;\;j=1,2.
\end{equation}
Then
\begin{equation}\label{21b}
\int_I \liminf_{j \rightarrow \infty} f_j(P) d\mu \leq  \liminf_{j \rightarrow \infty}\int_I f_j(P) d\mu,
\end{equation}
the integrand on the left existing when the right hand side is finite. If, instead of (\ref{21a}), we have
\begin{equation}\label{21c}
f_j(P) \leq g_2(P),\;\;\;\;j=1,2,
\end{equation}
then
\begin{equation}\label{21d}
\int_I \limsup_{j \rightarrow \infty} f_j(P) d\mu \geq  \limsup_{j \rightarrow \infty}\int_I f_j(P) d\mu,
\end{equation}
the integrand on the left existing when the right hand side is finite (i.e., $\neq - \infty$).
If both (\ref{21a}) and (\ref{21c}) are true, with $f(P) = \lim_{j\rightarrow \infty} f_j(P)$ existing almost everywhere, then
\begin{equation}\label{21e}
\int_I  f(P) d\mu =  \lim_{j \rightarrow \infty}\int_I f_j(P) d\mu,
\end{equation}
\end{theorem}
\proof
By Theorem \ref{20} (Corollary), applied $m-k$ times,
\begin{equation}\label{21f}
\min_{k \leq j \leq m} f_j(P)
\end{equation}
is integrable in $I$.
As (\ref{21f}) is monotone increasing for $m \rightarrow \infty$, and is bounded below by $g_1(P)$ so that $\int_I f_1(P)d\mu \geq \int_I g_1(P)d\mu$, we can apply Lebesgue's monotone convergence theorem to show that the following exist:
\begin{eqnarray}
\int_I \inf_{{{j \geq k}}}f_j(P) d\mu & \leq & \lim_{m \rightarrow \infty} \int_I \min_{k \leq j \leq m}f_j(P) d\mu \nonumber \vt
 & \leq &  \lim_{m \rightarrow \infty} \min_{k \leq j \leq m}\int_I f_j(P) d\mu \nonumber \vt
 &=& \inf_{{{j \geq k}}} \int_I f_j(P)d\mu, \label{21g}
 \end{eqnarray}
using $\int_I  \min_{k \leq j \leq m} f_j(P) d\mu \leq \int_I f_j(P) d\mu$ for any $j$ in $k \leq j \leq m$.
The first integrand in (\ref{21g}) is monotone increasing. If the right hand side of (\ref{21b}) is finite, then by (\ref{21g}), the first integral in (\ref{21g}) is bounded as $k \rightarrow \infty$, and Lebesgue's monotone convergence theorem shows that
\[
\liminf_{j \rightarrow \infty} f_j(P),\;\;=\;\;\lim_{k \rightarrow \infty} \left( \inf_{{{j \geq k}}} f_j(P)\right),
\]
exists as a finite limit almost everywhere, and is integrable, with
\[
\int_I \liminf_{j \rightarrow \infty} f_j(P) d\mu =\lim_{k \rightarrow \infty} \int_I \inf_{{{j \geq k}}} f_j(P) d\mu
\leq \liminf_{j \rightarrow \infty} \int_I f_j(P) d\mu.
\]
This gives (\ref{21b}). For (\ref{21d}) we replace $f_j$ by $g_2 - f_j$ in (\ref{21b}). (\ref{21e}), which is Lebesgue's majorised convergence theorem, then follows.
Fatou's lemma is (\ref{21a}) and (\ref{21b}). \nproof

Instead of using an integer parameter $j$ we can use a continuous parameter $y$ that takes all values in $y\geq 0$.
But in the first part we cannot prove from the integrability of $f(P,y)$ in $I$ for all $y\geq 0$ that
\begin{equation}\label{21h}
\mbox{g.l.b.} \left\{f(P,y)\;:\; Y \leq y \leq Z\right\}
\end{equation}
is also integrable. For let us take the dimension of the brick to be $n=1$. Let $I=[0,1]$ and let $X_1$ be a set with non-integrable characteristic function on $[0,1]$. Then for $P \in I$ and $0 \leq y \leq 1$ let
\[
f(P,y) = \left\{
\begin{array}{ll}
1 & (P=y \in X_1), \vt
0 & \mbox{otherwise};
\end{array}
\right.
\;\;\;\;\;
f(P,y+m) =f(P,y)\;\;\;(m=1,2,3, \ldots).
\]
Then $f(P,y)$ is bounded for $P \in I$ and $y \geq 0$ and is integrable. But
\[
\limsup_{y \rightarrow \infty} f(P,y) = \sup_{0\leq y <1} f(P,y) = \chi(X_1,P) \;\;\;\; \;\;\;\;(0\leq P<1)
\]
and $-f(P,y)$ falsifies the conditions. However, if to the hypotheses in the first part of Theorem \ref{21} (replacing $j$ by $y$) we add the integrability of (\ref{21h}) $\inf_{Y \leq y\leq Z} f(P,y)$ whenever $0 \leq Y \leq Z$, then the corresponding result (\ref{21b}) follows by a similar proof.

Similarly we can obtain the analogue of (\ref{21d}) if, for all $0 \leq Y<Z$,
\begin{equation}\label{21j}
\mbox{l.u.b.}\left\{ f(P,y)\;:\; Y \leq y \leq Z \right\}
\end{equation}
is integrable, with analogues of the hypotheses.

However the analogue of (\ref{21e}) is true with analogous conditions without the integrability of (\ref{21h}) and (\ref{21j}).
For if $\{y_j\}$ is any sequence tending to $+\infty$, then almost everywhere we have
\[
f(P) = \lim_{y \rightarrow \infty} f(P,y) = \lim_{j \rightarrow \infty} f(P,y_j).
\]
We first choose $\{y_j\}$ so that
\[
\liminf_{y\rightarrow \infty} \int_I f(P,y) d\mu = \lim_{j\rightarrow\infty} \int_I f(P,y)d\mu,
\]
and then (\ref{21b}) shows that this is
\[
\geq \int_I \lim_{j \rightarrow\infty}f(P,y_j) d\mu = \int_I \lim_{y\rightarrow \infty} f(P,y) d\mu.
\]
Similarly for the analogue of (\ref{21d}) using another $\{y_j\}$ so that the analogue of (\ref{21e}) follows.

Two theorems follow from the continuous parameter version of Theorem \ref{21}.

\begin{theorem}\label{t22}
If for some $\ve>0$ and $0<|y-a|\leq\ve$, $f(P,y)$, $g_1(P)$, $g_2(P)$ are real and integrable in $I$ with
\[
g_1(P) \leq f(P,y) \leq g_2(P),\;\;\;\;\;\; \lim_{y\rightarrow a} f(P,y) =f(P,a)
\]
almost everywhere, then
\[
\lim_{y\rightarrow a} \int_I f(P,y) d\mu =  \int_I f(P,a)d\mu,
\]
and the latter integral exists.
\end{theorem}

\begin{theorem}\label{t23}
If $\frac{\partial f(P,y)}{\partial y}$ exists almost everywhere in $I$ and if for some $\ve>0$ and all $y$ in $0<|y-a|\leq\ve$, $f(P,y)$, $g_1(P)$, $g_2(P)$ are real and integrable in $I$ with
\[
g_1(P) \leq \frac{f(P,y)-f(P,a)}{y-a} \leq g_2(P),\;\;\;\;\;\; 0<|y-a|\leq\ve,
\]
 then
\[
\frac d{dy} \int_I f(P,y) d\mu = \int_I \frac{\partial f(P,y)}{\partial y} d\mu,
\]
and the latter integral exists.
\end{theorem}

\begin{theorem}\label{t24}
Let $f_j(P)$ and $|f_j(P)|$ be integrable in $I$ for $j=1,2,3, \ldots$. If $\sum_{j=1}^\infty \int_I |f_j(P)|d\mu$ is a convergent series then $f(P) \equiv \sum_{j=1}^\infty f_j(P)$ is an absolutely convergent series almost everywhere and is integrable, and
\begin{equation}\label{24a}
\int_I \sum_{j=1}^\infty f_j(P)d\mu = \int_I fd\mu =
\sum_{j=1}^\infty \int_I f_j(P)d\mu.
\end{equation}
\end{theorem}
\proof
For $k=1,2,3,\ldots$,
\[
\int_I \sum_{j=1}^k |f_j(P)|d\mu  = \sum_{j=1}^k \int_I |f_j(P)|d\mu \leq \int_I \sum_{j=1}^\infty |f_j(P)|d\mu,
\]
the series being convergent. By Lebesgue's monotone convergence theorem, $f_0(P) \equiv \sum_{j=1}^\infty |f_j(P)|$ is convergent almost everywhere and is integrable in $I$, and $f$ exists almost everywhere. In Theorem \ref{21} we now replace $g_1$, $g_2$, $f_k$ by the present $-f_0$, $f_0$, $\sum_{j=1}^k \Re(f_j)$, and obtain the integrability of $\Re(f)$. Similarly for the integrability of $\Im(f)$, and so of $f$; giving (\ref{24a}).
\nproof

For a real point function $f$ and constants $a<b$, let
\[
X(f=a),\;\;\;X(f>a),\;\;\;X(f\geq a),\;\;\;X(f<b),\;\;\;X(f\leq b),\;\;\;X(a \leq f \leq b),\ldots ,
\]
denote the set of points for which the appropriate equality or inequality holds.

\begin{theorem} \label{t25}
Let $f(P)$ and $|f(P)|$ be real and integrable in $I$, and for real constant $b$ let $c(b;P)$ be the characteristic function $ \chi(X(f\geq b),P)$ of $X(f \geq b)$.Then $c(b;P)$ is integrable in $I$ with integral equal to $\mu^*(X(f\geq b)\cap I)$. Similar results hold for the other sets, and, by definition, these sets are measurable.
\end{theorem}
\proof
Since
\[
-|f| -|b| \leq f \leq |f| + |b|,\;\;\;\;-|f| -|b| \leq b \leq |f| + |b|,
\]
Theorem \ref{20} and its corollary show that $\max\{f,b\}$ and $\min\{f,b\}$ are integrable in $I$ for every constant $b$. Thus
\[
c(a,b;P), \equiv \max\{\min\{f,b\},a\},\;\;\;\;\;\;a<b,
\]
is integrable. It is $f$ when $a\leq f\leq b$, it is $a$ when $f \leq a$, and it is $b$ when $f \geq b$. As
\[
c(b;P) = \lim_{a \rightarrow b-} \frac{c(a,b;P) - a}{b-a},
\]
which is the limit of a monotone decreasing sequence bounded by $0$ and $1$, Levi's monotone convergence theorem shows that $c(b;P)$ is integrable. \nproof

\begin{theorem} \label{t26}
If $X,Y$ are two measurable sets, then so are $X \cap Y$ and $X \cup Y$.
\end{theorem}
\proof
As $\chi(X,P)$ and $\chi(Y,P)$ are integrable, so is their sum. Therefore, by Theorem \ref{25}, the set where $\chi(X,P) + \chi(Y,P) \geq 2$ is measurable, and this is the set $X \cap Y$. Also
\[
\chi(X \cup Y,P) = \chi(X,P) + \chi(Y,P) - \chi(X\cap Y,P)
\]
and so is integrable. \nproof

\noindent
\textbf{Corollary:}
If $f, |f|, g, |g|$ are real and integrable, the set where $a \leq f\leq b$ and where $r \leq g\leq s$ is measurable for constants $a<b$, $r<s$.

For real functions $f$ we say that $f$ is \emph{measurable} when $X(f\geq b)$ is measurable for each real constant $b$. Thus Theorem \ref{25} says that $f$ is measurable if $f$ and $|f|$ are integrable. A partial converse follows:
\begin{theorem} \label{t27}
If $g \geq 0$ is integrable in $I$, and if $f$ is real and measurable in $I$ with $|f|\leq g$, then $f$ and $|f|$ are integrable in $I$.
\end{theorem}
\proof
Replacing $f$ by $f+g$ if necessary, we can assume that $f \geq 0$. For each integer $m \geq 0$ we put $a_j = j2^{-m}$, $j=0,1,2, \ldots$, and
\[
f_m(P) = a_j\;\;\mbox{ for }\;\;a_j \leq f(P) < a_{j+1},\;\;\;\;f_m(P) =0\;\;\mbox{ for }\;\;f(P) \geq m2^{-m}.
\]
Then $f_m(P) \leq f(P) \leq g(P)$, and $f_m(P)$ is monotone increasing in $m$ for each fixed $P$. Also
\begin{equation} \label{*}
\int_I f_(P)d\mu = \sum_{j=1}^{m2^m} \int_I a_j \chi\left(X(a_j \leq f <a_{j+1}),P\right) d\mu
\end{equation}
which exists since $X(a \leq f<b)$ is measurable for each pair of constants $a<b$. As $f_m(P) \rightarrow f(P)$, Levi's monotone convergence theorem gives the integrability of $f$ when $f \geq 0$, otherwise the integrability of $f+g$, and so $(f+g)-g$, both of these being integrable. Similarly for $|f| \geq 0$. \nproof

\noindent
\textbf{Corollary:} Instead of $|f| \leq g$ we can assume that $f \geq 0$ and that the integral of $f_n$ is bounded as $n \rightarrow \infty$.

If $f$ is real and bounded, the proof of Theorem \ref{27} contains a more usual and more difficult definition of the integral.The right hand side of (\ref{*}) is equal to
\[
\sum_{j=1}^{m2^m} a_j \mu^*\left(X(a_j \leq f<a_{j+1}) \cap I\right).
\]
First $\mu^*(X \cap I)$ is defined in a more difficult way. Then the range of $f \geq 0$ is divided up by a finite number of $0=b_0 < b_1 < \cdots < b_N$, and then when the sets involved are ``measurable'' it is shown that sums
\[
\sum_{j=1}^N b_{j+1} \mu^* \left( X(b_{j+1} \leq f < b_j)\right)
\]
tend to a limit as $b_N \rightarrow \infty$ and $\max\{ b_j-b_{j-1}: 1 \leq j \leq N\} \rightarrow 0$. In the other theory, the limit is taken as the definition of the integral
of $f$.

By a proof similar to that of Theorem \ref{27}, the limit is equal to our $\int_I fd\mu$ provided Lebesgue's $\mu^*$ has the same values as ours.

\section{H\"{o}lder's and Minkowski's Inequalities for Integrals}
For a constant $p>1$ let us consider $\phi(x) = x^{p-1}$. If $y = \phi(x)$ then $x = \psi(y) = y^{-(p-1)}$. Putting $q = 1 +(p-1)^{-1}$ we can write $\psi(y) = y^{q-1}$. Here
\[
qp-q = (p-1)+1,\;\;\;\;p+q = pq,\;\;\;\; \frac 1p+\frac 1q =1.
\]
The functions $y=\phi(x)$ and $x=\psi(y)$ have the same graph, and the rectangle formed by joining the points $(0,0)$, $(0,y)$, $(x,y)$, $(x,0)$ is divided into two sections, $A$ and $B$ by this graph. A simple geometric argument\footnote{The inequalities of  H\"{o}lder and Minkowski are derived in J.E.~Littlewood's \emph{Lectures on the Theory of Functions}, Oxford University Press, 1944.}.....????

Thus, when $x>0$, $y>0$, we have
\begin{eqnarray}
\frac{x^p}p+\frac{y^q}q \geq xy, \label{mink1} \vt
\frac{x^p}p+\frac{y^q}q = xy\;\;\;\mbox{ if and only if } y=x^{p-1}\;\;\;(\mbox{i.e. }\; y^q = x^p. \label{mink2}
\end{eqnarray}
Clearly this holds also when $x \geq 0$, $y \geq 0$.

Now let $x_1, \ldots , x_m, y_1, \ldots , y_m, a,b$ all be non-negative Then by (\ref{mink1}),
\begin{equation} \label{mink3}
ab \sum_{j=1}^m x_j y_j \leq \frac{a^p}p \sum_{j=1}^m x_j^p + \frac{b^q}q \sum_{j=1}^m y_j^q.
\end{equation}
By (\ref{mink2}), there is equality in (\ref{mink3}) if and only if
\begin{equation} \label{mink4}
b^q y_j^q = a^px_j^p,\;\;\;\;1 \leq j \leq m.
\end{equation}
the right hand side of (\ref{mink3}) becomes $p^{-1}+q^{-1}=1$ on taking the special values
\[
a^{-p} = \sum_{j=1}^m x_j^p,\;\;\;\;b^{-q} = \sum_{j=1}^m y_j^q,
\]
and then Minkowski's inequality yields H\"{o}lder's inequality:
\begin{equation}\label{holder1}
\sum_{j=1}^m x_jy_j \leq \frac 1{ab} = \left(\sum_{j=1}^m x_j^p\right)^\frac 1p\left(\sum_{j=1}^m y_j^q\right)^\frac 1q
\end{equation}
with equality when
\begin{equation}\label{holder2}
\frac{y_j^q}{x_j^p} \;\;\mbox{ is constant in }\;\;j.
\end{equation}

\begin{theorem}\label{t28}
Let $f \geq 0$, $g \geq 0$ be point functions, let $p>1$ be constant, and let $q$ be defined by $p^{-1}+q^{-1}=1$. If $f^p$ and $g^p$ are integrable in $I$ then $fg$ is integrable, and
\begin{equation} \label{28a}
\int_I fgd\mu \leq \left(\int_I f^pd\mu\right)^\frac 1p \left(\int_I g^q d\mu\right)^\frac 1q.
\end{equation}
This is H\"{o}lder's inequality for integrals. If equality occurs in (\ref{28a}) then
\begin{equation} \label{28b}
\frac{f^p(P)}{g^q(P)}\;\;\;\mbox{ is constant almost everywhere.}
\end{equation}
\end{theorem}
\proof
For each $b>a \geq 0$, $s>r\geq 0$, the set
\[
X\left(a \leq f <b, r \leq g<s \right)= X\left(a^p \leq f^p <b^p, r^q \leq g^q<s^q \right)
\]
is measurable by Theorem \ref{26} (Corollary). For each integer $m$ we put
\begin{eqnarray*}
a_j &=& \frac j{2^m},\;\;\;\;j=0,1,2, \ldots, \vt
f_m(P)&=& \left\{
\begin{array}{lll}
a_j & \mbox{if} & a_j \leq f(P) < a_{j+1},\;\;\;j=0, 1,2, \ldots, m2^m,\vt
0 & \mbox{if} & f(P) \geq m+2^{-m};
\end{array} \right.
\vt
g_m(P)&=& \left\{
\begin{array}{lll}
a_j & \mbox{if} & a_j \leq g(P) < a_{j+1},\;\;\;j=0, 1,2, \ldots, m2^m,\vt
0 & \mbox{if} & g(P) \geq m+2^{-m};
\end{array} \right.
\end{eqnarray*}
Then $f_m(P)g_m(P)$ is integrable as $X\left(a \leq f <b, r \leq g<s \right)$ has an integrable characteristic function; also $f_m(P)g_m(P)$ is monotone increasing in $m$ for each fixed $P$ and $\rightarrow f(P)g(P)$ as $m \rightarrow \infty$. By (\ref{mink1})
\[
abf_m(P)g_m(P) \leq \frac{a^pf_m^p(P)}{p} + \frac{b^qg_m^q(P)}{q} \leq \frac{a^p}p f^p(P) + \frac{b^q}q g^q(P),
\]
and we are given that the final expression is integrable in $I$. Hence by Levi's monotone convergence theorem, $f(P)g(P)$ is integrable over $I$, and
\begin{equation}\label{28c}
ab \int_I f(P)g(P)d\mu \leq \frac{a^p}p\int_I f^p(P)d\mu + \frac{b^q}q \int_Ig^q(P)d\mu,
\end{equation}
which is the analogue of (\ref{mink3}).
As $f \geq 0$, if $\int_I f^pd\mu =0$ so that $V(f^p\mu;I)=0$, then $f=0$ almost everywhere. Thus we can assume that $\int_I f^pd\mu >0$, and similarly $\int_I g^qd\mu >0$.
We have $\int_I h(P) d\mu \geq 0$, where
\begin{equation}\label{28d}
h(P) = \frac{a^p}p f^p(P) + \frac{b^q}q g^q(P) - abf(P)g(P) \geq 0.
\end{equation}
Hence $\int_I h(P)d\mu = V(h\mu;I)$, and there is equality in (\ref{28c}) if and only if this is $0$. So, by (\ref{28d}) and (\ref{mink2}),
\begin{equation}\label{28e}
a^p f^p(P) = b^q g^q(P)\;\mbox{ almost everywhere.}
\end{equation}
This result holds for all values of $a>0, b>0$. Taking the special values given by
\[
a^{-p} = \int_I f^p(P)d\mu,\;\;\;\;\;\;b^{-q} \int_I g^q(P)d\mu,
\]
(\ref{28c}) becomes (\ref{28a}), and equality in  (\ref{28a}) occurs if and only if  (\ref{28e}) is true for the special values of $a$ and $b$, giving  (\ref{28b}).
\nproof

\begin{theorem}\label{t29}
Let $f \geq 0$, $g \geq 0$ be point functions and let $p>1$ be constant. If $f^p(P)$ and $g^q(P)$ are integrable in $I$ then $f(f+g)^{p-1}$ and $g(f+g)^{p-1}$ are integrable in $I$ and
\begin{equation}\label{29a}
\left(\int_I(f+g)^pd\mu\right){\frac 1p}\leq  \left(\int_If^pd\mu\right){\frac 1p}  \left(\int_Ig^pd\mu\right){\frac 1p},
\end{equation}
(Minkowski's inequality), with equality in (\ref{29a}) when
\begin{equation}\label{29b}
\frac fg \;\;\;\mbox{ is constant almost everywhere.}
\end{equation}
\end{theorem}
\proof
If $p^{-1}+g^{-1}=1$ then $(p-1)q=pq-q=p$. If $\int_I (f+g)^pd\mu = 0$ then, as $f\geq 0$, $g\geq 0$, we have $f=0$ and $g=0$ almost everywhere. Hence we can assume $\int_I(f+g)^pd\mu >0$.
Assuming the integrability of
\[
f(f+g)^{p-1},\;\;\;\;\;\;g(f+g)^{p-1},
\]
we have by H\"{o}lder's inequality
\begin{eqnarray*}
\int_I \left(f+g\right)^pd\mu &=& \int_I f\left(f+g\right)^{p-1}d\mu + \int_I g\left(f+g\right)^{p-1}d\mu  \vt
&\leq &  \left(\int_I f^p d\mu \right)^{\frac 1p} \left(\int_I (f+g)^p d\mu \right)^{\frac 1q}
+ \left(\int_I g^p d\mu \right)^{\frac 1p} \left(\int_I (f+g)^p d\mu \right)^{\frac 1q}.
\end{eqnarray*}
Dividing by $\left(\int_I (f+g)^p d\mu \right)^{\frac 1q}$ we have (\ref{29a}).
Equality occurs when $f^p(f+g)^p$ is constant, $g^p(f+g)^p$ is constant, so $fg^{-1}$ is constant almost everywhere. Hence we need only prove the integrability of $f(f+g)^{p-1}$ and $g(f+g)^{p-1}$ by the methods of Theorems \ref{27} and \ref{28} using the $f_m(P)$ and $g_m(P)$ of Theorem \ref{28}. By the inequality (\ref{mink1}),
\begin{eqnarray*}
f_m(f_m+g_m)^{p-1} & \leq & \frac{f_m^p}p + \frac{(f_m+g_m)^q}q \vt
&\leq & f^p + 2^p \max\{f^p,g^p\}\;\;\;\leq \;\;\;f^p + 2^p (f^p + g^p) \;\;\;<\;\;\;2^p.
\end{eqnarray*}
By Levi's monotone convergence theorem $f(f+g)^{p-1}$ is integrable, and similarly for $g(f+g)^{p-1}$. \nproof

\section{Mean Convergence}
Denote by $L^p$ the set of all measurable functions $f(P)$ with $|f|^p$ integrable in $I$. For norm we use
\[
||f|| = \left(\int_I|f(P)|^p d\mu \right)^\frac 1p.
\]
If $b$ is a constant then $||bf|| = |b|\, ||f||$.
If $f,g$ are measurable then so is $f+g$, as
\begin{eqnarray*}
|f+g| &\leq & |f| + |g|,\vt ||f+g|| &\leq & ||f|| + ||g||, \vt
\left(\int_I |f+g|^p d\mu \right)^\frac 1p & \leq & \left(\int_I (|f|+|g|)^p d\mu \right)^\frac 1p\vt &\leq & \left(\int_I |f|^p d\mu \right)^\frac 1p + \left(\int_I |g|^p d\mu \right)^\frac 1p .
\end{eqnarray*}
$||f-g||$ can be regarded as the ``distance'' between $f$ and $g$. There is equality in the triangle inequality if $fg^{-1}$ is constant almost everywhere. Using this distance we can say that $\{f_j(P)\}$ is a \emph{fundamental} or \emph{Cauchy} sequence in $L^p$ if $||f_j(P) - f_k(P)|| < \ve$ for all integers $j>k>K$ where $K$ depends on $\ve>0$. $\{f_j(P)\}$ is a \emph{convergent sequence in} $L^p$ if there is a function $f(P)$ in $L^p$ such that \[||f_j(P) - f(P)|| < \ve\] for all $j > K_1$ where $K_1$ depends on $\ve$.
We prove that $L^p$ is \emph{complete}, that is, every fundamental sequence is a convergent sequence.

\begin{theorem}\label{t30}
Let $p>0$ be fixed, and let $\{f_j(P)\}$ be a sequence of measurable point functions such that, to each $\ve>0$, there corresponds a $K=H(\ve)$ such that, for $j>k \geq K(\ve)$,
\begin{equation}\label{30a}
\int_I |f_j(P) - f_k(P)|^pd\mu < \ve,
\end{equation}
where we suppose that the integral exists for all $j,k$. Then there is a point function $f$, and also, given $\ve>0$, there is an integer $K_1 = K_1(\ve)$ such that, for $j \geq K_1$,
\begin{equation}\label{30b}
\int_I |f_j(P) - f(P)|^pd\mu < \ve,
\end{equation}
the integral existing for all $j$. Further, if there is another point function $g(P)$ satisfying the same conditions then $f=g$ almost everywhere. If each $f_j(P)$ is in $L^p$ (taking $p>1$) then $f(P) \in L^p$.
\end{theorem}
\proof
(If (\ref{30a}) is true it is said that $f_j(P)$ converges \emph{in mean} with index $p$. If (\ref{30b}) is true, $f_j(P)$ converges in mean to $f(P)$ with index $p$.)
To each integer $m$ there corresponds an integer $j_m$ such that, for $j>k\geq j_m$,
\[
\int_I\left|f_j(P) - f_k(P)\right|^p d\mu < \frac 1{2^{mp+m}}.
\]
If this works for a particular $j_m$ it will be true for $j_m+1$ or $j_m+2$ or $\ldots $ or $j_m+r$ replacing $j_m$. Hence we can suppose $j_m < j_{m+1}$ ($m=1,2,3, \ldots$).
(The measurability of $X(a \leq f_j \leq b)$  and $X(a \leq f_j \leq b, c\leq f_k\leq d)$ implies the measurability of
$X(a \leq |f_j-f_k| \leq b)$.) If $X_m$ is the set of $P$ where $|f_{j_m}(P) - f_{j_{m+1}}(P)| \geq 2^{-m}$, then, using Theorem \ref{25} for measurability,
\begin{eqnarray*}
\mu^*(X_m\cap I)2^{-mp} & = & \int_I 2^{-mp} \chi(X_m,P) d\mu \vt
& \leq & \int_I |f_{j_m}(P) - f_{j_{m+1}}(P)|^p d\mu \;\;\;<\;\;\;2^{-mp-m}.
\end{eqnarray*}
Therefore $\mu^*(X_m\cap I) < 2^{-m}$, and $Y_M = \bigcup _{m\geq M} X_m$ implies
\[
\mu^*(Y_M\cap I) \leq \sum_{m=M}^\infty \mu^*(X_m\cap I) < \sum_{m=M}^\infty 2^{-m} = 2^{1-M}.
\]
If $P$ is not in $Y_M$ then every one of the terms $|f_{j_m} (P) - f_{j_{m+1}} (P) | < 2^{-m}$ so that $\sum_{m=M}^\infty (f_{j_m}(P) - f_{j_{m+1}}(P))$ is absolutely convergent and so convergent.
As the partial sums
\[
\sum_{m=M}^N (f_{j_m}(P) - f_{j_{m+1}}(P)) = f_{j_M}(P) - f_{j_N + 1}(P)
\]
we can let $N \rightarrow \infty$ so that $f(P) = \lim_{N\rightarrow \infty} f_{j_N}(P)$ exists for all $P$ not in $Y_M$. Thus the set where the limit does not exist lies in $Y_M$ for every $M$ and, for all $M$,
\[
\mu^*(X \cap I) \leq \mu^*(Y_M \cap I) < 2^{1-M},
\]
so $\mu^*(X \cap I) =0$. We can now define $f(P)$ as $0$ for $P \in X$. By Fatou's lemma, with $g_1(P) =0$, (\ref{30a}) implies
\begin{eqnarray*}
\int_I |f_j(P) - f(P)|^p d \mu &=& \int_I |f_j(P) - \lim_{N \rightarrow \infty}f_{j_N}(P)|^p d \mu \vt
 &\leq & \liminf_{N \rightarrow \infty} \int_I |f_j(P) - f_{j_N}(P)|^p d \mu \;\;\;<\;\;\;\ve
 \end{eqnarray*}
 for $j \geq K(\ve)$. Hence we have (\ref{30b}). Further, if there is a $g(P)$ with the same properties then, for each $\ve>0$,
 \[
 \int_I |f-g|^p d\mu \leq \liminf_{N \rightarrow \infty} \int_I |f_{j_n} - g|^p d\mu \leq \ve,
 \]
 so $\int_I |f-g|^pd\mu =0$, $|f-g|^p = 0$ almost everywhere, and $f=g$ almost everywhere. Finally, if $p \geq 1$ and if $f_j(P) \in L^p$ for each $j$, then
 \[
 \left(\int_I |f(P)|^pd\mu\right)^\frac 1p \leq  \left(\int_I |f(P) - f_j(P)|^pd\mu\right)^\frac 1p  +   \left(\int_I |f_j(P)|^pd\mu\right)^\frac 1p \leq \ve + ||f_j||
 \]
 for $j \geq (\ve)$. Thus $||f||$ is finite and $f \in L^p$. \nproof

\begin{theorem} \label{t31}
Let $p>1$ be fixed with $p^{-1}+q^{-1}=1$. Let $g \geq 0$ and $f_j(P)$ be point functions such that $f_jg$, $g^q$, $|f_j-f_k|^p$ ($j,k=1,2,3, \ldots$) are all integrable in $I$ with (\ref{30a}) of Theorem \ref{30} satisfied. Then
\[
\lim_{j \rightarrow \infty} \int_I f_jg d\mu = \int_I fg d\mu.
\]
\end{theorem}
\proof
By H\"{o}lder's inequality and Theorem \ref{27} the first integral below exists and
\begin{equation} \label{31a}
\int_I |f_j-f_k| g \,d\mu \leq \left(\int_I |f_j-f_k| ^p g d\mu \right)^\frac 1p\left(\int_I g^q d\mu \right)^ \frac 1q.
\end{equation}
Therefore by Fatou's lemma we have
\[
\int_I |f_j-f| g \,d\mu \leq \liminf_{r \rightarrow \infty}\left(\int_I |f_j-f_{j_r}| ^p g d\mu \right)^\frac 1p\left(\int_I g^q d\mu \right)^ \frac 1q ,
\;\;\;\rightarrow 0
\]
as $j \rightarrow \infty$.
Thus to finish the proof we need only show that $fg$ is integrable, noting that $fg = \lim_{r \rightarrow \infty} f_{j_r}g$.
This follows since $f_jg$ is given integrable for each $j$, and since
\[
\sum_{r=1}^\infty \int_I \left| f_{j_r} - f_{j_{r+1}}\right| g d\mu \leq \left(\sum_{r \geq 1} 2^{-r -\frac rp}\right)\left( \int_I g^q d\mu \right)^\frac 1q < \left( \int_I g^q d\mu \right)^\frac 1q.
\] 

\section{The Cauchy Extension on the real Line}
\begin{theorem}\label{t32}
Let $\int_{[a,b]} f d\mu = F(a,b)$ exist for all $b$ in $a<b<c$, and let
\[
\lim_{b \rightarrow c-} F(a,b) \equiv F
\]
exist. Then $\int_{[a,c]}fd\mu$ exists and is equal to $F$.
\end{theorem}
\proof
Let us put $b_j = c - (c-a)2^{-j}$, $j=0,1,2, \ldots $. Then $F(b_{j-1}, b_j)$ exists for $j=1,2, \ldots$. Given $\ve>0$ there are functions $\delta_j(P)>0$ in $[b_{j-1}, b_j]$ such that, for all divisions $\D_j$ of $[b_{j-1},b_j]$ compatible with $\delta_j(P)$,
\[
\left|(\D_j)\sum f(P)\mu(J) - F(b_{j-1},b_j)\right| \leq \frac{\ve }{2^{j+1}}.
\]
Clearly we can suppose that
\[
\delta_j(P) \leq \min\left\{ \frac 12\left|P-b_j\right|,\frac 12 \left|P-b_{j-1}\right|\right\}
\]
for $b_{j-1}<P<b_j$. For points $P$ in $a \leq P <c$ we put
\[
\begin{array}{rll}
\delta(P)& =& \left\{
\begin{array}{lll}
\delta_j(P)& \mbox{when} & b_{j-1}<P<b_j, \vt
\min\left\{\delta_j(P), \delta_{j+1}(P), \frac{c-a}{2^{j+1}}\right\} & \mbox{when} &  P=b_j,\;\;j \geq 1,
\end{array}
\right.
\vt
\delta(a) &=& \min \left\{ \delta_1(a), \frac{c-a}4\right\}.
\end{array}
\]
Then any interval compatible with $\delta(P)$ that includes a $b_k$ has that $b_k$ as associated point and no other $b_k$ lies in the interval. There is a $\delta(c)>0$ such that if $c-\delta(c) \leq b<  c$ then
\[
|F(a,b) -F| < \frac \ve 4,\;\;\;\;\;|f(c) (c-b)| < \frac \ve 4.
\]
 If $\D$ is a division of $[a,c]$ compatible with $\delta(P)$, there is an interval $[v,c]$ in $\D$ with associated point $c$ since, for any other associated point the interval will not stretch as far as $c$; and each of the other intervals lies entirely within an interval $[b_{j-1},b_j]$ or else it has a $b_j$ as associated point. Thus we can split up the sum for $\D$ into a finite number $m$ of sums for divisions $\D_j$ of $[b_{j-1},b_j]$, $1 \leq j \leq m$, a sum for a partial division $\mathcal{Q}$ of $[b_m, b_{m+1}]$ that is a division of $b_m, v]$ and the term $f(c)(c-v)$. hence
 \begin{eqnarray*}
 |(\D)\sum f(P)\mu(J) -F| & \leq & \left| \sum_{j=1}^m \left((\D_j)\sum f(P)\mu(J) - F(b_{j-1}, b_j) \right) \right| \;+ \vt
 && \;\;\;\; +\; \left|(\mathcal{Q}) \sum f(P)\mu(J) - F(b_m,v)\right|\;+ \vt
 && \;\;\;\; +\; |f(c)(c-v)| + |F(a,b) -F| \vt
 &\leq & \sum_{j=1}^m \frac \ve{2^{j+1}} + \frac{2\ve}{2^{m+2}} + \frac \ve 4+ \frac \ve 4 \;\;\; =\;\;\;\ve,
 \end{eqnarray*}
 giving the result. \nproof

Similar results hold if $f$ is integrable in $[b,c]$ for all $b$ in $a<b<c$, and $\lim_{b \rightarrow a+}F(b,c) =F$; and also if $f$ is integrable in $[u,v]$ to $F(u,v)$ for all $u,v$ in $a<u<v<b$, and if
\[
\lim_{u \rightarrow a+} F(u, \frac{a+c}2)=F_1,\;\;\;\;\;\lim_{v \rightarrow c-} F(\frac{a+c}2,v)=F_2
\]
exist. This result does not necessarily hold if $u \rightarrow a+$,  $v \rightarrow c-$ simultaneously. For instance if $u-a =c-v$. Thus $\int_{-1}^1 \frac{dx}x$ does not exist even though $\int_{-1}^{-\ve} + \int_\ve^1 =0$ for all $\ve>0$ so that the $\lim_{\ve \rightarrow 0}$ exists.

When $a$ is finite it is usual to define $\int_a^\infty f d\mu$ as $\lim_{b \rightarrow \infty} \int_a^b fd\mu$. Similarly
\[
\int_{-\infty}^b fd\mu =\lim_{a \rightarrow - \infty}\int_a^b fd\mu,\;\;\;\;\;\;\mbox{ and }\int_{-\infty}^\infty fd\mu =     \lim_{a \rightarrow - \infty}\lim_{b \rightarrow  \infty}\int_a^b fd\mu
\]
where $a$ and $b$ are independent.

But we need not prove the integration theorems again for these intervals over an infinite range as we shall see. The ordinary proofs will suffice. However the theorems need not hold for the integral $\int_{-A}^A f d\mu$ which turns up in contour integration; for instance $f(x)=1$ for $x \geq 0$ and $f(x)=-1$ for $x<0$.

We can define integrals over an infinite range, without using an extra limit, by means of suitable divisions. To find out how to do this we map $(-\infty,\infty)$ onto $(-1,1)$ by using $x(1-x^2)^{-1}$. This mapping is strictly increasing, with derivative always $>0$, taking values from $-\infty$ to $\infty$ in $-1<x<1$. The intervals $(P-\delta(P), P+\delta(P))$ of $(-1,1)$ become intervals in $(-\infty,\infty)$ though $P$ does not transform onto the centre of the new interval. What are missing are the transforms of the associated points $-1,+1$ that often have to be used to obtain a division of $[-1,1]$. The corresponding intervals $[-1, \alpha]$, $[\beta,1]$ transform into infinite intervals $(-\infty, a]$, $[b,+\infty)$ for real values of $a,b$. This suggests the following scheme.

Let $\delta(P)>0$ be defined for each real number $P$, and let us use finite intervals compatible with $\delta(P)$ in the usual way. Also let $a<b$ be two real numbers. Then a division $\D$ of $(-\infty, \infty)$ compatible with $a,b, \delta(P)$ is defined to be $(-\infty,u]$, $[v, +\infty)$ and a division of $[u,v]$ compatible with $\delta(P)$, such that $u<a<b<v$. By convention, in $(\D)\sum f(P) \mu(J)$ we put $f(P)\mu(J)=0$ when $J=(-\infty,u]$ or $[v, +\infty)$.

The theory of the integral using a finite interval $[u,v]$ has an analogous theory using an infinite interval $(-\infty, \infty)$. In particular, the analogue of Theorem \ref{32} shows that
\[
\int_{-\infty}^\infty fd\mu = \lim_{a \rightarrow -\infty}\lim{b \rightarrow \infty} \int_a^b fd\mu
\]
which is the usual definition. But the usefulness of the definition of a division of $(-\infty, \infty)$ lies in the fact that earlier theorems are fairly general, i.e.~they do not depend tightly on the geometry, but are true whatever definition we give of a division, within broad limits.

The new scheme for $(-\infty,\infty)$ is within these limits so that we have, for example:
\emph{
If $g(P)$, $h(P)$ are integrable in $(-\infty, \infty)$, and if $\{f_j(P)\}$ is a sequence of point functions each integrable in $(-\infty,\infty)$ and satisfying}
\[
g(P) \leq f_j(P) \leq h(P)
\]
\emph{for all $j$ and a.e.~in $P$, then}
\begin{eqnarray*}
\int_{-\infty}^\infty \liminf_{j \rightarrow \infty} f_j(P)d\mu
&\leq &\liminf_{j \rightarrow \infty} \int_{-\infty}^\infty f_j(P)d\mu \vt
&\leq &\limsup_{j \rightarrow \infty} \int_{-\infty}^\infty f_j(P)d\mu \vt
&\leq &\int_{-\infty}^\infty \limsup_{j \rightarrow \infty} f_j(P)d\mu.
\end{eqnarray*}
\emph{If also $\lim_{j\rightarrow \infty} f_j(P) =f(P)$ exists a.e.~then $\int_{-\infty}^\infty f(P)d\mu$ exists and is equal to $\lim_{j\rightarrow\infty}\int_{-\infty}^\infty f_j(P) d\mu$, which also exists.
}

\section{Integration by Parts}
Again we restrict the integration to be on the real line. Let $g(x)$ be a point function, then define
\[
\mu(g;a,b) = g(b)-g(a) = \Delta g = [g]_a^b.
\]
Then for the integral of $f$ with respect to $g$ over the interval $I$, $\int_Ifdg$, we use sums $(\D) \sum f(P) \mu(g;J)$ over divisions $\D$ of $I$. In other words if $I=[u,v]$ then we have $f(P)(g(v)-g(u))$.

\begin{theorem} \label{t33}
If $\int_a^b fdg$ and $\int_a^b gdf$ both exist with sum
\begin{equation}\label{33a}
\int_a^b fdg + \int_a^b gdf = f(b)g(b) - f(a)g(a)
\end{equation}
for each $[a,b] \subseteq I$, then
\begin{equation}\label{33b}
V(\Delta f \Delta g;I)=0
\end{equation}
and $\int_a^b fdg =  f(b)g(b) - f(a)g(a)-\int_a^b gdf$. Conversely, if (\ref{33b}) is true, and if the first integral exists then so does the second, and their sum is (\ref{33a}).
\end{theorem}
\proof
First, note that we can take the associated points at the ends of the intervals. For if $P$ is the associated point of $[u,v]$ and if $u<P<v$ then
\[
P-\delta(P) <U<P<v<P+\delta(P)
\]
so that $P$ can also be the associated point of $[u,P]$ and $[P,v]$. We can repeat this for every associated point that lies inside an interval, and then every associated point will lie at the end of an interval. We can then use the identity
\[
\begin{array}{r}
f(u)\left(g(u) - g(v)\right)+
g(u)\left(f(u) - f(v)\right)-
f(u)g(u) - f(v)g(v) \vt
=\;\;\;
\left(f(u)-f(v)\right)\left(g(u) - g(v)\right).
\end{array}
\]
In the first part of the theorem,
\[
V\left(f\Delta g - \int f\,dg;I\right)=0,\;\;\;\;\;V\left(g\Delta f - \int g\,df;I\right)=0,
\]
by Theorem \ref{Henstock Lemma}.
Hence
\begin{eqnarray*}
V(\Delta f \Delta g;I) &=& V(f\Delta g + g\Delta f - \Delta (fg);I)\vt
&=& V\left(\int f d g + \int g d f - \Delta (fg);I\right) \;\;=\;\;V(0;I)\;\;=\;\;0.
\end{eqnarray*}
Conversely, if (\ref{33b}) holds and if the first integral exists then $V(f\Delta g - \int fdg;I)=0$. Hence
\begin{eqnarray*}
0&=& V(\Delta f\Delta g;I) \;\;=\;\; V(f \Delta g + g \delta f - \Delta (fg) ;I) \vt
&=& V \left( \int fdg + g\Delta f - \Delta (fg) ;I\right) \;\;=\;\;V(g \Delta f - H;I)
\end{eqnarray*}
where $H(a,b) = [fg]_a^b - \int_a^b fdg$. Here $H$ is finitely additive. Therefore $\int_a^b gdf$ exists and is equal to $\Delta (fg) - \int_a^b fdg$ for all $[a,b] \subseteq I$, giving the converse. \nproof

If $f$ is continuous and $g$ of bounded variation (i.e.~$V(\Delta g;I)$ is finite) then (\ref{33b}) is true. For, by uniformity of continuity, $|\Delta f|<\ve$ for $\delta(x)$ small enough. Thus, for divisions $\D$ of $I$ compatible with $\delta(x)$ we have
\[
(\D) \sum |\Delta f\Delta g| \leq (\D) \sum \ve |\Delta g|.
\]
In turn, by choice of a suitable $\delta(P)>0$, this last is $\leq \ve\left(V(\Delta g;I)+ \ve\right)$. This holds for all $\ve>0$, so (\ref{33b}) is true. It can be proved that $\int_I fdg$ exists. Hence $\int_I gdf$ exists and we have the formula for integration by parts.

All conditions in Theorem \ref{33} are independent. Let
\[
\begin{array}{rrrrrr}
F(x) &=&\left\{
\begin{array}{lll}
1& \mbox{for}& x \leq 0, \vt
0 & \mbox{for}& x>0,
\end{array}\right.
& \;\;\;\;\;G(x) &=& 1-F(x).
\end{array}
\]
We can arrange that $0<\delta(x) <|x|$ when $x \neq 0$ and then every $\delta$-compatible division $\D$ of $[-1,1]$ contains an interval $[u,v]$ either with $u<0<v$ or two intervals $[u,0]$ and $[0,v]$, and each inter val has associated point $0$. As $F(0)=1$,
\[
(\D) \sum F(x) \Delta g = F(0) \left(G(v) - G(u)\right) = 1 \times \left(F(u) - F(v)\right) =1,
\]
and $\int_{-1}^1 F dG =1$. Similarly
\[
(\D) \sum G(x) \Delta F = G(0) \left(F(v) - F(u) \right) =0,\;\;\;\;\;\;\;\;\;\;\int_{-1}^1 GdF=0.
\]
Also $[F(x)G(x)]_{-1}^1 =0$. Thus the two integrals exist, but there is no integration by parts.
In fact, by a similar proof,
\begin{eqnarray*}
(\D) \sum |\Delta F \Delta G| &=& |F(v)-F(0)|\,|Gv) - G(0)| + |F(0) - F(u)|\,|G(0)-G(u)| \vt
&=& \left( F(v) - F(0)\right)^2 + \left(F(0)-F(u)\right)^2 \;\;=\;\;1+0,
\end{eqnarray*}
so $V(\Delta F \Delta G;[-1,1])=1$.

On the other hand, let
\begin{eqnarray*}
F(x)& =& \frac 1{x\ln x}\;\;\;\;(0<x<1),\;\;\;\;\;\;F(0)=0,\;\;\;\;\;\;G(x)=x,\vt
h(\ve)& =&  \sup\left\{|F'(x)|\;:\;\ve \leq x \leq \frac 12\right\}\;\;\;\mbox{ for each }\;\; \ve>0.
\end{eqnarray*}
Let
\[
\delta(0) = \exp\left({-\frac 2 \ve}\right) >0,\;\;\;\;\;\;\;\;\;\;\;\;\delta(x) = \min \left\{ \frac x2, \frac \ve{h\left(\frac 12 x\right)}\right\} \;\;\;\;(0 \leq x \leq \frac 12).
\]
Let $\D$ be a $\delta$-compatible division of $[0, \frac 12]$. Then $\D$ contains an interval $[0,u]$ for some $u>0$. For $x \geq 0$, $\delta(x) \leq \frac 12 x$. Hence the associated point must be $0$. Using the mean value theorem,
\[
(\D) \sum | \Delta F \Delta G| = (\D) \sum^* |F'(x) \mu(J)^2 + u |F(u)|
\]
where $\sum^*$ means ``omit $[0,u]$'', and where $\xi$ is some point in $J$. Let $J$ have associated point $\eta$. Then
\[
|F'(\xi)| \mu(J)^2 \leq \frac{h(\xi) 2 \ve}{h(\frac 12 \eta)} \mu(J).
\]
If $J = [v,w]$ we have $|\eta - v| \leq \delta(\eta) \leq \frac 12 \eta$, and $0 \leq \eta - v$, so mod sign unnecessary; then we have $v \geq \frac 12 \eta$ and so $h(v) \leq h(\frac 12 \eta)$. Hence
\begin{eqnarray*}
F'(\xi) \mu(J)^2 &\leq &\frac{2h(v) \ve}{h\left( \frac 12 \eta\right)} \mu(J) \;\;\leq \;\;2 \ve \mu(J), \;\;\;\;\mbox{ and}\vt
(\D) \sum \Delta F \Delta G &\leq & 2\ve\left(\frac 12 -u \right) + \frac 1{|\ln u |}.
\end{eqnarray*}
As $\ve \rightarrow 0+$ we have $u \rightarrow 0+$ and so $V(\Delta F\Delta G;[0, \frac 12]) =0$.
But, using $y = \ln x$,
\begin{eqnarray*}
\int_\ve^{\frac 12} F dG &=& \int_\ve^{\frac 12} \frac{dx}{x\ln x} \;\;=\;\;\int_{\ln \ve}^{-\ln 2} \frac{e^y dy}{e^y y} \vt
&=& -\left[ \ln y\right]_{\ln \frac 1\ve}^{\ln 2} \;\;=\;\; \ln \ln \frac 1\ve - \ln \ln 2 \;\;\rightarrow\;\;+\infty
\end{eqnarray*}
as $\ve \rightarrow \infty$. If $\int_0^{\frac 12} F dG$ exists, $\int_\ve^{\frac 12} F dG$ should tend to it as $\ve$ tends to $0$.
Hence $\int_0^{\frac 12} F dG$ cannot exist. Further, as the integration by parts works in $[\ve,\frac 12]$ for each $0<\ve<\frac 12$, we have
\begin{eqnarray*}
\int_\ve^{\frac 12} G dF &=&
\frac{2 \times \frac 12}{-\ln 2} - \frac 1\ve \frac \ve{\ln \ve} - \ln\ln \frac 1\ve + \ln\ln 2 \vt
&=& \frac 1{\ln \frac 1\ve} - \frac 1{\ln 2} -\ln\ln \frac 1\ve + \ln\ln 2 \;\;\rightarrow \;\;-\infty
\end{eqnarray*}
as $\ve \rightarrow 0$. Hence $\int_0^{\frac 12} GdF$ cannot exist.
Thus if $V(\Delta F\Delta G;I) =0$ it does not follow that $\int_IFdG$ exists.

\section{Fubini's Theorem}
Denote by $E^j$ the Euclidean space of dimension $j$ and let $m$ and $n$ be positive integers with sum $N$, $N=m+n$. Then the co-ordinates $(x_1, \ldots ,x_N)$ of points in $E^N$ can be separated into two collections giving points $x=(x_1, \dots , x_m)\in E^m$ and $y=x_{m+1}, \ldots , x-N) \in E^n$ so that we can write $(x,y)$ for the points of $E^N$. Note that we could take any $m$ co-ordinates out of the $N$ co-ordinates provided we take the same co-ordinates each time. Then we could re-arrange the co-ordinates so that the chosen $m$ come first. So the choice of the first $m$ co-ordinates is quite general.

Let $X \subseteq E^m$, $Y \subseteq E^n$ be sets in their respective spaces, and let $Z=X\times Y$ denote the set in $E^N$ formed of all $(x,y)$ with $x \in X$, $y \in Y$. We call $X \times Y$ the Cartesian product of $X$ and $Y$. Thus $E^N = E^m \times E^n$.

\begin{theorem} \label{t34}
Let $I,J$ be bricks in $E^m,E^n$ respectively, and let $\delta(P,Q)>0$ be defined at all points of $K=I \times J$, a brick in $E^N$. Then to each $x$ of $I$ and each division $\D(x)$ of $J$, with a given set of associated points, which is compatible with $\frac 12 \delta(x,y)$, one division for each point $x$, there is a $\delta_1(x)>0$ on $I$ such that if $I^* \subseteq E^m$, with associated point $x$, is compatible with $\delta_1(x)$, and if $J^* \in \D(x)$ with associated point $y$, then $I^* \times J^*$, with associated point $(x,y)$, is compatible with $\delta(x,y)$.
\end{theorem}
\proof
For the fixed $x$, $\delta(x,y)>0$ is a function of the $y\in J$. Let $\D(x)$, a division of $J$ compatible with $\frac 12 \delta(x,y)$, consist of $J_1, \ldots , J_r$ with associated points $y_1, \ldots , y_r$ respectively. Define
\[
\delta_1(x) = \min \left\{ \frac 12 \delta(x,y_j): 1 \leq j \leq r\right\} >0.
\]
If $I^*$, with associated point $x$, is compatible with $\delta_1(x)>0$, and if $J_j \in \D(x)$ with associated point $y_j$, then $I^* \times J_j$, with associated point $(x,y_j)$, is compatible with
\[
\sqrt{\delta_1(x)^2 + \frac 14 \delta(x,y_j)^2} \leq \sqrt{\frac 12 \delta (x,y_j)^2} < \delta(x,y_j)
\]
as required. \nproof

\begin{theorem} \label{t35} \textbf{(Fubini)}
Let $\mu(U), \nu(V)$ be the volume functions in $E^m, E^n$, and let $\lambda(U \times V) =\mu(U) \nu(V)$ be the volume function in $E^N$. Let $f(x,y)$ be a point function that is integrable in $K = I \times J$ to the value $F$, where $I,J$ are bricks in $E^m, E^n$ respectively. then the point functions
\begin{equation} \label{35a}
g(x)=\int_I f(x,y)d\nu,\;\;\;\;\;\;\;\;h(y) = \int_I f(x,y)d\mu
\end{equation}
exist almost everywhere in $I$, $J$ respectively. Putting $g(x)=0$, $h(y)=0$ where the integrals do not exist, then
\begin{equation} \label{35b}
\int_I g(x) d\mu = \int_J h(y) d\nu = F = \int_K f(x,y) d\lambda.
\end{equation}
\end{theorem}
\proof
Let $F(U \times V)$ be the integral of $f$ in $U\times V$. Given $\ve>0$ there is a function $\delta(x,y)>0$ defined in $K$ such that every $\delta$-compatible division $\D$ of $K$ satisfies
\begin{equation}\label{35c}
(\D) \sum |f(x,y) \lambda(U \times V) - F(U \times V)| < \ve.
\end{equation}
First we prove that the $g(x)$ exist for almost all $x \in I$. Let $\delta_2(y)>0$ be defined in $J$, and for each fixed $x$ let $S(x, \delta_2)$ be the set of sums $(\D_1)\sum f(x,y) \nu(V)$ for all $\delta_2$-compatible divisions $\D_1$ of $J$. Let $X_p$ be the set of all $x$ in $I$ for which $S(x, \delta_2)$ has diameter $>\ p^{-1}$ for all $\delta_2(y)$. Then there are two sums $S_1(x)$ and $S_2(x)$ in $S(x, \delta_2)$ such that
\begin{equation}\label{35d}
\left| S_1(x) - S_2(x) \right| > \frac 1p.
\end{equation}
For each $x \in I$ and the given $\delta(x,y)$ we choose a division $\D_2(x)$ of $J$ that is compatible with $\frac 12 \delta(x,y)$ with sum
\[
S_3(x) = (\D_2(x))\sum f(x,y)\nu(V).
\]
When $x \in X_p$ we can arrange that $S_3(x)$ is the $S_1(x)$ of (\ref{35d}) when $\delta_2(y) = \frac 12 \delta(x,y)$. Also we can choose another division $\D_3(x)$ of $J$ with sum $S_4(x)$, the $S_2(x)$ of (\ref{35d}) when $\delta_2(y) = \frac 12 \delta(x,y)$.
From $\D_2(x)$ when $x \in J$ we can define a $\delta_1(x)$ as in Theorem \ref{34}. Also from $\D_2(x)$ when $x \in I \setminus X_p$, and from $\D_3(x)$ when $x \in X_p$, we can define a similar $\delta_3(x)>0$ Then
\[
\delta_4(x) = \min\{\delta_3(x), \delta_4(x)\}>0
\]
has the property of Theorem \ref{34} relative to both $\D_2(x)$ and $\D_3(x)$ when $x \in X_p$. Let $\D_4$ be a division of $I$ compatible with $\delta_4(x)$. Then if $U \in \D_4$ with associated point $x \in X_p$, and if $V \in \D_2(x)$ or $\D_3(x)$, then $U \times V$ is compatible with $\delta(x,y)$. Thus, by (\ref{35c}) and (\ref{35d}),we have
\[
\begin{array}{r}
(\D_4 \times \D_2(x))\sum |f(x,y) \mu(U)\nu(V) - F(U \times V)| < \ve, \vt (\D_4 \times \D_3(x))\sum |f(x,y) \mu(U)\nu(V) - F(U \times V)| < \ve, \vt
(\D_4) \sum \chi(X_p,x) |S_1(x) \mu(U) - F(U \times V)| < \ve, \vt  (\D_4) \sum \chi(X_p,x) |S_2(x) \mu(U) - F(U \times V)| < \ve, \vt
(\D_4)\sum \chi(X_p,x) \mu(U) \;\; \leq \;\; p(\D_4) \sum \chi(X_p,x) |S_1(x) - S_2(x) | \mu (V)\;\;< \;\;2 \ve p,
\end{array}
\]
giving $V (\mu, I \cap X_p; \delta_4) \leq 2 \ve p$. As $\ve>0$ is arbitrary,
\[\mu^* (X_p \cap I) =0,\;\;\;\;\;\;\mu^* \left(\bigcup_{p=1}^\infty \left(X_p \cap I\right)\right) =0.\]
Outside $X=\bigcup_{p=1}^\infty X_p$, (\ref{35d}) is false for each $p>0$ and some $\delta_2(y)>0$ depending on $p$. The situation is as in Theorem \ref{integrability on sub-interval}.
Thus $g(x)$ exists for each $x \notin X$. Further, by Theorem \ref{13: measure integral null set}, $V(S_1(P) \mu(U);I;X) =0$, so that there is a function $\delta_5(P)>0$ defined in $I$ for which
\[
(\D_5)\sum \chi(X,P) \left|S_1(P)\right| \mu(U) \;\;<\;\;\ve
\]
for each $\delta_5$-compatible division $\D_5$ of $I$. If $\D_5$ happens also to be $\delta_4$-compatible, where we take $S_1(P)$ within $\ve$ of $g(P)$, $P \notin X$, then (\ref{35c}) gives
\begin{eqnarray}\label{35e}
(\D_5)\sum \chi(X,P) \left|S_1(P) \mu(U) - F(U \times V)\right| &< &\ve,\nonumber \vt
(\D_5)\sum \chi(X,P) \left| F(U \times V)\right| &<& 2\ve.
\end{eqnarray}
Then, from (\ref{35c}) again,
\begin{eqnarray*}
(\D_5)\sum \chi(\setminus X,P) | S_1(P) \mu(U)-F(U\times V)| &<& \ve, \vt
(\D_5)\sum \chi(\setminus X,P) | g(P) \mu(U)-F(U\times V)| &< &\ve + \ve \mu(I) \ve,
\end{eqnarray*}
so $\left|(\D_5)\sum \chi(\setminus X,P) g(P) \mu(U)-F\right|\;\;<$
\begin{eqnarray*}
 &<& \ve+ \ve \mu(I) + \left|(\D_5)\sum \chi(\setminus X,P)F(U\times V) -F\right| \vt
  &=& \ve+ \ve \mu(I) + \left|(\D_5)\sum \chi(\setminus X,P)F(U\times V) \right| \vt
  &<& 3\ve + 3\mu(I).
\end{eqnarray*}
As we have defined $g(P)=0$ in $X$ we have $g(P) \chi(\setminus I,P) = g(P)$, and so $g(P)$ is integrable in $I$ with integral $F$. Similarly $h(Q)$ exists almost everywhere as an integral, and is integrable to the value $F$. \nproof

\section{Inner Variation and Differentiation}
In $E^n$ we avoid ``needle-like'' bricks.
Thus, if $I$ is a brick and $J$ denotes cubes contained in $I$, the \emph{coefficient of regularity} of $I$ is
\[
r(I) \equiv \sup\left\{ \frac {\mu(I)}{\mu(J)}\;:\; J \subseteq I\right\}.
\]
If the longest edge of $I$ has length $d$ then $r(I) = \mu(I) d^{-n}$.
We say a set $\mathcal{C}$ of bricks $I$ is an $\alpha$-\emph{inner cover} of a set $X \subseteq E^n$, where $\alpha>0$, if to each point $P \in X$ there corresponds at least one brick $I \in \mathcal{C}$ with $r(I) \geq \alpha$ and with $P$ as associated point.

Let us take a fixed brick $I$ and function $\delta(P)>0$ of the points $P\in I$. If $\mathcal{C}$ is an $\alpha$-inner cover of a set $X \subseteq I$ such that each brick in $\mathcal{C}$ lies in $I$ and is compatible with $\delta(P)$, we say that $\mathcal{C}$ is compatible with $\delta(P)$ and $I$, or $(\delta,I)$\emph{-compatible}. Let
\[
IV_1\left(\mu;\mathcal{C}\right) := \sup\left\{\sum \mu(J) \right\},
\]
the $\sup$ for all finite, non-overlapping collections of $J \in \mathcal{C}$.  Let
\[
IV_2\left(\mu;I;X;\alpha;\delta\right) := \inf\left\{IV_1\left(\mu;\mathcal{C}\right) \right\}
\]
the $\inf$ for all $(\delta;I)$-compatible inner covers $\mathcal{C}$ of $X$. If $0 < \delta_1(P) \leq \delta_2(P)$ for $P \in I$, then a $(\delta_1,I)$-compatible $\mathcal{C}$ is also  $(\delta_2,I)$-compatible. Hence
\[
IV_2(\mu;I;X;\alpha;\delta_1) \geq IV_2(\mu;I;X;\alpha;\delta_2).
\]
Thus we can define the $\alpha$-\emph{inner variation} of $X$ to be
\[
IV(\mu;I;X;\alpha) = \sup\left\{ IV_2(\mu;I;X;\alpha;\delta)\;:\; \delta(P>0),\;P \in I\right\}.
\]
Clearly we have
\begin{eqnarray}
0\;\;\leq \;\; IV_1(\mu;\mathcal{C}) \leq \mu(I),\;\;\;\;\;\;\;\;\;\;\;\; 0\;\;\leq \;\;IV(\mu;I;X;\alpha)\;\;\leq \;\;\mu(I), \label{IV1}&&\vt
IV_1(\mu;\mathcal{C})\;\;\leq \;\; V(\mu\chi(X;P);I;\delta)\;\;=\;\;V(\mu;I;X;\delta),&&\nonumber \vt
IV(\mu;I;X;\alpha)\;\; \leq \;\; V(\mu;I;X). \label{IV2}&&
\end{eqnarray}
Now consider \textbf{differentiation}. Let $F(J)$ be a function of bricks $J \subseteq I$, let $P \in I$ and let a number $D$ exist with the following property. Given $\ve>0$, $\alpha>0$, there is a number $\delta(P, \ve, \alpha) >0$ such that for all bricks $J \subseteq I$ with associated point $P$, having $r(J) \geq \alpha$ and being compatible with $\delta(P, \ve,\alpha)$, we have
\begin{equation}\label{diff3}
\left| \frac{F(J)}{\mu(J)} -D \right| < \ve.
\end{equation}
Then $D$ is called the \emph{derivative of $F$ with respect to $\mu$, $I$ at $P$} and we denote $D$ by $D(F, \mu,I,P)$. If $P$ is on the border of $I$ the choice of bricks $J$ is restricted. But if $P$ is in the interior of $I$, then $D$ does not depend on $I$.

\begin{theorem} \label{t36}
If $F(J)$ is the integral of $f(P)$ in $J$ for all $J \subseteq I$, then
\begin{equation}\label{36a}
D(F;\mu;I;P) = f(P)
\end{equation}
except possibly at points $P$ of a set $X=\bigcup\{X_\alpha\,:\, \alpha>0\}$ with
\begin{equation}\label{36b}
IV(\mu;I;X_\alpha;\alpha) =0 \;\;\;\mbox{ all }\;\alpha >0.
\end{equation}
\end{theorem}
\proof
If the derivative does not exist at $P$, or if it exists but does not equal $f(P)$, then for some fixed $\ve_1(P)>0$, some fixed $\alpha>0$, and each $\delta(P)>0$, there is a $\delta(P)$-compatible brick $J \subseteq I$ with associated point $P$ that has $r(J) \geq \alpha(P)$, with
\begin{equation}\label{36c}
|F(J) - f(P)\mu(J)| \geq \ve_1(P)\mu(J).
\end{equation}
Let $X_\alpha$ be the set of all $P$ with $\alpha(P) \geq \alpha$, and let $X_{\alpha k}$ be the set of $P \in X_\alpha$ such that
\[
\ve_1(P) \geq \frac 1{2^k},\;\;\;\;k=m+1, m+2, \ldots ,
\]
with $X_{\alpha m}$ the set with $\ve_1(P) \geq 2^{-m}$. Take a suitable $\delta(P)>0$. By (\ref{HL2}) of Theorem \ref{Henstock Lemma} there is a $\delta_k(P)>0$ with
\[
(\D)\sum |F(J)-f(P)\mu(J)| \geq \ve_1(P)\mu(J)
\]
for all $\delta_k$-compatible divisions $\D$ of $I$. For  $P \in X_{\alpha k}$ let
\[
\delta^*(P) = \min \left\{ \delta(P), \delta_k(P) \right\}.
\]
Then for an $\alpha$-inner cover $\mathcal{C}$ of $X_\alpha$ that is compatible with $\delta^*(P)$ and $I$, and that comes from (\ref{36c}), and for all non-overlapping $J in \mathcal{C}$,
\begin{eqnarray*}
\sum \mu(J) &\leq & \sum_{k=m}^\infty 2^k \sum_{P \in X_{\alpha k}}\ve_1(P) \mu(J) \vt
& \leq & \sum_{k=m}^\infty 2^k \sum_{P \in X_{\alpha k}}\left| F(J) - f(P) \mu(J) \right| \vt
& \leq & \sum_{k=m}^\infty 2^k 4^{-k} \;\;= \;\;2^{1-m},
\end{eqnarray*}
so $IV_1(\mu; \mathcal{C}) \leq 2^{1-m}$ and $IV_2(\mu;I;X_\alpha; \alpha; \delta) \leq 2^{1-m}$. As $\delta \rightarrow 0$ we can take $m \rightarrow \infty$,
 so $IV(\mu;I;X;\alpha) =0$. \nproof

We now have a theorem due to Lebesgue in its original form with $n=1$.

\begin{theorem} \label{t37}
Let the real $f$ and $|f|$ be integrable in  brick $J$. Then
\[
\lim_{\mu(J) \rightarrow 0} \frac 1{\mu(J)} \int_J \left| f(P) - f(Q) \right| d\mu =0
\]
for each fixed $Q \in I$ and all $J \subseteq I$ with $r(J) \geq \alpha$ and associated point $Q$ except for a set $\bigcup \left\{X_\alpha\;:\; \alpha>0\right\}$ of points $Q$ with $IV(\mu;I;X_\alpha;\alpha)=0$.
\end{theorem}
\proof
The integrability of $f(P)$ and $|f(P)|$ imply the integrability of $f(P) - \beta$ and $|f(P) - \beta|$ for all constants $\beta$, since
\[
|f(P) - \beta| = \max\{f(P)-\beta, \beta-f(P)\} \leq |f(P)| + |\beta|.
\]
Thus by Theorem \ref{36} we have
\[
\lim_{\mu(J) \rightarrow 0} \frac 1{\mu(J)} \int_J \left| f(P) - \beta \right| d\mu =|f(Q) - \beta|
\]
for all $Q \in I$ and $r(J) > \alpha$, $J \subseteq I$, where the associated point of $J$ is $Q$, except for a set $X_\alpha(\beta)$ of $Q$ with $IV(\mu;I;X_\alpha(\beta))=0$.
We now take $\beta = \pm m 2^{-k}$ for positive integers $m,k$ and we can put such $\beta$ in the form of a sequence $\{\beta_j\}$. Let
\[
Y_1 = X_\alpha(\beta_1),\;\;\;\;\; Y_m = X_\alpha (\beta_m) \setminus \bigcup_{k<m} X_\alpha (\beta_k).
\]
Then
\[
0 \leq IV(\mu;I;Y_m;\alpha) \leq IV(\mu;I;X_\alpha(\beta_m);\alpha) =0.
\]
Thus, for each $\delta_m(P)>0$ and some $\alpha$-inner cover $\mathcal{C}_m$ of $Y_m$ compatible with $\delta_m(P)$, $IV(\mu;\mathcal{C}_m) \leq \ve 2^{-m}$.
Let
\[
X_\alpha = \bigcup_{m=1}^\infty X_\alpha (\beta_m)= \bigcup_{m=1}^\infty Y_m,
\]
and $\delta (P) = \delta_m(P)$ when $P \in Y_m$ ($m=1,2, 3, \ldots $).
Then
\[
\mathcal{C} = \bigcup_{m=1}^\infty \mathcal{C}_m
\]
is an $\alpha$-inner cover of $X$ compatible with $\delta(P)$, and
\[
IV_1(\mu; \mathcal{C}) \leq \sum_{j=1}^\infty IV_1(\mu; \mathcal{C}_m) \leq \ve.
\]
As $\delta(P)$ is as small as we please by choice of the $\delta_m(P)$, we have
\[
IV_1(\mu; I;X_\alpha; \alpha; \delta) \leq \ve,\;\;\;\;\;
IV_1(\mu; I;X_\alpha; \alpha) \leq \ve,\;\;\;\;\;
IV_1(\mu; I;X_\alpha) = 0.
\]
We now prove that for $Q$ in $\setminus \bigcup X_\alpha$, and for all real $\beta$,
\begin{equation} \label{37beta}
\lim_{\mu(J) \rightarrow 0} \frac 1{\mu(J)} \int_J |f(P) - \beta| d\mu = |f(Q)-\beta|.
\end{equation}
Let $\beta_m$ be one of the special $\beta$; then (\ref{37beta}) is true for $\beta_m$. Also
\begin{eqnarray*}
\int_J|f(P) - \beta| d\mu &\leq & \int_J|f(P) -\beta_m|d\mu + |\beta_m - \beta| \mu(J), \vt
\int_J|f(P) - \beta| d\mu  & \geq & \int_J|f(P) -\beta_m|d\mu - |\beta_m - \beta| \mu(J), \vt
|f(Q) - \beta| - 2|\beta_m - \beta| & \leq & |f(Q) - \beta_m| - |\beta_m - \beta| \vt
&\leq & \liminf_{\mu(J) \rightarrow 0} \frac 1{\mu(J)} \int_J |f(P) - \beta| d\mu \vt
&\leq & \limsup_{\mu(J) \rightarrow 0} \frac 1{\mu(J)} \int_J |f(P) - \beta| d\mu \vt
&\leq & |f(Q) - \beta_m| + |\beta_m - \beta| \vt
&\leq &|f(Q) - \beta| + 2 |\beta_m-\beta|.
\end{eqnarray*}
By construction we can take $\beta_m$ as near as we please to $\beta$ for suitable values of $m$. Hence in $\setminus \bigcup \{X_\alpha: \alpha>0\}$, and for all $\beta$, (\ref{37beta}) is true. Taking the special value $\beta=f(Q)$ we have the theorem. \nproof

To find a relation between $IV$ and $V$, we first need a continuity result.

\begin{theorem}\label{t38}
Let $J \subseteq I$ be bricks in $E^n$. Then $V(\mu;I;J) = \mu(J)$, and similarly for a finite union of bricks in $I$. Also, if $X \subseteq E^n$ then
\[
V(\mu;I;J\cap X)=V(\mu;J;X),
\]
i.e.~the boundary of $J$ has outer measure zero.
\end{theorem}
\proof
Let $B$ be the boundary of $J$. If $P \in I \setminus B$ we can define $\delta(P) <$ half the distance from $P$ to $B$. Then a $\delta$-compatible division $\D$ of $I$ is such that every brick that intersects $B$ has associated point on $B$. If also we have $\delta(P) \leq \ve$ on $B$ then the sum of the volumes of bricks in $\D$ with associated points on $B$ must be $\leq 2 \ve$(area $A$ of surface $B$). Thus
\begin{eqnarray*}
\mu(J) - 2\ve A &\leq & (\D) \sum \chi(J,P)\mu (J) \;\;\leq \;\;\mu(J) + 2 \ve A,\vt
\mu(J) - 2 \ve A &\leq & V(\mu;I;J;\delta) \;\;\leq \;\; \mu(J) + 2 \ve A,
\end{eqnarray*}
giving the first result. Simlarly for the second.
\nproof

\begin{theorem} \label{t39} \textbf{(Vitali)}
For fixed $\alpha>0$ and $X \subseteq I \subseteq E^n$ and for arbitrarily small $\delta(P)>0$, defined in $I$, let there be an $\alpha$-inner cover of $X$ that is compatible with $\delta(P),I$. Let $\mathcal{F}$ be the collection of all bricks in such $\alpha$-inner covers. Then there is a finite or countable sequence $\{I_p\}$ of mutually disjoint bricks of $\mathcal{F}$ such that
\begin{equation}\label{39a}
V\left(\mu;I;X \setminus \bigcup_{p=1}^\infty I_p\right)=0.
\end{equation}
\end{theorem}
\proof
$\{I_p\}$ is defined by induction. Let $I_1$ be an arbitrary brick of $\mathcal{F}$. When $I_1, \ldots ,I_p$ have been defined, no two with common points, (*) if their union includes all of $X$ we need go no further. Otherwise there is a point $P$ of $X$ that does not lie in the finite union of bricks---supposed closed, so their union is also closed. Hence there is a sphere $S$ with centre $P$ and radius $\delta(P)>0$ that has no point in common with the union. Now there is an $\alpha$-inner cover compatible with $\delta(P),I$ whatever the definiton of $\delta$ at other points of $I$, and so there is a brick of $\mathcal{F}$ with associated point $P$ that is cmpatible with $\delta(P)$. Thus there are bricks $J$ of $\mathcal{F}$ that are disjoint from $\bigcup_{j=1}^\infty I_j$. Let
\[
d_p = \sup\{\mbox{diameters of all such $J$}\}
\]
and let $I_{p+1}$ be any one of these $J$ with diameter $> \frac 12 d_p$. And so on. If (*) never occurs for any $p$ we put
$
Y=X \setminus \bigcup_{p=1}^\infty I_p.
$
If
\begin{equation}\label{39c}
V(\mu;I;Y)>0
\end{equation}
then, as the coefficient of regularity of each $I_p$ is $\geq \alpha$, we can associate with each $I_p$ a cube $K_p$ such that
\begin{equation}\label{39d}
I_p \subseteq K_p,\;\;\;\;\;\;\;\mu(I_p) \geq \alpha \mu(K_p).
\end{equation}
 Let $L_p$ be the cube with the same centre as $K_p$ but with diameter $4n+1$ times as big. ($n$ is the dimension of $E^n$.) The series in (\ref{39e}) must be convergent since
 \begin{equation}\label{39e}
 \sum_{p=N+1}^\infty \mu(L_p) < V(\mu;I;Y).
 \end{equation}
 Hence, by (\ref{39c}) there is an integer $N$ for which
 \[
 \sum_{p=N+1}^\infty \mu(L_p) < V(\mu;I;Y),\;\;\;\;\;\;\;\;\;\;\;\;\;\;V\left(\mu;I;\bigcup_{p=N+1}^\infty L_p\right)< V(\mu;I;Y)
 \]
by Theorem \ref{38} and $V(\mu;I; \bigcup X_j) \leq \sum_{j=1}^\infty V(\mu;I;X_j)$ (Theorem \ref{properties of variation} (3)).
Hence there is a point $P \in Y$ not belonging to any $L_p$ ($p>N$). By construction of $Y$, $P \notin \bigcup_{p=1}^\infty I_p$ and the $I_p$ are closed. Hence
\begin{equation}\label{39f}
\mbox{there is a sphere $S_1$ with centre $P$ and radius $\delta_1(P)>0$}
\end{equation}
that is disjoint from the closed set $\bigcup_{p=1}^N I_p$ and a brick $J \in \mathcal{F}$ with associated point $P$ and compatible with $\delta_1(P),I$.
But diam$(J) >0$. Now series (\ref{39e}) is convergent, hence $\mu(K_p) \rightarrow 0$ and diam$(K_p)\rightarrow 0$ as each $K_p$ is a cube. So we cannot have
\[
\mbox{diam}(J) \leq d_p \leq \mbox{diam}(I_{p+1}) \leq 2\mbox{diam}(K_{p+1})
\]
for all $p$. Since diam$J \leq d_p$ fails, then, for some $p$, $d_p < \mbox{diam}(J)$. By the definition of $d_p$, $J$ must have points in common with $\bigcup_{j=N+1}^p I_j$. Let $p_0$ be the smallest $p$ for which $J$ and $I_p$ have points in common. Then $J$ is disjoint from $I_p$ for $p=1,2, \ldots, p_0+1$  and
\begin{equation}\label{39g}
\mbox{diam}(J) \leq d_{p_0-1},
\end{equation}
by definition of function $d_p$. By (\ref{39f}), $p_0>N$. Therefore by definition of $P$, $P \notin L_{p_0}$. Thus $J$ contains points outside $L_{p_0}$ and also points of $I_{p_0} \subseteq K_{p_0}$. The dimensions of $L_{p_0}$ are $4n+1$ times the dimensions of $K_{p_0}$. Thus
\[
\mbox{diam}(J)  >  4n\times(\mbox{half edge of } K_{p_0}) > 2\times\mbox{diam}(K_{p_0}) \geq 2\mbox{diam}(I_{p_0}) > d_{p_0-1},
\]
contradicting (\ref{39g}). Thus (\ref{39c}) gives a contradiction, so (\ref{39c}) is false, (\ref{39a}) is true, and the theorem is proved. \nproof

\begin{theorem} \label{t40}
Let $IV(\mu;I;X;\alpha)=0$ Then $\mu^*(X \cap I) = V(\mu;I;X)=0$.
\end{theorem}
\proof
By definition, given $\ve>0$ there are $\delta_m(P)>0$ and $\alpha$-inner covers $\mathcal{C}_m$ compatible with $\delta_m(P),I$ such that
$IV(\mu; \mathcal{C}_m) \leq \ve 2^{-m}$. We can also suppose that $\delta_m(P) \leq m^{-1}$ at all point $P$. Then $\bigcup_{m=1}^\infty \mathcal{C}_m$ is a suitable $\mathcal{F}$ in Theorem \ref{39}, and by definition of $IV_1$, the $I_p$ constructed from $\mathcal{F}$ satisfy
\[
\sum_{p=1}^\infty \mu(I_p) \leq \sum_{m=1}^\infty \frac \ve{2^m} = \ve.
\]
Also by Theorem \ref{39} (\ref{39a}) and Theorem \ref{properties of variation},
\begin{eqnarray*}
V(\mu;I;X) & \leq & V\left(\mu;I;X \setminus \bigcup_{p=1}^\infty I_p\right) + V\left(\mu;I;X  \cap \bigcup_{p=1}^\infty I_p\right) \vt
&= &\lim_{N\rightarrow \infty} V\left(\mu;I;X  \cap \bigcup_{p=1}^N I_p\right).
\end{eqnarray*}
By Theorem \ref{38} $V\left(\mu;I;X  \cap I_p\right) = V(\mu; I_p;X)$ so
\[
\lim_{N \rightarrow \infty} \sum_{p=1}^N V(\mu;I_p;X) \leq \sum_{p=1}^\infty V(\mu;I_p) = \sum_{p=1}^\infty \mu(I_p) \leq \ve.
\]
Thus $V(\mu;I;X) \leq \ve$, giving the result since $\ve>0$ is arbitrary. \nproof

\noindent \textbf{Corollary:} The exceptional sets of Theorems \ref{36} and \ref{37} are of measure zero.

\section{Limits of Step Functions}
We now show that if limit is taken to mean limit almost everywhere, the the set of step functions is dense (in the \emph{``ess sup''} sense) in the set of integrable functions.

First, a function $S(P)$ of the points $P \in I \subseteq E^n$ is a \emph{step function} if there are a division $\D$ of $I$ consisting of $I_1, \ldots , I_m$ and constants $v_1, \ldots v_m$ such that $S(P) = v_j$ in the interior of $I_j$ ($1 \leq j \leq m$). Then whatever finite values $S(P)$ takes on the boundaries of the $I_j$, $S(P)$ is integrable in $I$ with
\[
\int_I S(P)d\mu = \sum_{j=1}^m v_j \mu(I_j).
\]
For proof see proof of Theorem \ref{38}.

\begin{theorem}\label{t41}
Let $f$ be integrable in $I$ to $F$. Then $f$ is the limit a.e.~of a sequence of step functions, the integral of each of these over $I$ being equal to $F$.
\end{theorem}
\proof
We divide the brick $I$ into $2^{kn}$ equal bricks $I_j$ ($j=1, \ldots , 2^{n}$) by continued bisection with respect to each co-ordinate. Let us put
\[
S_k(P) = \frac{\int_{I_j} f d\mu}{\mu(I_j)}
\]
for all $P$ in the interior of $I_j$. Then
\begin{equation} \label{41a}
\int_I S_k(P) d\mu = \sum_{j=1}^{kn} \frac{\int_{I_j} f d\mu}{\mu(I_j)} \mu(I_j) = \sum \int_{I_j} = \int_I fd\mu.
\end{equation}
Also, by Theorems \ref{36}, \ref{40}, $F(J)$ is differentiable with value $f$ except in a set $X$ of measure zero. Here $r(I_j) = r(I) >0$. Therefore if $P \notin X$, $P$ not on any boundary,
\[
S_k(P) = \frac{\int_{I_j} f d\mu}{\mu(I_j)} = \frac{F(I_j)}{\mu(I_j)} \rightarrow f(P).
\]
(Note that a countable number of boundaries of measure zero gives measure zero.)
\nproof

A further result on the differentiation of finitely additive brick functions follows similarly.

\begin{theorem}\label{t42}
Let the real $h(J)$ be finitely additive in the bricks $J \subseteq I$ and of bounded variation in $I$. Then
\begin{equation}\label{42a}
f(P) =D(h;\mu;I;P)
\end{equation}
exists everywhere in $I$, and if $f(P)=0$ in the exceptional set then
\begin{equation}\label{42b}
f \;\mbox{ and }\;|f|\;\mbox{ are integrable in }\;I\;\mbox{ with }\;\int_I |f|d\mu \leq H(h;I).
\end{equation}
\end{theorem}
\proof
First, if $\D$ is any division of a brick $J$, finite additivity gives
$
|h(J)| = |(\D)\sum h| \leq (\D) \sum |h|$, so
\begin{equation}\label{42c}
|h(J)| \leq V(h;J) \;\;\mbox{ for all }\;\; J \subseteq I.
\end{equation}
Thus we can write $h$ as the difference between two non-negative finally additive brick functions,
\begin{equation}\label{42d}
h(J) = \frac 12\left(V(h;J) +h(J)\right) - \frac 12\left(V(h;J) - h(J)\right).
\end{equation}
Therefore it is enough to prove the result when $h \geq 0$.
Given numbers $M>0$, $\alpha>0$, suppose that a set $X$ has the property that for arbitrary $\delta(P)>0$ in $I$ and some $\alpha$-inner cover of $X$ that is compatible with $\delta(P)$ and $I$ we have $h(J) \geq M\mu(J)$. We then prove
\begin{equation} \label{40e}
V(h;I;X) \geq MV(\mu;I;X).
\end{equation}
Given $\ve>0$ there is a $\delta(P)>0$ depending on $\ve$ such that
\[
V(h;I;X;\delta) < V(h;I;X) + \ve.
\]
By Theorems \ref{38} and \ref{39} (Vitali) there is a finite collection $\mathcal{Q}$ of disjoint bricks from the various $\alpha$-inner covers, with union $U$, such that
\begin{eqnarray*}
M \left(V(\mu;I;X) - \ve \right) & \leq & M.V(\mu;I;X\cap U)\;\;\leq \;\;\ M.V(h;I;U) \vt
&= &M. (\mathcal{Q})\sum \mu(J)\;\; \leq \;\; M. (\mathcal{Q})\sum H \vt
& \leq &V(h;I;X;\delta)
\;\;<\;\; V(h;I;X) + \ve.
\end{eqnarray*}
Thus (\ref{40e}) follows as $\ve \rightarrow 0$. Therefore if for $X$ we have $h \geq M\mu$ for one $(\delta(P),I)$-compatible $\alpha$-inner cover and $h \leq N\mu$ for another, where $M>N$, and where $\delta(P)$ is arbitrary, we have
\[
M.V(\mu;I;X) \leq V(h;I;X) \leq N.V(\mu;I;X).
\]
Therefore, with $M>N$,
\begin{equation}\label{42f}
V(\mu;I;X) = 0 .
\end{equation}
Therefore each set for which $M=(m+1)2^{-p}$, $N= m2{-p}$ ($m,p=0,1,2, \ldots$) satisfies (\ref{42f}). Therefore the union does for each fixed $\alpha>0$. Taking $\alpha = q^{-1}$ for $q=1,2,3, \ldots$, the union of the corresponding sets also satisfies (\ref{42f}). In the complement of this union (no oscillation) the derivative exists. As in Theorem \ref{41}, $\frac{h(J)}{\mu(J)} \rightarrow f(P)$. As in Theorem \ref{41} we divide $I$ up into $2^{kn}$ equal bricks $I_j$ ($j=1,2, \ldots , 2^{kn}$), putting
\[
s_k(P) = \frac{h(I_j)}{\mu(I_j)} \;\mbox{ for }\;P\;\mbox{ inside }\;I_j\;\;\;j=1,2, \ldots , 2^{kn}.
\]
Then
\[
\int_I s_k(P) d\mu = \sum_{j=1}^{2^{kn}}  \frac{h(I_j)}{\mu(I_j)} \mu(I_j) = \sum_{j=1}^{2^{kn}}  {h(I_j)} = h(I).
\]
Also $s_k(P) \geq 0$ and tends to $f(P)$ outside $X$, except on the boundaries of the bricks, i.e.~ $s_k(P) \rightarrow f(P)$ almost everywhere. By Fatou's lemma,
\[
\int_I f d\mu = \int_I \lim_{k \rightarrow \infty} s_k(P) d\mu \leq \liminf_{k \rightarrow \infty} \int_I s_k(P) d\mu = h(I).
\]
This result occurs when $h \geq 0$. For arbitrary $h$, $D(h;\mu;I;P)$ is a difference $f_1-f_2$ of two non-negative integrable functions, existing save in an exceptional set of outer measure zero, while
\[
\int_I|f|d\mu \leq \int_I f_1d\mu + \int_I f_2d\mu \leq \frac 12\left(V(h;I) +h(I)\right) + \frac 12\left(V(h;I) -h(I)\right) = V(h;I),
\]
giving the result. \nproof

\section{Absolutely Continuous Functions}
Let $h(J)$ be a function of bricks $J \subseteq I$. We say that $h$ is \emph{absolutely continuous} in $I$ if $V(h;I;X) \rightarrow )$ when $V(\mu;I;X) \rightarrow 0$. (The latter is the outer measure of $X$.) More exactly, we put $H(\ve) = \sup \{V(h;I;X)\}$, the $\sup$ being taken for all sets $X \subseteq I$ with $V(\mu;I;X) \leq \ve$. Then $h$ is absolutely continuous in $I$ if $H(\ve) \rightarrow 0$ as $\ve \rightarrow 0$.

\begin{theorem}\label{t43}
Let $f$ be integrable to $F(J)$ for bricks $J \subseteq I$, and also let $|f|$ be integrable in $I$. Then $F$ is absolutely continuous.
\end{theorem}
\proof
As $M + |f|$ is integrable for any constant $M$, Theorem \ref{20} implies that
\[
f_M \equiv \max \{ M, |f|\} - M
\]
is integrable in $I$.
Let $X_M$ be the set where $f_M \neq 0$. Then
\[
f_M = \left( |f| - M\right) \chi(X_M;P) \rightarrow 0
\]
as $M \rightarrow \infty$, and is bounded by $|f|$. Therefore by Levi's monotone convergence theorem, $\int_I f_M d\mu \rightarrow 0$ as $M \rightarrow \infty$. Also $0 \leq |f|-f_M \leq M$. Therefore, given $\ve>0$, we can choose $M$ so that
\[
V(f_M \mu;I) = \int_I f_M d\mu < \frac 12 \ve.
\]
(The equality occurs since $f_M\mu \geq 0$ and since the integral exists.) Then for each $X$ with $V(\mu;I;X) < \ve 2^{-m}$ we have
\[
V(f\mu;I;X) = V\left((|f| - f_M)\mu+f_M \mu;I;X\right).
\]
Therefore
\begin{eqnarray*}
V(f\mu;I;X) & \leq &  V\left((|f| - f_M)\mu;I;X\right) +  V\left(f_M \mu;I;X\right) \vt
&\leq & M.V(\mu;I;X) + V(f_M \mu;I) \;\;<\;\;\frac \ve 2+\frac \ve 2 \;\;=\;\;\ve.
\end{eqnarray*}
Since $\ve>0$ is arbitrary, $V(f\mu;I;X)=0$. Further, $V(f\mu-F;I)=0$ so that
\[
V(F;I;X) \leq V(F - f\mu;I;X) + V(f\mu;I;X) < 0+\ve.
\]
This proves that $F$ is absolutely continuous. \nproof

We also prove a converse.

\begin{theorem} \label{t44}
Let $h(J)$ be finitely additive for $J \subseteq I$ and absolutely continuous. Then $h$ is the integral of its derivative, the modulus of which is also integrable.
\end{theorem}
\proof
First we show that $h$ is of bounded variation. We have $V(h;I;X) < \delta$ as soon as $V(\mu;I;X) < \ve$, where $\delta>0$ depends on $\ve>0$. Dividing up $I$ into $2^{kn}$ bricks $I_j$ in the usual way, we have
\[
V(\mu;I_j;X) = \mu(I_j) = \frac{\mu(I)}{2^{kn}} < \ve
\]
for $k$ sufficiently large. Therefore $h$ is of bounded variation. Hence, by Theorem \ref{42}, $D(h;\mu;I;P)$ exists almost everywhere, and is integrable, to the value $F(J)$, say, over brick $J$, and $|D|$ is also integrable. By Theorems \ref{36} \ref{40}, $F(J)$ is differentiable a.e.~to $D(h;\mu;I;P)$. Hence $D(h-F;\mu;I;P)=0$ a.e., i.e.
\[
|h(J) -F(J) |<\ve\mu(J)
\]
for all $J$ compatible with some $\delta(P)>0$ depending on $\ve>0$ and almost everywhere in $P$. As in Theorem \ref{42}, taking $X$ as the exceptional set,
\[
V(h-F;I; \setminus X) \leq \ve V(\mu;I;\setminus X) \leq \ve \mu(I)
\]
for all $\ve>0$. Hence $V(h-F;I;\setminus X) =0$. Also $V(\mu;I;X)=0$ so $V(h;I;X)=0$. Further, by Theorem \ref{43} $V(F;I;X)=0$. Hence $V(h-F;I;X)=0$, and so $V(h-F;I)=0$. Both $h(J)$ and $F(J)$ are finitely additive, so
\[
|h(J)-F(J)| \leq V(h-F;;J) \leq V(h-F;I)
\]
for all $J \subseteq I$, so $h(J) = F(J)$. Thus $h(J) = \int_J D(h;\mu;I;P)d\mu$. \nproof

\section{Connections between Variation and Leb\-es\-gue Outer Measure}
In Lebesgue theory we first have to find the measure of an open set $G$, and this is done as follows. We begin with a set of $(n-1)$-dimensional hyperplanes, each with one variable equal to a constant, say $x_j =c$ for some $j$ in $1 \leq j \leq n$. All the planes with $c$ an integer can be used to cut up the $n$-dimensional space $E^n$ into bricks $J$.

Let $\mu_1(G)$ be the sum of the volumes of the $J$ for all such $J \subseteq G$, the union $H_1$ of the $J$ being taken closed.
If $\mu_1(G) = +\infty$ we put $\mu(G)=+\infty$. Otherwise $G \setminus H_1$ is open.

Next we consider all planes with $2c$ an integer, obtaining smaller bricks $J$. Let $\mu_2(G)$ be the sum of the volumes of all such $J \subseteq G \setminus H_1$, the union $H_2$ of the $J$ being taken closed. If $\mu_2(G) = +\infty$ we put $\mu(G) = +\infty$. Otherwise $G \setminus (H_1 \cup H_2)$ is open and we consider all planes with $2c$ an integer, and so on.

We thus define $\mu_j(G)$ when $\mu_m(G)$ is finite for all $m < j$. Putting
\[
\mu(G) = \sum_{j=1}^\infty \mu_j(G)
\]
when all $\mu_j(G)$ are finite, we have $\mu(G) $ for all open sets $G$. By convention we can assume that $\mu(G) = \sum_{j=1}^\infty \mu_j(G)$ for all $G$.

Since at each stage the $J$'s can be put in sequence, it follows that all the $J$'s can be put in sequence, and if the sequence is $\{J_j\}$ we have
\[
\mu(G) = \sum_{j=1}^\infty \mu(J_j),
\]
as this is just a rearrangement of the expanded form of the series $\sum \mu_j(G)$ using the $\mu(J)$. If $G \subseteq I$ then
\[
V(\mu;I;G) \leq \sum_{j=1}^\infty V(\mu;I;I_j) = \sum_{j=1}^ \infty \mu(J_j) = \mu(G).
\]
Conversely,
\[
\sum_{j=1}^N \mu(J_j) = V\left(\mu;I;\bigcup_{j=1}^NJ_j\right)\leq V(\mu;I;G).
\]
Letting $N \rightarrow \infty$ we see that, for $G \subseteq I$,
\[
\mu(G) = V(\mu;I;G).
\]
More generally, if $I^0$ is the interior of $I$, then $V(\mu;I;G) = \mu(G \cap I^0)$
whether or not $G \subseteq I$.

Let $X$ be a set in $E^n$. Then the \emph{Lebesgue outer measure} $\mu^*_L(X)$ of $X$ is the infimum of $\mu(G$ for all $G \supseteq X$. As
\[
V(\mu;I;X) \leq V(\mu;I;G) = \mu(G \cal I^0)
\]
we see that
\[
V(\mu;I;X) \leq \mu^*_L (X \cap I^0).
\]
To prove the converse we need Vitali's theorem for Lebesgue outer measure and the proof of this is the same as the proof of Theorem \ref{35}.
Thus let $\delta(P)>0$ be given in $I$, and let $\alpha>0$. Then there is a finite or countable sequence $\{I_p\}$ of mutually disjoint bricks $I_p \subseteq I$ with associated points in $X$, and compatible with $\delta(P)>0$, and with $r(I_p) \geq \alpha$, such that
\[
\mu^*_L\left(X \cap I^0 \setminus \bigcup_{p=1}^ \infty I_p \right) =0.
\]
Then
\[
\mu^*_L (X \cap I^0) \leq \sum_{p=1}^\infty \mu(I_p),\;\;\mbox{ and }\;\;\sum_{p=1}^N \mu(I_p) \leq V(\mu;I;X;\delta)
\]
for each $N$, independent of $\delta$. Hence $\mu^*_L(X \cap I^0) \leq V(\mu;I;X;\delta)$,
\[
\mu^*_L(X \cap I^0) \leq V(\mu;I;X) ,\;\;\mbox{ so }\;\;\mu^*_L(X \cap I^0) = V(\mu;I;X).
\]
To illustrate in another way the connection between outer measure and variation, let us call $f(P)$ a \emph{null} function if there is a sequence $\{X_j\}$ of sets with $\mu^*_L(X_j)=0$ such that
\[
|f(P)| \leq 2^j\;\;\;\;(P \in X_j);\;\;\;\;\;\;\;\;\;\;f(P) =0\;\;\;\;(P \notin X=\bigcup_{j=1}^\infty X_j).
\]
For example, in $E^n$ take $n=1$, and let $f_1(P)$ be the characteristic function of the rationals. These can be put in a sequence $\{r_j\}$ by writing
\[
0,1,-1, \frac 12, -\frac 12, \frac 21, -\frac 21, \frac 13, -\frac 13, \frac 22, - \frac 22, \frac 31, -\frac 31, \frac 14, -\frac 14, \frac 23, -\frac 23, \ldots.
\]
Crossing out the second and later appearances of each rational, we have the required sequence $\{r_j\}$. Let
\[
I_j = \left( r_j - \frac \ve{2^{1+j}}, \;r_j + \frac \ve{2^{1+j}}\right),\;\;\;\;\;\;\;\;\;\;
G=\bigcup_{j=1}^\infty I_j.
\]
Then
\[
\mu(G) \leq \sum_{j=1}^\infty \mu(I_j) = \sum_{j=1}^\infty \frac \ve{2^j} = \ve,
\]
so the set of rational numbers has Lebesgue (outer) measure $0$, and $f_1(P)$ is a null function.

\begin{theorem}\label{t45}
The integral of a null function is zero.
\end{theorem}
\proof
(We could use the result $V(\mu;I;X) \leq \mu^*_L(X \cap I^0)$, but instead we prove it directly.)
Let $\ve>0$. As $\mu^*_L(X_k)=0$ there is an open set $G_k \supseteq X_k$ with $\mu(G_k) < \ve 4^{-k}$. Then
\[
G_k = \bigcup_{j=1}^\infty I_{kj}\;\;\;\;\mbox{ with }\;\;\;\;\sum_{j=1}^\infty \mu(I_{kj}) < \frac \ve{4^k}.
\]
Let $1<\delta(P)\leq 1$ for $P \notin X$; and for $P \in Y_k = X_k \setminus \bigcup \{X_m: m<k\}$. Then
\[
J \subseteq \left(P-\delta(P), P+\delta(P) \right) \subseteq G_k = \bigcup_{j=1}^\infty I_{kj},
\]
the $I_{kj}$ being open intervals. As $J$ is closed, Borel's covering theorem shows that $J$ lies in the union of a finite number of the $I_{kj}$ for the fixed $k$, so
\[
\mu(J) \leq \sum_{j=1}^\infty \mu\left(I_{kj} \cap J\right).
\]
As the only non-zero $f(P)$ have $P\in X = \bigcup_{k=1}^\infty Y_k$, we have
\begin{eqnarray*}
\left| (\D)\sum f(P)\mu(J) -0\right| & \leq & \sum_{k=1}^\infty 2^k \left(\sum_{j=1}^\infty (\D) \sum \mu\left(I_{kj} \cap J\right)\right) \vt
&\leq & \sum_{k=1}^\infty 2^k \left( \sum_{j=1}^\infty\mu\left(I_{kj}\right)\right) \;\;\leq\;\; \sum_{k=1}^\infty 2^k \frac \ve{4^k} \vt
 &=&\sum_{k=1}^\infty \frac \ve{2^k} \;\;=\;\; \ve,
\end{eqnarray*}
giving the result. \nproof

\section{Lebesgue Integration is Included in Ours}
The proof is from a paper by R.O.~Davies and Z.~Schuss  in which they prove a rather stronger result.

Let $f$ be a finite real function, Lebesgue integrable over a brick $I$ with $(L)\int_I f(P)d\mu =F$. Then the Lebesgue integral is absolutely continuous, so that, given $\ve>0$, there is an $\eta>0$ such that for measurable sets $A \subseteq I$,
\[
\mu(A)<\ve\;\;\;\;\mbox{ implies }\;\;\;\;(L)\int_A |f|d\mu < \ve 2^{-1}.
\]
For $m=0, \pm1, \pm2, \ldots$, and $\xi= \ve\left(3(\eta + \mu(I))\right)^{-1}$, let $X_m$ be the set of $P$ in $I$ where
\[
(m-1)\xi < f(P)\leq m\xi.
\]
 Then we can choose for each $m$ an open set $G_m \supseteq X_m$ such that
 \[
 \mu^*\left(G_m \setminus X_m\right) < \frac \eta{2^{|m|+2}(|m|+1)}.
 \]
 Finally we can choose $\delta(P)>0$ in such a way that if $P \in X_m$ then the sphere with centre $P$ and radius $\delta(P)$ lies in $G_m$.

 Let $\D$ be a $\delta$-compatible division of $I$ and let the brick $J \in \D$ have associated point $P \in X_m$ for some $m$, say $m(P)$. Then
 \[
 J \subseteq G_{m(P)}\;\;\;\;\mbox{ and }\;\;\;\;J \setminus X_{m(P)} \setminus G_{m(P)} \setminus X_{m(P)},
 \]
 and, with $P$ fixed in $J \in \D$, and $Q$ the point-variable in $(L)\int_J$,
 \begin{eqnarray*}
 \left|(\D)\sum f(P)\mu(J) -F\right| &=&\left|(\D)\sum (L)\int_J \left(f(P) -f(Q)\right) d\mu\right| \vt
 &\leq & (\D)\sum (L)\int_J \left| f(P)-f(Q)\right| d\mu \vt
 & \leq & (\D)\sum (L)\int_{J \cap X_{m(P)}} \left| f(P) - f(Q)\right|d\mu\;\;+ \vt
 &&\;\;\;\;\;+\;\;(\D)\sum (L)\int_{J \setminus X_{m(P)}} \left| f(P) \right|d\mu \;\;+\vt
&&\;\;\;\;\;\;\;\;\;\;+\;\; (\D)\sum (L)\int_{J \setminus X_{m(P)}} \left|  f(Q)\right|d\mu \vt
 &=& R+S+T,\;\;\;\mbox{ say}.
 \end{eqnarray*}
In $R$ we have $Q \in J \cap X_{m(P)}$, and so $f(Q)$ and $f(P)$ lie in the same interval $((m(P)-1)\xi, m(P)\xi]$, so
\[
R \leq (\D)\sum (L)\int_{J\cap X_{m(P)}} \xi d\mu \leq (\D)\sum \mu(J) = \xi \mu(I) < \frac \ve 3.
\]
In $S$ we collect those terms (if any) for which $m(P)$ has a given value $m$, and we write
\[
S = \sum_{m=-\infty}^\infty \left( (\D) \sum_{m(P)=m} (L)\int_{J \setminus X_m} |f(P)| d\mu\right),
\]
the inner sum being empty for all but a finite number of $m$, and
\begin{eqnarray*}
S & \leq &
\sum_{m=-\infty}^\infty (\D) \sum _{m(P)=m} \left(|m|+1\right) \xi \mu_L (J \setminus X_m)  \vt
&\leq &
\sum_{m=-\infty}^\infty (|m|+1) \xi \mu_L(G_m \setminus X_m) \vt
&<& \sum_{m=-\infty}^\infty \frac{(|m|+1)\ve}{3(\eta + \mu(I))} \times \frac \eta{2^{|m|+2}(|m|+1)} \vt
&<& \sum_{m=-\infty}^\infty \ve{3\left(2^{|m|+2}\right)} \;\;=\;\; \frac \ve 3\left( \frac 14+ \frac 14 + \frac 14\right) \;\;<\;\;\ve.
\end{eqnarray*}
Finally, $T = (L)\int_A |f(Q)| d\mu$ where $A= \bigcup\{J \setminus X_{m(P)}\;:\; J \in \D\}$,
\[
\mu(A)  \leq \sum_{m=-\infty}^\infty \left( (\D)\sum \mu(J \setminus X_m)\right) \leq \sum_{m=-\infty}^\infty \mu(G_m \setminus X_m) < \eta.
\]
Thus $T < \frac 13 \ve$, and $0 \leq R+S+T<\ve$. So $F=(L)\int_I fd\mu$ is the value of the generalised Riemann integral over $I$.

\section{The Denjoy Extension}
Two methods are used to extend the definition of the Lebesgue integral to become the special denjoy integral. First there is the Cauchy extension that produces the H\"{o}lder-Lebesgue integral, which is non-absolute.

In Euclidean one-dimensional space $E^1$, if the Lebesgue integral over $[b,c]$ does not exist, but the Lebesgue integral over $[u,v]$ exists for all $u,v$ in $b<u<v<c$, we define the integral over $[b,c]$ to be
\[
\lim_{u \rightarrow b+} \lim_{v \rightarrow c-} (L)\int_u^v f(x)dx.
\]
Since $(L)\int_u^v f(x)dx = \int_u^v f(x)dx$, this new integral has value $\int_b^c f(x)dx$.

There is another extension, due to Denjoy, and we shall again restrict the discussion to $E^1$, in which an open set $G$ is the union of a sequence of disjoint open intervals $(b_j,c_j)$. When the Lebesgue integral exists over $[b,c] \setminus G$, a closed set, and when the H\"{o}lder-Lebesgue integrals over the $(b_j,c_j)$ exists, we can sometimes define the special Denjoy integral to be
\[
(L)\int_{[b,c]\setminus G} fd\mu + \sum_{j=1}^\infty (HL)\int_{(b_j,c_j) \cap [b,c]} fd\mu.
\]
We now turn to the corresponding theorem in generalised Riemann integration.

\begin{theorem}\label{t46}
In $E^1$, if $G = \bigcup_{j=1}^\infty (b_j,c_j)$ is an open set in a finite interval $[b,c]$, the $(b_j,c_j)$ being disjoint, if $f(x)\mu(J) \chi(\setminus G,x)$ is integrable in $[b,c]$, if $f(x)\mu(J)$ is integrable over each $[b_j,c_j]$ ($j=1,2$), and if, given $\ve>0$, there is an integer $N$ such that for every finite collection $\mathcal{Q}$ of intervals $[u,v]$, each contained in a $[b_j,c_j]$ with $j \geq N$, no two intervals $[u,v]$ lying in the same $[b_j,c_j]$, we have
\begin{equation} \label{46(1)}
\left|(\mathcal{Q}) \sum \int_u^v f(x)d\mu \right| < \ve,
\end{equation}
then there exists
\begin{equation} \label{46(2)}
\int_b^c f(x)d\mu =\int_b^c \chi(\setminus G,x) f(x) d\mu + \sum_{j=1}^\infty \int_{b_j}^{c_j} f(x) d\mu.
\end{equation}
\end{theorem}
\proof
In (\ref{46(2)}), subtracting the first integral on the right from the left hand side, we have to prove that there exists
\begin{equation} \label{46(3)}
\int_b^c \chi(G,x) f(x)d\mu = \sum_{j=1}^\infty \int_{b_j}^{c_j} f(x)d\mu.
\end{equation}
Therefore define
\begin{equation} \label{46(4)}
H_1(u',v') = \sum_{j=1}^\infty \int_{[b_j,{c_j}] \cap [u',v']} f(x)d\mu
\end{equation}
(putting $\int =0$ if the intersection ${[b_j,{c_j}] \cap [u',v']}$ is empty or a single point). First show that the series is convergent. Here ${[b_j,{c_j}] \cap [u',v']}$ is either empty or $[b_j,c_j]$ for all but at most two values of $j$ for which $[b_j,c_j]$ contains $u'$ or $v'$ or both. From (\ref{46(1)}), the sequence of partial sums for (\ref{46(4)}) is fundamental, and is convergent, and $H(u,v)$ exists for all $u,v$ in $b \leq u<v \leq c$. Next there is a $\delta_j(x)>0$ defined in $[b_j,c_j]$ for which every $\delta$-compatible division $\D_j$ of $[b_j,c_j]$ satisfies
\begin{equation} \label{46(5)}
\left| (\D_j) \sum f(x)\mu(J) - H_1(b_j,c_j)\right| < \frac \ve{2^j} \;\;\;\;\;\;\;\;(j=1,2,3, \ldots),
\end{equation}
and so, by a previous theorem (Theorem \ref{Henstock Lemma}),
\begin{equation} \label{46(6)}
(\D_j) \sum \left|f(x)\mu(J) - H_1(b_j,c_j)\right| < \frac {4\ve}{2^j} \;\;\;\;\;\;\;\;(j=1,2,3, \ldots).
\end{equation}
From the separate $\delta_j(x)$ we construct a $\delta(x)$ defined in $[b,c]$ in the following way. Let $N$ be an integer satisfying (\ref{46(1)}). If $x \in G$ then, for some $j$, $b_j <x<c_j$ and we take $\delta(x)>0$ satisfying
\begin{equation} \label{46(7)}
(x-\delta(x),x+\delta(x)) \subseteq (b_j,c_j),\;\;\;\;\;\; \delta(x) \leq \delta_j(x).
\end{equation}
On the other hand, if $x \in \setminus G$ we take a $\delta(x)>0$ such that
\begin{eqnarray}
&&\!\!\!\!\!\!\!\!\!\!\!\!\!\!\!\!\!\!\!\!\!\!\!\!\!\!\!\mbox{if}\;x =c_j\;\mbox{for some}\;j\;\mbox{then}\; (x-\delta(x),x) \subseteq (b_j,c_j),\;\delta(x) \leq \delta_j(x); \label{46(8)}\vt
&&\!\!\!\!\!\!\!\!\!\!\!\!\!\!\!\!\!\!\!\!\!\!\!\!\!\!\!\mbox{if}\;x \neq c_j\;\mbox{for any}\;j,\;\mbox{take}\; (x-\delta(x),x) \cap (b_k,c_k) = \emptyset,\;\;k=1,\ldots,N-1; \label{46(9)}
\vt
&&\!\!\!\!\!\!\!\!\!\!\!\!\!\!\!\!\!\!\!\!\!\!\!\!\!\!\!\mbox{if}\;x =b_j\;\mbox{for some}\;j\;\mbox{then}\; (x,x+\delta(x)) \subseteq (b_j,c_j),\;\delta(x) \leq \delta_j(x); \label{46(10)} \vt
&&\!\!\!\!\!\!\!\!\!\!\!\!\!\!\!\!\!\!\!\!\!\!\!\!\!\!\!\mbox{if}\;x \neq b_j\;\mbox{for any}\;j,\;\mbox{take}\; (x,x+\delta(x)) \cap (b_k,c_k)=\emptyset,\;\;k=1,\ldots,N-1. \label{46(11)}
\end{eqnarray}
Let $\D$ be a division of $[b,c]$ compatible with $\delta(x)$. If $x$ is an associated point of an interval $[u,v]$ of $\D$ and if $u<x<v$ then
\[
f(x)(v-u) = f(x)(v-x) + f(x)(x-u),
\]
so we can assume that the associated point is an end point of its interval. If $x=v=c_j$ then $b_j \leq u<c_j$ by (\ref{46(8)}). If $x=v \neq c_j$, and $x \in \setminus G$, then,  by (\ref{46(9)}), $[u,v]$ can only overlap with $[b_j,c_j]$ when $j \geq N$. If $x=v \in G$ then, for some $j$, $[u,v] \subseteq [b_j,c_j]$ by (\ref{46(7)}). Similarly when $x=u$.
Using (\ref{46(9)}) and (\ref{46(11)}) and all $[u,v] \in \D$ with $x \in \setminus G$, $x=u \neq b_j$ or $x=v \neq c_j$ (all $j$), (*) the sum of $H_1(u,v)$ is the limit of a sequence of sms over various $\mathcal{Q}$ satisfying (\ref{46(1)}), and the modulus of the sum is $\leq \ve$. Here $x \notin G$ so
$
H_1(J) - f(x) \chi(G,x) \mu(J) = H_1(J).
$
If $x \in \setminus G$ and $x=u=b_j$ or $x=v=c_j$ for some $j$, (**) or if $x \in G$, then by (\ref{46(7)}), (\ref{46(8)}), (\ref{46(10)}), $[u,v]$ is compatible with $\delta_j(x)$. Thus we see that $[b_1, c_1], \ldots ,[b_{N-1}, c_{N-1}]$ are each divided by partial divisions of $\D$, using the $[u,v]$ of (**), while the parts of $[b_j,c_j]$ ($j \geq N$) are covered. Using (\ref{46(4)}), (\ref{46(5)}), (\ref{46(6)}) and the remarks of (*) and (**) we have
\begin{eqnarray*}
\left| (\D)\sum \chi(G,x) f(x) \mu(J) - H_1(b,c) \right| &=& \left| (\D)\sum \left(\chi(G,x) f(x) \mu(J) - H_1(J)\right) \right| \vt
&\leq & \ve + \sum_{j=1}^{N-1} \frac \ve{2^j} + \sum_{j=N}^{\infty} \frac {4\ve}{2^j}\;\;< \;\;5 \ve,
\end{eqnarray*}
which proves the theorem. \nproof

\section{The Radon-Nikodym Theorem}
Let $h_1(I,P)$ and $h(I,P)$ be brick-point functions.  Then $h_1$ is \emph{absolutely continuous} in an elementary set $E$ with respect to $h$ if, given $\ve>0$, there is a $\delta>0$ such that all sets $X$ with $V(h;E;X)<\delta$ also have  $V(h_1;E;X)<\ve$. We look for a point function $f$ with the property that $V(h_1 - fh;E) =0$.

\begin{theorem}\label{t47}
If $h(I,P)$ is real or complex, integrable with integral $H$ in $E$, and of bounded variation, then $|h(I,P|$ and $|H(I)|$ are integrable to $V(E_1) \equiv V(h;E_1)$ for each each elementary set $E_1 \subseteq E$. Also $V(E_1) = V(H,E_1)$ and $|H(E_1)| \leq V(E_1)$.
\end{theorem}
\proof
For divisions $\D$ of elementary sets $E_1 \subseteq E$,
\begin{equation}\label{47(1)}
\left|(\D)\sum |h(I,P)| -(\D) \sum |H(I)|\right| \leq (\D) \sum \left| |h| - |H| \right| \leq (\D) \sum |h-H|.
\end{equation}
This is as small as we please by choice of $\delta(P)>0$ at each point of $E$ and for all divisions $\D$ of $E_1$ compatible with $\delta(P)$. Hence
\begin{equation}\label{47(2)}
V(H;E_1) = V(h;E_1) \equiv V(E_1).
\end{equation}
As $H$ is finitely additive, $|H|$ is finitely subadditive and so is integrable, and its integral is $V(H;E_1)$ since $|H| \geq 0$. Finally, by subadditivity,
$|H(E_1)| \leq \int_{E_1} |H| = V(E_1)$. \nproof

\begin{theorem}\label{t48}
Let $h(I,P)$ be real or complex, integrable with integral $H$ in $E$, and of bounded variation. For each division $\D(I)$ of each elementary set $E_1 \subseteq E$, let there be disjoint sets $J \subseteq I$ with union $E_1$ (note: the edges of the $I$ might overlap but the $J$ do not), such that if $g(\D;P)$ is a function of points that is constant at a value $g(I)$ in the set $J \subseteq I$ for each $I \in \D$, then we have
\begin{equation}\label{48(3)}
\int_{E_1} g(\D;x) dV = (\D) \sum g(I) V(I)
\end{equation}
where $V(I) = V(h;E;I)$.
Also, for each elementary set $E_1 \subseteq E$ let
\begin{equation}\label{48(4)}
V(h;E_1) = V(h;E;E_1).
\end{equation}
then there is an $f(P)$ in $E$ with $|f(P)|=1$ such that, for all bricks $I^* \subseteq E$,
\[
H(I^*) = \int_{I^*} f(x)dV.
\]
\end{theorem}
\proof
Using
\[
\mbox{sgn} (z) = \left\{
\begin{array}{lll}
\frac z{|z|} & \mbox{for} & z \neq 0,\vt
1 & \mbox{for} & z = 0,
\end{array}
\right.
\]
and $g(I) = \mbox{sgn}(H(I))$, then, for $H=0$, $|H-gV| = |V|$, and for $H \neq 0$,
\[
|H-gV| = \left|H - \frac H{|H|} V\right| = \left|\frac H{|H|}|H| - \frac H{|H|} V\right|= \left| |H|-V \right|.
\]
Hence, by Theorem \ref{47}, $V(|H| -V;E) =0$ and
\begin{equation}\label{48(5)}
H(I) = \int_I g dV,\;\;\;\;\;\;|g| =1
\end{equation}
for all bricks $I \subseteq E$. We now have to replace $g$ by a point function $f(x)$ independent of $I$, with $|f(x)| = 1$. The notation $\D'' \leq \D'$ means that every $I' \in \D'$ is a union of $I'' \in \D''$. First there is a sequence $\{\D_j\}$ of divisions of $E$ where $\D_j$ is compatible with $\delta_j(P)>0$ suitable chosen so that
\begin{eqnarray}
&&\delta_{j+1}(P) \;\; \leq \;\; \delta_j(P), \label{48(6)} \vt
&& \D_{j+1} \;\;\leq \;\;\D_j \;\;\mbox{ and every }\;I^* \in \D_{j+1}\;\mbox{ lies in an }\;I \in \D_j, \label{48(7)} \vt
&& (\D_j)\sum |H| \;\;>\;\;V(E) - \frac 1{2^{4j}}, \label{48(8)} \vt
&& (\D') \sum |gV -H| \;\;<\;\; \frac 1{2^j} \label{48(9)}
\end{eqnarray}
for all divisions $\D'$ of $E$ compatible with $\delta_j(P)>0$ ($j \geq 1$). Using (\ref{48(8)}), (\ref{48(9)}) with $\D'=\D_j$ and $f_j(P) = g(I)$ ($P \in J \subseteq I$ and all $I \in \D_j$) we have
\begin{equation} \label{48(10)}
|f_j(P)|=1,\;\;\;\;\;\;\;\left| \int_E f_j(P) dV - H(E) \right| < \frac 1{2^j}\;\;\;\;\;(j \geq 1).
\end{equation}
Let $E_j$ be the union of those $I \in \D_j$ with
\begin{equation} \label{48(11)}
|H| \;\;<\;\;\left(1-\frac 1{2^{3j}}\right) V,
\end{equation}
these $I$ forming a collection $\mathcal{Q}_j$. By (\ref{48(8)}), the finite additivity of the variation, (\ref{48(4)}), and Theorem \ref{47},
\begin{eqnarray}
V(E) - \frac 1{2^{4j}} &<& (\D_j)\sum |H| \;\;=\;\;(\D_j \setminus \mathcal{Q}_j ) \sum |H| + (\mathcal{Q}_j ) \sum |H| \nonumber \vt
& \leq & (\D_j \setminus \mathcal{Q}_j ) \sum V \;\;=\;\;V(E) -\frac 1{2^{3j}} V(E_j),\nonumber \vt
V(h;E_j;E) &=& V(E_j) \;\;\leq \;\;\frac 1{2^j},\nonumber\vt
V\left(h;E;\bigcup_{j \geq k} E_j \right) &\leq & \frac 1{2^{k-1}} . \label{48(12)}
\end{eqnarray}
If a division $\D'(L)$ of $E$ has $\D'\leq\D$, let $\mathcal{U}$ be the family of those $L$ which are not subsets of $ E_j$ with
\begin{equation} \label{48(13)}
\Re \left(\frac{H(L)}{g(I)} \right) \leq \left(1-\frac 1{2^{2j}}\right)|H(L)|\;\;\;\;\;\;(L \subseteq I \in \D_j).
\end{equation}
As (\ref{48(11)}) is false for the $I \notin \mathcal{Q}_j$, and by definition of $g(I)$, and Theorem \ref{47},
\begin{eqnarray}
\left(1-\frac 1{2^{3j}} \right) V(I) &\leq & |H(I)| \;\;=\;\; \Re \left( (\D' \cap I)\sum \frac{H(I)}{g(I)} \right)  \nonumber \vt
&\leq & (\mathcal{U} \cap I) \sum |H(L)| \left(1 - \frac 1{2^{2j}}\right) + (\D' \cap I \setminus \mathcal{U}) \sum |H(L)|  \nonumber \vt
& \leq & (\mathcal{U} \cap I) \sum V(L) \left( 1 - \frac 1{2^{2j}} \right)  + (\D' \cap I \setminus \mathcal{U}) \sum V(L) \nonumber  \vt
&=& V(I) - \frac 1{2^{2j}}(\mathcal{U} \cap I) \sum V(L). \nonumber \vt
(\mathcal{U} \cap I) \sum V(L)&\leq & \frac{V(I)}{2^j},\;\;\;\;\;\;\;\;\;\;(\mathcal{U} ) \sum V(L)\;\;\leq\;\;\frac{V(E)}{2^j}. \label{48(14)}
\end{eqnarray}
For the union $W$ of the $L \in \mathcal{U}$, with $W=W_j$ when $\D' = \D_{j+1}$, then, by (\ref{48(4)}), (\ref{48(14)}).
\begin{equation}\label{48(15)}
V(h;E;W) =V(W) < \frac{V(E)}{2^j},\;\;\;\;\;\;\;\; V\left(h;E;\bigcup_{j\geq k}W_j \right) < \frac{V(E)}{2^{k-1}}.
\end{equation}
For the $L$ which are not subsets of $E_j \cup W_j$, (\ref{48(13)}) is false; that is,
\begin{eqnarray*}
&&\Re \left(\frac{f_{j+1}(x)}{f_j(x)} \right) = \Re \left(\frac{H(L)}{g(I)} \right) > 1 - \frac 1{2^j}\vt
&&(x \in J\mbox{-set for } L,\;L \subseteq I,\;L \in \D_{j+1},\;I \in \D_j), \vt
&& \left|\mbox{arg}\left( \frac{f_{j+1}(x)}{f_j(x)} \right)\right| < \theta_j < \pi \sin \left(\frac{\theta_j}2\right) = \frac \pi{2^ \frac{2j+1}2} < \frac 1{2^{j-2}}.
\end{eqnarray*}
\small{(``arg'' is the \emph{argument} or angle in polar co-ordinates;  $\frac\phi{\sin \phi} \leq \frac \pi 2$ used.)}

\noindent
Hence by (\ref{48(12)}), (\ref{48(15)}),
there is a set $X_j$ with
\begin{equation} \label{48(16)}
V(h:E;X_j) < \frac 1{2^{j-1}} \left(1+V(E)\right),
\end{equation}
and if $k \geq j$, $x \notin X_j$, $|\mbox{arg}(f_{k+1}(x))-\mbox{arg}(f_{k}(x))| < 2^{2-k}$, then, as $j \rightarrow \infty$,
\[
\left|\mbox{arg}(f_{k}(x))- \mbox{arg}(f_{j}(x))\right|< \frac 1{2^{j-3}} \;\;\;\mbox{ and }\;\;\;\mbox{arg}(f_{j}(x)) \rightarrow \theta(x)
\]
as $j \rightarrow \infty$.
If $f(x) \equiv \exp\left(\iota \theta(x)\right)$ then $|f(x)|=1$, $\mbox{arg}(f(x)) = \theta(x)$, and
\begin{eqnarray}
\left| \mbox{arg}(f(x))- \mbox{arg}(f_j(x))\right| & \leq & \frac 1{2^{j-3}},\nonumber \vt
\left|f-f_j\right|^2 &=& \left| \frac f{f_j} -1\right|^2 \;\;=\;\; 4 \sin^2\left( \frac 12 \mbox{arg}\left(\frac f{f_j}\right)\right)\nonumber \vt
&\leq &
\left( \mbox{arg}\left(\frac f{f_j}\right)\right)^2 , \nonumber \vt
\label{48(17)}
\left|f(x)-f_j(x)\right| &<& \frac 1{2^{j-3}}\;\;\;\;\;\;\;\;\;\;\;\;(x \notin X_j).
\end{eqnarray}
From (\ref{48(3)}), (\ref{48(16)}), (\ref{48(17)}), $f_j(x)$ is bounded, integrable in $E$, and convergent almost everywhere to $f$, so that by dominated convergence theorem $f$ is integrable with respect to $V$ and
\[
\left| \int_E f(x)dV-H(E)\right| \leq \frac 1{2^j} + \frac 1{2^{j-3}}V(E) + \frac 1{2^{j-1}} \left(1+V(E)\right).
\]
Hence the theorem for $I^*=E$. For any other $I^*$ there is a $\delta_j^*(P)>0$ such that $\delta_j^*(P) \leq \delta_j(P)$ and such that all $\delta_j^*$-compatible divisions $\D'$ of $E$ satisfy the geometric conditions of Theorem 5. We only have to take $I^*$ in such a division $\D'$. Then except in a set $W$ satisfying (\ref{48(15)}), there is an $f_j^*(x)$ constant in a suitable collection of $J$-sets corresponding to the bricks of $\D'$, such that
\[
\left| f_j^*(x)-f_j(x)\right| < \frac 1{2^{j-2}},\;\;\;\;\;\;\;\;\lim_{j\rightarrow \infty} f_j^*(x) =f(x).
\]
Hence the theorem is true for $I^*$. \nproof

\begin{theorem}\label{t49}
Let $W$ be a real, finitely additive function of bounded variation of elementary sets. If, for all elementary sets $E_1 \subseteq E$,
\[
\overline{W}(E) = \inf\{W(E_1)\},\;\;\;\;\;\;\;\;\underline{W}(E) = \sup\{W(E_1)\},
\]
then
\begin{eqnarray}
0&\leq& \overline{W}(E)\;\;\leq \;\; V(W;E),\nonumber \vt
0&\geq& \underline{W}(E) \;\;\geq\;\; -V(W;E)
\label{49{18}} \vt
W(E) &=& \overline{W}(E) \;\;+\;\;\underline{W}(E). \label{49{19}}
\end{eqnarray}
\end{theorem}
\proof
(\ref{49{18}}) follows from Theorem \ref{47} and the fact that if $E_1$ is empty $W(E_1)=0$. For (\ref{49{19}}), given $\ve>0$ there is an elementary set $E_1 \subseteq E$ with
\[
W(E_1) > \overline{W}(E_1) - \ve,\;\;\;\;\;\;W(E) = W(E_1) + W(E \setminus E_1) > \overline{W}(E) - \ve + \underline{W}(E),
\]
so $W(E) \geq \overline{W}(E) + \underline{W}(E)$. For the opposite inequality we have an elementary set $E_2 \subseteq E$ with
\[
W(E_2) < \underline{W}(E) + \ve,\;\;\;\;\;\;\;W(E) = W(E_2) + W(E \setminus E_2) < \underline{W}(E) + \ve + \overline{W}(E),
\]
so $W(E) \leq \overline{W}(E) + \underline{W}(E)$. \nproof

\begin{theorem}\label{t50}
If an elementary set $E_1 \subseteq E$ has $W(E_1) > \overline{W}(E) + \ve$ then $\underline{W}(E_1)>-\ve$ and $\overline{W}(E\setminus E_1)<\ve$.
\end{theorem}
\proof
By Theorem \ref{49}
\[
 \overline{W}(E) + \underline{W}(E) = W(E_1) > \overline{W}(E)-\ve \geq \overline{W}(E_1)-\ve,
 \]
 hence $\underline{W}(E_1)>-\ve$.
Further $ \overline{W}(E\setminus E_1) + \underline{W}(E\setminus E_1) =$
\[
= W(E \setminus E_1)
 =W(E) - W(E_1) \geq
W(E) -  \overline{W}(E)+\ve < \underline{W}(E)+\ve < \underline{W}(E\setminus E_1)+\ve.
 \]
Therefore $ \overline{W}(E\setminus E_1)< \ve$. \nproof

\begin{theorem}
If $h$, $h_1$ obey the conditions of Theorem \ref{48}, with $h_1$ absolutely continuous in $E$ with respect to $h$, then for every brick $J \subseteq E$ and a point function $f$ independent of $J$ and integrable with respect to $h$, with $|f|$ integrable with respect to $V \equiv V(h;J)$, we have
\[
\int_J h_1 = \int_J f(x) dh.
\]
\end{theorem}
\proof
If $V_1 \equiv V(h_1;J)$ then, using Theorem \ref{Henstock Lemma},
\[
V(h_1;J;X) = V(|h_1|;J;X) = V(|H_1|;J;X) = V(V_1;J;X),
\]
and similarly $V(h;J;X) = V(V;J;X)$. Hence $V_1$ is absolutely continuous in $E$ with respect to $V$.
In Theorem \ref{50} we can put
\[
W = V_1 - bV,\;\;\;\;\;\; V(W;E) \leq V_1(E) + bV(E),
\]
where $b>0$ is a constant. For $\ve>0$ there is an elementary set $E_1 \subseteq E$ with
\[
V_1(E_1) < bV(E_2) +\ve,\;\;\;\;\;\;V_1(E_3) > bV(E_3) - \ve,
\]
for all elementary sets $E_2 \subseteq E_1$, $E_3 \subseteq E \setminus E_1$. For $b = 2^{-j}$, $\ve = 2^{-3-j}$, let $E_1$ be $E_4$. Suppose $E_4, \ldots , E_m$ have been defined. We replace $E$ by $E \setminus \left(E_4 \cup \cdots \cup E_m\right)$, we take $b = (m-2)2^{-j}$, $\ve = 2^{-m-j}$, and then the $E_1$ is called $E_{m+1}$, continuing the induction.
For all elementary sets $E^* \subseteq E$ and $m \geq 4$,
\begin{eqnarray}
\label{51(20)}
\!\!\!\!\!\!&&V(E^*)\frac{m-4}{2^{j}} - \frac 1{2^{m+j-2}} \;\;<\;\;V_1(E^*) \;\;<\;\;V(E^*) \frac{m-3}{2^j} + \frac 1{2^{m+j-1}} , \vt
\label{51(21)}
\!\!\!\!\!\!&& V\left(E\setminus E_m\right) \;\;<\;\; \frac{2^j\left(V_1\left(E \setminus E_m\right) + \frac 1{2^{m+j-1}}\right)}{m-3}\;\;\leq \;\; \frac{2^j\left(V_1(E) +1\right)}{m-3}.
\end{eqnarray}
Taking $f_{jN}(P) = (m-3)2^{-j}$ in the sets $J$ contained in the bricks of $E_m$ $(4 \leq m \leq N)$, and $0$ in the sets $J$ for $E \setminus E_m$, then, for $E^* \subseteq E$,
\begin{eqnarray*}
\int_{E^*} f_{jN}(P) dV &=& \frac{m-3}{2^j} V(E^*), \vt
\left|\int_{E^*} f_{jN}(P) dV - V_1(E^*)\right| & \leq & \frac{V(E^*)}{2^j} + \frac 1{2^{j+m-1}}.
\end{eqnarray*}
For all elementary sets $E^* \subseteq E$
we therefore have
\[
\left|\int_{E^*} f_{jN}(P) dV - V_1(E^*)\right| \;\;\leq\;\; \frac{V(E^*) +1}{2^j} + V_1\left( E^* \cap(E \setminus E_N)\right).
\]
As $f_{jN}(P)$ increases with $N$ and as its integral is bounded, Lebesgue's monotone convergence theorem shows that $f_j(P) \equiv \lim_{N\rightarrow \infty}f_{jN}(P)$ is integrable,
and by the absolute continuity of $V_1$, and (\ref{51(20)}), (\ref{51(21)}),
\begin{eqnarray}
\label{51(22)}
\left|\int_{E^*} f_{j}(P) dV - V_1(E^*)\right| &\leq & \frac{V(E^*) +1}{2^j}, \vt
\label{51(23)}
\lim_{j \rightarrow \infty} \int_{E^*} f_jdV &=& V_1(E^*).
\end{eqnarray}
From (\ref{51(22)}), for all $j>k$ and all elementary sets $E^* \subseteq E$,
\[
\left| \int_{E^*} \left(f_j-f_k\right) dV\right| \;\;\leq\;\; \frac{V(E^*) +1}{2^{k-1}}.
\]
Take a division $\D$ of $E$ and take those $I \in \D$ with $\int_I(f_j-f_k)dV \geq 0$. This gives an elementary set $E^*$. The other $I \in \D$ form another $E^*$. Hence
\begin{eqnarray*}
(\D) \sum \left| \int_J(f_j-f_k)dV\right| & \leq & \frac{V(E)+1}{2^{k-2}}, \vt
V(|f_j-f_k|V;E) &=& V\left(\left| \int_J(f_j-f_k) dV\right|;E\right) \;\;\leq\;\;\frac{V(E)+1}{2^{k-2}},\vt
V \left(\sum_{k=1}^\infty |f_{k+1} - f_k|V;E\right)
& \leq &
\sum_{k=1}^\infty V \left( |f_{k+1} - f_k|V;E\right) \;\;\leq\;\;4(V(E)+1).
\end{eqnarray*}
Hence $\sum_{k=1}^\infty |f_{k+1}-f_k|V$ is finite excepte in a set $X$ with $V(V;E;X) =0$. That is, $V(h;E;X)=0$. Hence $\lim_{j\rightarrow \infty}f_j(P)$ exists $h$-almost everywhere, and Lebesgue's majorised convergence theorem gives $V_1 = \int fdV$. Thus
\[
\int h_1 = \int f_1dV_1,\;\;\;\;\;\;\;\;\;\;\;\;\int h = \int f_2 dV,\;\;\;\;\;\;\;\;\;\;\;\; \int h_1 = \int f_1 f dV = \int \frac{f_1f}{f_2} dh,
\]
completing the proof. \nproof

\section{Abstract Division Space}
We proceed to a general theory for a space $T$ in which there are certain non-empty sets $I$ called \emph{generalized intervals}. We denote the family of these intervals by $\mathcal{T}$.

A set $E \subseteq T$ is called an \emph{elementary set} if $E$ is an interval or a finite union of mutually disjoint intervals. A \emph{division} $\mathcal{D} \equiv \mathcal{D}(I)$ of a set $E \subseteq T$ is the family of an interval $I=E$ or a finite number of mutually disjoint intervals $I$ with union $E$, so that $E$ is necessarily an empty set.

A subfamily $\mathcal{T}_1 \subseteq \mathcal{T}$ \emph{divides} $E$ if a division $\mathcal{D}$ of $E$ exists with $\mathcal{D} \subseteq \mathcal{T}_1$, and then we say that $\mathcal{D}$ comes from $\mathcal{T}_1$.

In the integration theory we often associate one or more points $t$ of $T$ with each interval $I \in \mathcal{T}$ in order to consider elements $f(t)\mu(I)$ used to construct the integral.

To generalize this idea we suppose given a family $\mathcal{T}_t$ of some interval-point pairs $(I,t)$, $I \in \mathcal{T}$, $t \in T$, saying that $t$ is an associated point of $I$. Let $\bS \subseteq \T_t$. Then $\bS$ \emph{divides} $E$ if the sets $\T_2$ of intervals $I$, for all $(I,t) \in \bS$, divides $E$.

If $\D$ is a division of $E$ from $\T_2$ we write $\D_t$ as a general notation for the division $\D$, together with any choice of $t$, one for each $I \in \D$, such that $(I,t) \in \bS$, and we say that $\D_t$ \emph{comes from} $\bS$.

The properties of the integration process obtained by using a certain family $\A$ of the $\bS$, depend greatly on the properties of $\A$.
\begin{example}
$\A$ contains all $\bS$ that are families of those $(I,t)$ that are compatible with some $\delta(t)>0$ in some elementary set $E$. (In other words $\delta$-fine --- P.M.)
\end{example}
So now we have to consider what kinds of general properties of the $\A$ are needed in order that our integration process behaves decently. 

First, we need divisions. Thus we say that $\A$ \emph{divides all elementary sets} if, for each elementary set $E \subseteq T$, there is an $\bS \in \A$ that divides $E$.

$\A$ is \emph{directed} in the sense of divisions if, given $\bS_1, \bS_2 \in \A$, both dividing $E$, there is an $\bS_3 \in \A$ that divides $E$ with $\bS_3 \subseteq \bS_1 \cap \bS_2$. 

The integration will be a Moore-Smith kind of limit with respect to this direction, as we shall see.

To go from a main interval to certain intervals contained in it, we suppose that $\A$ has the \emph{restriction property}. This is explained as follows.

 If $E_1$, $E_2$ are disjoint elementary sets, and if $\bS$ divides $E_1 \cup E_2$, then a \emph{restriction} of $\bS$ to $E_1$ is defined to be a family of some of the $(I,t) \in \bS$ that have $I \subseteq E_1$.

If, for each pair $E_1, E_2$ of disjoint elementary sets and each $\bS \in \A$ that divides $E_1 \cup E_2$, there is a restriction of $\bS$ to $E_1$ that is in $\A$ and divides $E_1$, we say that $\A$ has the restriction property.

If all preceding properties of $\A$ hold, we call $(T, \T, \A)$ a \emph{non-additive division space}.

Such a space has some useful properties, but it does not behave very well when we can integrate over each of two disjoint elementary sets and wish to consider integration over the union of the two sets. To deal easily with this, we say that $\A$ is \emph{additive} if, given disjoint elementary sets $E_j$, and $\bS_j \in \A$ dividing $E_j$ with $I \subseteq E_j$ for all $(I,t) \in \bS_j$ ($j=1,2$), there is an $\bS \in \A$ dividing $E_1 \cup E_2$ with $\bS \subseteq \bS_1 \cup \bS_2$.

If $\A$ is additive, and if all previous properties of $\A$ hold, we say that $(T,\T, \A)$ is a \emph{division space}.

A \emph{partial division} $\Q$ of an elementary set $E$ is a collection of none, or some, or all $I$ in a division $\D$ of $E$ that comes from some $\bS \in \A$. We then denote by $\Q_t$ the set of all $(I,t) \in \D_t$ with $I \in \Q$.

A \emph{partial set} $P$ of $E$ is the union of the $I$ of a partial division $\Q$ of $E$ from $\D$, while $E\setminus P$ is the union of the $I \in \D \setminus \Q$, and $E\setminus P$ is also a partial set, and $P$, $E\setminus P$ are disjoint elementary sets with union $E$.

\begin{theorem}\label{Theorem 1}
Let $(T,\T,\A)$ be a non-additive division space.

\noindent
(a) If $P$ is a partial set of an elementary set $E$, and if $\bS \in \A$ divides $E$, then $\bS$ divides $P$.

\noindent
(b) Let $P$ be a partial set of $E$, and $\Q$ a partial division of $E$
whose corresponding partial set is $E\setminus P$. If $\Q$ is from an $\bS \in \A$ that divides $E$, and if $\bS_0 \in \A$ divides $P$, then there is a division $\D$ of $P$ from $\bS_0$ such that $\Q\cup \D$ is a division of $E$ from $\bS$.
\end{theorem}
\textbf{Proof:} For (a), as $\A$ has the restriction property, as as $P$, $E\setminus P$ are disjoint elementary sets with union $E$, there is a restriction of $\bS$ to $P$ that is in $\A$ and divides $P$. Hence $\bS$ itself divides $P$. For (b) let $\bS_1 \in \A$ be a restriction of $\bS$ to $P$ that divides $P$, and let $\bS_2 \in \A$ divide $P$ with $bS_2 \subseteq \bS_0 \cap \bS_1$ (as $\A$ is directed with respect to divisions. If $\D$ is a division of $P$ from $\bS_2$ then $\D$ is from $\S_0$, and $\Q \cup \D$ is a division of $E$ from $\bS$ as required. \nproof

\begin{theorem}\label{Theorem 2}
Let $(T, \T, \A)$ be a division space\footnote{That is, additive}.

\noindent
(a) If $\D$ is a division of $E$ from an $\bS \in \A$ that divides $E$ there is an $\bS^* \in \A$ that divides $E$ such that $\D'\leq \D$ (i.e.~there is a division of each $I \in \D$ formed of those $J \in \D'$ with $J \subseteq I$; i.e.~$\D'$ is a \emph{refinement} of $\D$) for each division $\D'$ of $E$ from $\bS^*$.

\noindent
(b) Let $P_1, P_2$ be disjoint elementary sets that are partial sets of $E$. Then $P_1 \cup P_2$ is a partial set of $E$.
\end{theorem}
\textbf{Proof:} For (a), as $\A$ has the restriction property there is an $\bS(I) \in \A$ that divides $I$ for each $I \in \D$ while the $(J,t) \in \bS(I)$ has $J \subseteq I$. As $\A$ is additive there is an $S^* \in \A$ that divides $E$, with 
\[
\bS^* \subseteq \bigcup_{I \in \D} \bS(I).
\]
Let $\D'$ be a division of $E$ from $\bS^*$. Then each $J \in \D'$ has $(J,t) \in S^*$ for some $t \in T$, so that $(J,t) \in \bS(I)$ for some $I \in \D$. Then $J \subseteq I$, and $J \cap I'$ is empty for all other $I' \in \D$. Now
\[
\bigcup_{J \in \D} J = E \supseteq I.
\]
Hence $I$ is the union of those $J \in \D'$ with $J \supset I$. Therefore $|d' \leq \D$.

\noindent
For (b), let $P_j$ be formed from a partial division obtained from $\D_j$, a division of $E$, while $\D_j$ comes from an $\bS_j \in \A$ dividing $E$, say. By result (a) there is an $\bS^*_j \in \A$ that divides $E$ such that $\D' \leq \D_j$ for every division $\D'$ of $E$ from $\bS_j^*$ ($j=1,2$). As $\A$ is directed in the sense of divisions there is an $\bS^{**} \in \A$ that divides $E$ and satisfies
\[
\bS^{**} \subseteq \bS_1^* \cap \bS_2^*.
\]
Hence every division $\D'$ of $E$ from $\bS^{**}$ satisfies
\[
\D' \leq \D_1,\;\;\;\;\;\;\D' \leq \D_2.
\]
Taking such a $\D'$, those $I \in \D'$ that lie in the intervals of $\D_1$ used for $P_1$ have union $P_1$, and they are disjoint from those $I \in \D'$ that lie in the intervals of $\D_2$ used for $P_2$, which intervals have union $P_2$. The union of all these intervals from $\D'$ is $P_1 \cup P_2$ which is therefore a partial sum\footnote{Should be ``partial set''?---P.M.}.  \nproof

\section{The space $K$ of values}

Values of integrated functions are to be in an additive topological group $K$, or perhaps some very special $K$ such as the real line or complex plane. We need some definitions:

First $K$ is a \emph{semigroup} if there is a mapping
\[
m:\;\;\;K \times K \rightarrow K
\]
called \emph{multiplication}, and usually written
\[
m(x,y) = x.y,
\]
the product of $x,y$, such that
\begin{enumerate}
\item
$x.(y.z) = (x.y).z$
(all $x,y,z \in K$). $K$ is a \emph{group} if it is a semigroup such that
\item
there is an element $u \in K$ (the \emph{identity} or \emph{unit} element) such that $x.u=u.x =x$ (all $x \in K$) and
\item
to each $x \in K$ there is an element $x^{-1} \in K$ (the \emph{inverse} of $x$) such that
\[
x.x^{-1} = x^{-1}.x = u.
\]
\end{enumerate}
The unit $u$ and inverse $x^{-1}$ can be proved to be uniquely defined, and
\[
(x.y)^{-1} = y^{-1}.x^{-1}, \;\;\;\;\;\;(x^{-1})^{-1}=x.
\]
The group $K$ is \emph{commutative} (or \emph{abelian}) if
\[
x.y=y.x\;\;\mbox{ for all }\;\; x, y \in K,
\]
in which case we usually write $x.y$ as $x+y$, and we call the group \emph{additive}.

There is also a topology attached to $K$.
A family $\G$ of subsets $G$ of $K$ is called a \emph{topology} in $K$ if
\begin{enumerate}
\item[(a)]
The empty set and $K$ are in $\G$;
\item[(b)]
If $\mathcal{H}$ is a subfamily of $\G$, the union of the $G \in \mathcal{H}$ is also a set in $\G$;
\item[(c)]
If $G_1, G_2$ are in $\G$, then so is $G_1 \cap G_2$.
\end{enumerate}
The sets of $\G$ are called \emph{open} sets.
An \emph{open neighbourhood} of $x \in K$ is any $G$ satisfying $x \in G \in \G$. If there are several topologies being used, we say a $\G$-set and a $\G$-neighbourhood, replacing ``open'' by ``$\G$-''.

There is a connection between the group and the topology which is expressed as \emph{continuity}.  Let $(K_j, \G_j)$ ($j=1,2$) be topological spaces, and let $f: K_1 \rightarrow K_2$ be a function. Then $f$ is \emph{continuous at a point} $x \in K_1$ \emph{relative to} $\G_1, \G_2$ if, for $v=f(x)$, given a $\G_2$-neighbourhood $G_2$ of $v$, there is a $\G_1$-neighbourhood $G_1$ of $x$ such that $f(w) \in G_2$ for all $w \in G_1$.

We write this as $f(G_1) \subseteq G_2$.

If true for all $x \in K_1$, let $G_2 \in \G_2$ and let $x \in f^{-1}(G_2)$. (This last means that $f(x) \in G_2$.) By continuity there is a $\G_1$-neighbourhood $G_1$ of $x$ with $f(G_1) \subseteq G_2$, and so
\[
G_1 \subseteq f^{-1}(G_2).
\]
Hence every point of $f^{-1}(G_2)$ is in a $\G_1$-nbd (-neighbourhood) of the point that lies in $f^{-1}(G_2)$. Hence $f^{-1}(G_2)$ is the union of $\G_1$-sets and so is a $\G_1$-set; i.e.,
\[
f^{-1}(G_2) \in \G_1.
\]
Conversely,
 let an arbitrary $x \in K_1$, and let $G_2$ be an arbitrary $\G_2$-nbd of $f(x)$. If $G_3 = f^{-1}(G_2)$ is a $\G_1$-set then it contains $x$, so it is a $G_1$-nbd of $x$ with $f(G_3) \subseteq G_2$. Hence $f$ is continuous at $x$ and so at all points of $K_1$.
 
$K$ is a \emph{Hausdorff space} if, given $x \neq y$, and $x,y in K$, there are disjoint open nbds of $x$ and $y$. 

We can now define the \emph{additive topological group} $K$ as an additive group with a Hausdorff topology (i.e.~$K$ is a Hausdorff space) $\G$ for which $f(x,y):=x-y$ is continuous in $(x,y)$. 

Then $y+G$ is open for all open sets $G$ and all $y \in K$. We denote the identity of the group by $z$ and call it the \emph{zero}.

If $v \in K$, $X \subseteq K$, $Y \subseteq K$, then $v+X$, $X+Y$, $X-Y$ denote the respective sets of points $v+x$, $x+y$, $x-y$ for all $x \in X$, $y \in Y$. 
Thus $X-X$ is not empty; it is the set of $x-y$ for all $x,y \in X$.

A sequence $(u_j)$ in $K$ is \emph{fundamental} (or a \emph{Cauchy sequence}) if, given any open nbd $G$ of the zero $z$, there is an integer $N$ such that $u_j -u_k \in G$ for all $j,k \geq N$.

A sequence $(u_j)$ in $K$ is \emph{convergent} with \emph{limit} $v \in K$ if, given any open nbd $G$ of $z$, there is an integer $N$ such that $u_j-v \in G$ for $j \geq N$. Then the space $K$ is \emph{complete} if every fundamental sequence is convergent.

Generalized Riemann integrals are defined by Moore-Smith limits of a general kind. For each $\bS \in \A$ that divides $E$, let there exist a non-empty set $V(\bS) \subseteq K$ that satisfies
\begin{enumerate}
\item[(d)]
 $V(\bS_1) \subseteq V(\bS)$ when $\bS_1 \subseteq \bS$, and $\bS_1, \bS \in \A$ and divide $E$; i.e.~$V$ is monotone increasing in $\bS$.
\end{enumerate}
Then $V$ is fundamental $(\A,E)$ if, given open nbd $G$ of $z$, there is an $\bS \in \A$ that divides $E$ with
\begin{enumerate}
\item[(e)]
$V(\bS) - V(\bS) \subseteq G$; i.e.~every difference of elements of $V(\bS)$ lies in $G$.
\end{enumerate} 
Also, $V$ is \emph{convergent} $(\A,E)$ with \emph{limit} $v \in K$ if, given any open nbd $G$ of $z$, there is an $\bS \in \A$ dividing $E$ such that
\begin{enumerate}
\item[(f)]
$V(\bS) \subseteq v+G$.
\end{enumerate}
\begin{theorem}\label{Theorem 3}
If $\A$ is directed in the sense of divisions with $V$ monotone increasing in $\bS$, and convergent $(\A,E)$, then $V$ is fundamental $(\A,E)$ and has only one limit.
\end{theorem}
\textbf{Proof:} Let $v,w$ be limits of $V$. By (f), and continuity of $x-y$ at $(z,z)$, given any open nbd $G$ of $z$, there are an open nbd $G_1$ of $z$ with $G_1-G_1 \subseteq G$, and $\bS_j \in \A$ ($j=1,2$) dividing $E$, with
\[
V(\bS_1) \subseteq v+G_1,\;\;\;\;\;\; V(\bS_2) w +G_1.
\]
By (d), as $\D$ is directed in the sense of divisions, there is an $\bS_3 \in \A$ dividing $E$, that satisfies both conditions, so that
\[
v-w = (x-w) - (x-v) \in G_1 - G_1 \subseteq G\;\;\;\mbox{ for all }\;\;x \in V(\bS_1)
\]
i.e.~ $v-w \in G$ for all open nbds $G$ of the zero $z$; i.e.~we cannot have disjoint nbds of $v-w$ and $z$. As $\G$ is a Hausdorff topology,
\[
v-w=z,\;\;\;\;\;\;v=w+z =w,
\]
and $V$ has only one limit. Also
\[
V(\bS_1) -V(\bS_1) =\left\{V(\bS_1)-v\right\} - \left\{V(\bS_1)-v\right\}
\subseteq G_1 - G_1 \subseteq G,
\]
so $V$ is fundamental $(\A,E)$.  \nproof

We assume that $K$ is \emph{complete} $(\A,E)$ for all elementary sets $E$; i.e.~every fundamental $(\A,E)$ $V(\bS)$ that satisfies 4 is convergent $(\A,E)$.

A \emph{local base} $\G_0$ of $\G$ at $z$ is a collection of non-empty open sets such that for all $\G$ satisfying $z \in G\in \G$, there is a $G_0 \in \G$ with $z \in G_0 \subseteq G$.

\begin{theorem}\label{Theorem 4}
If $\G$ has a countable local base $\G_0$ at $z$, and if $K$ is complete, then $K$ is complete $(\A,E)$.
\end{theorem}
\textbf{Proof:}
Let the countable local base at $z$ be put in the sequence $(G_j)$. Then, as addition is continuous, there is an open nbd $G_{1j}$ of $z$ such that \[
G_{1j} +G_{1j} \subseteq G_j,\;\;\;\;\;j=1,2,3, \ldots.
\]
As $G_{1j}$ and $G_1 \cap G_2 \cap \cdots \cap G_{j+1}$ are open nbds of $z$, by definition there is an integer $k = k(j)$ such that\footnote{Note $G_k \subset G_{1j}$.}
\begin{eqnarray}
G_k +G_k \subseteq G_j, \label{1} \vt
G_k \subseteq G_l,\;\;\;\;1 \leq l \leq j+1. \label{2}
\end{eqnarray}
Hence $k(j) >j$, so that by taking a subsequence of $(G_j)$, which by (\ref{2}) is also a local base at $z$, we can assume that, by (\ref{1},\ref{2}), 
\begin{equation}
\label{3}
G_j+G_j \subseteq G_{j-1},\;\;\;\;j=2,3, \ldots.
\end{equation}
Let $\bS_j \in \A$ be such that
\begin{eqnarray}
V(\bS_j)-V(\bS_j) \subseteq G_j \label{4}\vt
V(\bS_{j+1}) \subseteq V(\bS_j). \label{5}
\end{eqnarray}
(5) is possible by (d)  and the fact that $\A$ is directed in the sense of divisions\footnote{If $\bS_{j+1}$ is not contained in $ \bS_j$, we replace $\bS_{j+1}$ by $\bS_{j+1} \cap \bS_j$ or an $\bS_{j+1}^*$ contained in the union [[\emph{should be intersection? - P.M.}]] and a member of $\A$.}
By (3,4,5), $z_k \in V(\bS_j)$ for $k \geq j$ (where $z_j$ is a point of $V(\bS_j)$),
\[
z_j - z_k \in V(\bS_j) - V(\bS_j) \subseteq G_j.
\]
As $(G_j)$ is a local base at $z$, and as $K$ is complete, $(z_j)$ is fundamental, and so is convergent to some point $v$. Hence, by (3,4,5), there is an $l\geq j$ such that $z_l -v \in G_j$,
\[
\begin{array}{lll}
V(\bS_j) -v &=& V(\bS_j) - z_l +z_l -v \vt
& \subseteq & V(\bS_j) - V(\bS_j) + G_j \vt
&\subseteq & G_j + G_j \subset G_{j-1}.
\end{array}
\]
$(G_j)$ being a local base at $z$, $V(\bS_j)$ is convergent $(\A,E)$ to $v$, and $K$ is complete $(\A,E)$.  \nproof

We need some further topological results. First, a \emph{closed set} $F$ is the complement of an open set, while the \emph{closure} $\bar X, = \mbox{Cl}X$ of a set $X \subseteq K$ is the smallest closed set $F \supseteq X$.

Similarly, a \emph{cover} $\mcC$ of a set $C \subseteq K$ is a family of open sets whose union contains $C$. Then $C$ is \emph{compact} if, given an arbitrary cover $\mcC$ of $C$, there is a cover $\mcC_0$ of $C$ consisting of a finite number of open sets (i.e.~a finite cover) so that if $x \in C$, there
are $G \in \mcC$, $G_0 \in \mcC_0$, with $x \in G_0 \subseteq G$. By choice of $G \in \mcC$ with $G \supset G_0$, we can assume that $\mcC_0 \subseteq \mcC$.

The following result follows directly from the definitions.
\begin{theorem}\label{Theorem 5}
(a) As $K$ has a Hausdorff topology, a compact set in $K$ is closed. (b) If $F \subseteq C$, $F$ closed, $C$ compact, then $F$ is compact.
\end{theorem}

Not all closed sets are compact, unless $K$ is compact. A family $\mathcal{F}$ of sets has the \emph{finite intersection property} if, for each finite collection from $\mathcal{F}$, their is a point in their intersection. Taking complements in the definition of compact sets, we see that a set $C$ is compact if and only if every family $\mathcal{F}$ of closed sets in $C$ with the finite intersetion property has a non-empty intersection in $C$.

\begin{theorem}\label{Theorem 6}
If $C$ is compact and contains the sequence $(x_j)$, there is a point $v \in C$ such that if $v \in G \in \G$, then $x_j \in G$ for an infinity of $j$.
\end{theorem}
\begin{theorem}\label{Theorem 7}
Let $f$ be a real-valued continuous function in $K$ and let $C$ be compact. Then $f(C)$ (the set of $f(c)$ for all $c \in C$) is also compact.
\end{theorem}
\begin{theorem}\label{Theorem 8}
If $X \subseteq K$ and $G \in \G$ then $G-X = G- \bar X$.
\end{theorem}
\textbf{Proof:} $G-X \subseteq G-\bar X$. If $y \in G- \bar X$ then $y=g-x$ ($g \in G, x \in \bar X$), and $z \in G-g \in \G$. Thus there is a point $v \in X \cap (G-g+x)$ (where $G-g+x$ is an open nbd of $x \in X$), and $y=h-v$ where 
\[
h=y+v = g-x+v \in g+G-g =G;
\]
i.e.~$y=h-v$ with $h \in G$, $v \in X+y \in G-X$.  \nproof

\vspace{15pt}
\noindent\textbf{Theorem 8A}
If $z in G \in \G$ there is another $G_1 \in \G$ such that $z \in G_1$ and $ \bar{G_1} \subseteq G$. [The topological group is \emph{regular}.]

\vspace{15pt}
\noindent\textbf{Theorem 8B}
$\bar{X+Y} \supseteq \bar X + \bar Y$. [Continuity of $+$.]

\section{The Integral}
We integrate functions $h(I,x)$ of the $I \in \T$, $x \in T$, with values in a set $K$. Sometimes $h(I,x) = f(x) \mu(I)$ where $f$ or $\mu$ or both are real or complex. To integrate $h$ we consider sums
\begin{equation}
\label{6}
(\D_x)\sum h(I,x),
\end{equation}
the notation meaning that we add up the values of $h(I,x)$ for the finite number of $(I,x) \in \D_x$.
Thus we need some operation $x\circ y$ in $K$ (corresponding to addition $\sum$).

Let $\cS(h;\bS;E)$ be the set of sums (\ref{6}) for all divisions $\D_x$ of $E$ from an $\bS \in \A$ that divides $E$, for all re-arrangements of the sums (\ref{6}) by changing the order of the $(I,x)$ in $\D_x$, and for all associations by brackets of values obtained by the operation $\circ$.

If $K$ is a semigroup, the value of (\ref{6}) will be the same whatever association by brackets is used. If $K$ is a communtative semigroup, so that we can use $+$ for $\circ$, the value of (\ref{6}) will only depend on $\D_x$, and not on the order of the $(I,x) \in \D_x$.

To make things easier still, we assume that $K$ is an additive group.

Clearly $\cS(h;\bS;E)$ is a $V(\bS)$ satisfying (d). 

If $K$ is an additive topological group and if this $V(\bS)$ is convergent $(\A,E)$ with limit $H$, we say that $h$ is integrable, and we put 
\[
H=\int_E h,
\]
calling $H$ the \emph{generalized Riemann integral of $h$ in $E$, relative to $\A$}.

If $h(I,x)=f(x)k(I,x)$ we often write the integral as
\[
\int_E h = \int_E f(x)\, dk(I,x).
\]

\begin{theorem}\label{Theorem 9}
Let $C$ be a compact set in $K$ with $\cS(h;\bS;E) \subseteq C$ for some $\bS \in \A$ dividing $E$. If $\A$ is directed in the sense of divisions and if $G$ is an open nbd of the zero, there is an $\bS \in \A$ dividing $E$ with\footnote{$\mbox{Cl}\cS(h;\bS_1;E) \subseteq \left(\bar{\cS}(h;\A;E)+G\right)\cap C$}
\begin{equation}
\label{7}
\cS(h;\bS_1;E) \subseteq \left(\bar{\cS}(h;\A;E)+G\right)\cap C,
\end{equation}
where 
$\bar{\cS}(h;\A;E) = \bigcap \mbox{Cl} \cS(h;\bS;E)$, the intersection being over all $\bS \in \A$ that divide $E$. If the integral $H$ exists then 
$\bar{\cS}(h;\A;E)$ is the set containing the single point $H$.
\end{theorem}
\textbf{Proof:} 
As $\A$ is directed, $\cS(h;\bS;E)$ has the finite intersection property for all $\cS \in \A$ that divide $E$.
As $\cS(h;\bS;E) \subseteq C$, $C$ closed, then $\mbox{Cl}\cS(h;\bS;E) \subseteq C$ and has the finite intersection property, and $\bar{\cS}(h;\A;E)$ is not empty.
 If $z \in G \in \G$ then for all $y \in K$, $y+G \in \G$. Let $G^* = \bar{\cS}(h;\A;E)+G$, $F = C\setminus G^*$. Then $G^*$ is a non-empty union of open sets and so is open, and $F$ is closed. As $F \subseteq C$, $F$ is compact by Theorem 5 (b). If $x \in F$ then $x \notin G^*$, so $x \notin \bar{\cS}(h;\A;E)$ (since $z \in G$). Thus there are an $\bS(x) \in \A$, dividing $E$, and n open nbd $G(x)$ of $x$, with
 $G(x) \cap \cS(h;\bS(x);E)$ empty. As the $G(x)$ cover $F$, a finite number, say $G(x_1), \ldots ,G(x_l)$. Hence $F \cap \mbox{Cl}\cS(h;\bS_1;E)$ is empty, and (\ref{7}) is true. If the integral $H$ exists, then, by its definition, $H \in \mbox{Cl}\cS(h;\cS;E)$. Hence $H\in \bar{\cS}(h;\A;E)$. Also, for each nbd $G_1$ of the zero, and some $\bS\in \A$ dividing $E$,
 \begin{eqnarray*}
 \cS(h;\bS;E) & \subseteq H+G_1, \vt
  \mbox{Cl}\cS(h;\bS;E) & \subseteq H+\bar{G_1}, \vt
   \cS(h;\A;E) & \subseteq H+\bar{G_1}.
 \end{eqnarray*}
Given any nbd $G$ of $z$, we can find another nbd $G_1$ of $z$ such that $\bar{G_1} \subseteq G$ (Theorem 8A). Hence as $\G$ is a Hausdorff topology, $\bar{\cS}(h;\A;E)$ is the set containing the single point $H$.  \nproof

More generally, if the integral does  not exist, and if $\cS(h;\bS;E) \subseteq C$ for some $\bS \in \A$ dividing $E$, then in some sense $\bar{\cS}(h;\A;E)$ characterises the limit points of integration.

A case where we do not have $\cS(h;\bS;E)\subseteq C$ for any $\bS \in \A$ dividing $E$, is one where $\cS(h;\bS;E)$ is always contained in $[n, \infty)$ (i.e.~ $\geq n$) for some integer $n$ depending on $\bS$ and arbitrarily large. E.g.~$T=(-\infty, \infty)$ and $\A$ as the collection of all $\cS$ defined by functions $\delta(x)>0$. Take $g(0)=0$, $g(x) = 1/x$ ($x>0$), $E=[0,1]$, and
\[
\begin{array}{lll}
h([u,v),x) &=& g(v)-g(u),\;\;\;u>0,\vt
h([0,v),x) &=& 0.
\end{array}
\]
If a division $\D$ over $[0,1)$ contains $[0,v)$ and other intervals, the sum
\[
(\D)\sum \Delta g = g(1) - g(v) = -1 + \frac 1v
\]
which tends to $\infty$ as $v$ tends to $0$. Then $\bar{\cS}(h;\A;E)$ is empty (or, by convention, we can say it contains $+\infty$).

A function $h$ of partial sets $P$ of $E$ is \emph{finitely additive} if, for all pairs $P_1,P_2$ of disjoint partial sets of $E$ we have
\[
h(P_1)+h(P_2) = h(P_1 \cup P_2).
\]
\begin{theorem}\label{Theorem 10}
Let $(T, \T,\A)$ be a division space and let $K$ be complete $(\A,E)$. If $\int_Eh$ exists then $\int_P h$ exists as a finitely additive function of the partial sets $P$ of $E$. Also, if $z \in G \in \G$ and
\begin{equation}
\label{8}
\cS(h;\bS;E) \subseteq \int_E h + G
\end{equation}
for some $\bS \in \A$ that divides $E$, then for each partial set $P$ of $E$,
\begin{equation}
\label{9}
\cS(h;\bS;P) \subseteq \int_P h + G-G.
\end{equation}
\end{theorem}
\textbf{Proof:} If $P$ is a partial set of $E$, and if $\bS\in \A$ divides $E$, then, by Theorem 1 (a), $\bS$ divides $P$ and $E \setminus P$. If $s_1,s_2$ are two sums over arbitrary divisions of $P$ from $\bS$, and if $s_3$ is a sum over a division of $E \setminus P$ from $\bS$, then $s_1+s_3$ and $s_2 + s_3$ are sums over two divisions of $E$ from $\bS$. Given $G_1$ with $z \in G_1 \in \G$, we can choose $G$ with $z \in G \in \G$ so that $G-G \subset G_1$ by continuity of $x-y$ at $(z,z)$. Then there is an $\bS\in \A$ dividing $E$ with
\[
\begin{array}{rll}
\cS(h;\bS;E)& \subseteq & \int_E h +G, \vt
s_1-s_2 &=& \left(s_1+s_3 - \int_E h\right)-\left(s_2+s_3 - \int_E h\right)
\in G-G, \vt
\cS(h;\bS;P)-\cS(h;\bS;P)& \subseteq & G-G \subseteq G_1,
\end{array}
\]

and $\cS(h;\bS;P)$ is fundamental $(\A,P)$ and hence convergent $(\A,P)$. Let $s_2 \rightarrow v$. Then by Theorem 8, whether or not $G-G \subseteq G_1$, 
\[
\cS(h;\bS;P)-\int_p h \subseteq  G-\bar G =G-G.
\]
Thus (8) gives (9). To show that the integral is finitely additive over partial sets we use Theorem 2 (b) to prove that partial sets are an additive family, and then we take two disjoint partial sets $P_1,P_2$ of $E$, writing the integral over these sets as $H(P_1), H(P_2)$, respectively. There are $\bS_j \in \A$ dividing $P_j$ such that all sums $s_j$ over $P_j$ from $\bS_j$ satisfy
\[
s_j - H(P_j) \in G_2,\;\;j=1,2, \;\;\mbox{ where we choose }\;\; G_1 + G_2 \subseteq G.
\]
As $\A$ has the restriction property and is additive, there is an $\bS \in \A$ dividing $P_1 \cup P_2$ with $\bS \subseteq \bS_1 \cup \bS_2$. By construction of the $\bS$, every sum $s$ from $\bS$ and over $P_1 \cup P_2$ is $s_1+s_2$ where $s_j$ is a sum over $P_j$ from $\bS$ and so from $\bS_j$ ($j=1,2$), and
\[
s-H(P_1)-H(P_2) = \left(s_1 - H(P_1) \right) + \left(s_2 - H(P_2) \right)
\in G_1 + G_2 \subseteq G,
\]
i.e.~$h$ is integrable over $P_1 \cup P_2$ to the value $H(P_1) + H(P_2)$, and
\[
H(P_1 \cup P_2) = H(P_1) + H(P_2).
\]

\section{The Variation Set}
The sum over a partial division $\Q_x$ of $E$ from an $\bS \in \A$ dividing $E$, the finite sum $(\Q_x)\sum h(I,x)$ is called a \emph{partial sum} over $E$ from $\bS$.
For fixed $E, \bS$ the set of all such partial sums is called the \emph{variation set} $\vs(h;\bS;E)$ of $h$ in $E$ using $\bS$. By construction, if $\Q_x$ is empty the partial sum is taken to have the value $z$, the zero of $K$. Then
\begin{equation}
\label{10}
z \in \vs(h;\bS;E) \supseteq \cS(h;\bS;E).
\end{equation}
We can also define 
\[
\VS(h; \A;E) = \bigcap \mbox{Cl} \vs(h;\bS;E),
\]
the intersection being taken over all $\bS\in \A$ that divide $E$. Clearly 
$\vs(h;\bS;E)$ is a $V(\bS)$ satisfying (d). As $\G$ is Hausdorff, (\ref{10}) shows that $z$ is the only possible limit of $\vs(h;\bS;E)$. When it has this limit we say that $h$ is of \emph{variation zero} in $E$ relative to $\A$.

Further, for each $X \subseteq T$, let
\[
h(X;I,x) := \left\{
\begin{array}{ll}
h(I,x), & x \in X, \vt
z, & x \notin X.
\end{array}
\right.
\]
Then
\begin{eqnarray}
\vs((h;\bS;E;X) &=& \vs(h(X;I,x); \bS;E) \subseteq \vs(h;\bS;E),\label{11.}\vt
\VS((h;\A;E;X) &=& \VS(h(X;I,x); \A;E) \subseteq \VS(h;\A;E),\label{12.}\vt
\vs((h;\bS;E;T) &=& \vs(h; \bS;E),\label{13.}\vt
\VS((h;\A;E;T) &=& \VS(h; \A;E).
\end{eqnarray}
If $\vs(h;\bS;E;X)$ has the limit $z$, we say that $h$ is of \emph{variation zero} in $X$ relative to $E,\A$, or that $X$ is of $h$-\emph{variation zero} in $E$.

\begin{theorem}\label{Theorem 11}
Let $(T, \T, \A)$ be a division space.
\begin{enumerate}
\item[14.]
If $h$ is of variation zero in $E$, its integral is $z$.
\item[15.]
A $h$ of variation zero in $E$ is the same in any partial set of $E$.
\item[16.]
If $h$ is of variation zero in $X$ relative to $E, \A$, then $h$ is of variation zero in any $X1 \subseteq X$.
\item[17.]
If the integral $H(P)$ of $h$ is $z$ for every partial set $P$ of $E$, then $h$ is of variation zero in $E$. Similarly for $h(X;I,x)$.
\item[18.]
If $h_1,h_2$ are of variation zero in $E$, so are $h_1 \pm h_2$.
\item[19.]
If $K$ is the real line, if $h(I,x) \geq 0$ ($x \in X, I\in \T$), if $P$ is a partial set, and if $h$ is integrable over $E$, then $\int_P h \leq \int_E h$.
\item[20.]
If $h \geq 0$ and $\int_E h=0$ then $h$ is of variation zero in $E$.
\end{enumerate}
\end{theorem}
\textbf{Proof:}
(\ref{10}) gives 14; Theorem 1 (a) gives 15; and 16 is clear from the definition.
For 17, let $z \in G_1\in \G$. By continuity and by hypothesis there are a $G \in \G$ with $z \in G$ and
\begin{equation}
\label{21}
G-G \subseteq G_1,\;\;\;\;\;\;\cS(h;\bS;E) \subseteq G.
\end{equation}
Let $\D$ be a division of $E$ from $\bS$, and $\Q$ a partial division, with $P$ as the union of those $I \in \D \setminus \Q$. By Theorem 1(b), if $\bS_1 \in \A$, dividing $P$, is such that
\begin{equation}
\label{22}
\cS(h;\bS_1;P) \subseteq G_1,
\end{equation}
then there is a division $\D_1$ of $P$ from $\bS_1$ such that $\Q\setminus \D_1$ is a division of $E$ from $\bS$. Thus by (\ref{21}),(\ref{22}),
\[
(\Q_x)\sum h(I,x) = (\Q_x \cup \D_{1x})\sum h(I,x) -(\D_{1x})\sum h(I,x)  \subseteq G-G \subseteq G_1,
\]
so $\vs(h;\bS;E) \subseteq G_1$, and $h$ has variation zero in $E$. For 18
\[
\begin{array}{rll}
(\Q_x)\sum \left\{h_1(I,x) \pm h_2(I,x) \right\} &=&
(\Q_x)\sum \left\{h_1(I,x)\right\} \pm \left\{ h_2(I,x) \right\} \vt
& \in & \vs(h_1; \bS;E) \pm \vs(h_2; \bS;E).
\end{array}
\]
For suitable $G_2, \bS$, this is in $G_2\pm G_2 \subseteq G$. In 19, the integrability of $h$ over $P$ follows from Theorems \ref{Theorem 10} an \ref{Theorem 4}, while 20 follows from 19 and 17.  \nproof

Two functions $h,h^*$ of pairs $(I,x)$ are \emph{variationally equivalent} in $E$ relative to $\A$ if $h-h^*$ is of variation zero in $E$ relative to $\A$. If $\Pa$ is the family of partial sets of an elementary set $E$, and if $H(P)$ is finitely additive on $\Pa$ with $H(I)$ variationally equivalent to $h(I,x)$ in $E$, we say that $H(P)$ is the \emph{indefinite variational integral} of $h$ in $E$ relative to $\A$.

\begin{theorem}
\label{Theorem 12}
Let $(T,\T, \A)$ be a division space, and let $K$ be complete $(\A,E)$ for all elementary sets $E$. If $\int_E h$ exists for one of them, then $H(P) = \int_P h$ is the indefinite variational integral of $h$ in $E$ relative to $\A$. Conversely, if the latter integral exists, so does the generalized Riemann integral, and the two are equal.
\end{theorem}
\textbf{Proof:}
Given $z \in G_1 \in \G$, let $z \in G in \G$ with $G-G \subseteq G_1$. For $\bS \in \A$ dividing $E$ let 
\[
\bS(h;\bS;E) \subseteq H(E) +G.
\]
Then, by Theorem \ref{Theorem 10} (\ref{9}),
\begin{equation}
\label{23}
\vs(h;\bS;P) \subseteq H(P) + G-G \subseteq H(P) + G_1
\end{equation}
for each partial set $P$ of $E$. Let $\Q$ be a partial division of $E$ from $\bS$, with union a partial set $P$. Then (\ref{23}) and the finite additivity of $H$ give
\[
(\Q_x)\sum \{h(I,x) - H(I)\} - (\Q_x)\sum h(I,x) - H(P) \in \vs(h;\bS;P) -H(P) \subseteq G_1,
\]
so $\vs(h-H;\bS;E) \subseteq G_1$, $h-H$ is of variation zero, and $H$ is the variational integral.
Conversely, if $H$ is the variational integral, and if $z \in G\in \G$, there is an $\bS \in \A$ dividing $E$ such that for all divisions $\D$ over $E$ from $\bS$,
\[
(\D_x)\sum h(I,x) -H(E) = (\D_x) \sum (h-H) \in \vs(h-H;\bS;E) \subseteq G,
\]
so $\cS(h;\bS;E) \subseteq H(E) +G$, the generalized Riemann integral exists,and is equal to $H(E)$. Similarly for all partial sets of $E$. \nproof

We now consider properties of $\VS(h;\A;E)$.

\begin{theorem}
\label{Theorem 13}
Let $C$ be a compact set with $\vs(h;\bS;E) \subseteq C$ for some $\bS \in \A$ dividing $E$. If $\A$ is directed in the sense of divisions and if $z \in G \in \G$, there is an $\bS_1 \in \A$ dividing $E$, with
\[
\vs(h;\bS_1;E) \subseteq \left(\VS(h;\A;E) +G\right) \cap C.
\]
\end{theorem}
\textbf{Proof:}
Follow the proof of Theorem \ref{Theorem 9}. \nproof

\begin{theorem}
\label{Theorem 14}
Let $(T,\T,\A)$ be a division space. If $E_1,E_2$ are disjoint elementary sets,
\begin{equation}
\label{24}
\VS(h;\A;E_1 \cup E_2) \supseteq \VS(h;\A;E_1)+\VS(h;\A;E_2).
\end{equation}
If, for some compact set $C \subseteq K$ and some $\bS \in \A$ dividing $E_1 \cap E_2$ we have
\begin{equation}
\label{25}
\vs(h;\bS;E_1 \cup E_2) \subseteq C,
\end{equation}
then there is equality in (\ref{24}).
\end{theorem}
\textbf{Proof:}
If $u_j \in \vs(h;\bS_j;E_j)$ for $j=1,2$ then $u=u_1+u_2$ is a sum for a partial division of $E_1 \cup E_2$ from $\bS_1 \cup \bS_2$. Thus if $\bS \in \A$ divides $E_1 \cup E_2$, there is a restriction $\bS_j$ to $E_j$ that lies in $\A$ and divides $E_j$ ($j=1,2$), and then $\bS_1 cup \bS_2 \subseteq bS$. Also the operation $+$ is continuous, so that if $X,Y$ are sets in $K$, Theorem 8B gives $\bar X + \bar Y \subseteq \overline{X+Y}$, and 
\[
\begin{array}{rll}
\VS(h;\A;E_1)+ \VS(h;\A;E_2) & \subseteq &
\mbox{Cl}\vs(h;\bS_1;E_1)+ \mbox{Cl}\vs(h;\bS_2;E_2) \vt
 & \subseteq &
\mbox{Cl}\left\{\vs(h;\bS_1;E)+ \vs(h;\bS_2;E)\right\} \vt
& \subseteq &
\mbox{Cl}\left\{\vs(h;\bS_1\cup\bS_2;E)\right\} \vt
& \subseteq &
\mbox{Cl}\vs(h;\bS;E).
\end{array}
\] 
Taking the intersection for all $\bS in \A$ dividing $E$, we see that (\ref{24}) is true.
By the restriction property there is a restriction $\bS_j \in \A$ of an $\bS in \A$ dividing $E$, such that $\bS_j$ divides $E_j$, that $I \subseteq E_j$ for all $(I,x) \in \bS_j$, ($j=1,2$). By the additive property of $\A$, there is an $\bS_3 \in \A$ that divides $E_1 \cup E_2$ with $\bS_3 \subseteq \bS_1 \cup \bS_2$. Then for each $(I,x) \in \bS_3$, either $I \subseteq E_1$ or $I \subseteq E_2$ so that if $u$ is a sun over a partial division of $E_1 \cup E_2$ from $\bS_3$, we have $u=u_1 + u_2$ where $u_j$ is the value of a sum over a partial division of $E_j$ from $\bS_3$, and so from $\bS_j$ ($j=1,2$), i.le.~$u_j \in \vs(h;\bS_j;E_j)$. Hence
\[
\vs(h;\bS_3; E_1\cup E_2) \subseteq \vs(h;\bS_1; E_1) + \vs(h;\bS_2;  E_2).
\]
But we do not necessarily have $\overline{X+Y} \subseteq \bar X + \bar Y $; e.g.~let $K$ be the real line, $X$ the set of negative integers, and $Y$ the set of points $n + \frac 12 + \frac 1n$ ($n=3,4, \ldots $). Then $X=\bar X$, $Y=\bar Y$, while $\overline{X+Y}$ contains $\frac 12$ which is not in $\bar X+\bar Y$. However, if $z \in X$, $z in Y$, $X+Y \subseteq C$ (compact), let $w \in \overline{X+Y}$ with $w\in G\in \G$. Then $G \cap (X+Y)$ is not empty and contains certain points $x+y$ with $x \in X$, $y \in Y$. Let $X(G)$ be the set of $x \in X$ with $x+y \in G$. Then, as $z \in Y$, we have $X \subseteq C$, $\overline{X(G)} \subseteq C$, and the family of closed sets $\overline{X(G)}$ has the finite intersection property for all $G$ with $w \in G\in \G$. Hence, as $C$ is compact, there is a point $u \in \overline{X(G)}$ for all such $G$. Hence, given $u \in G_1 \in \G$, there is an $x \in X(G) \cap G_1$ and so a $y \in Y$ with $x+y \in G$. Considering the sets $Y(G,G_1)$ of all such $y$ we prove that there is a $v \in \overline{Y(G,G_1)}$ for all $G,G_1$ with $w \in G \in \G$, $w \in G_1 \in \G$. Hence, if $v \in G_2 \in \G$, there is a point $y \in Y(G,G_1) \cap G_2$ and we have $x+y \in G_1 + G_2$, which is contained in an arbitrary nbd $G_3$ of $u+v$ by choice of $G_1, G_2$. Also $x+y \in G$. Thus $G \cap G_3$ is not empty. As $K$ is Hausdorff, 
\[
u+v=w,\;\;\;\;u \in \overline{X(G)},\;\;\;\;v \in \overline{Y(G,G_1)} \subseteq \bar Y,
\]
and we have 
\[
\overline{X+Y} = \bar X + \bar Y .
\]
Now let (\ref{25}) be true. Then we have 
\[
\VS(h;\A;E_1\cup E_2)
\subseteq \mbox{Cl} \vs(h;\bS_3;E_1\cup E_2)
\subseteq
\mbox{Cl}\vs(h;\bS_1;E_1)+ \mbox{Cl}\vs(h;\bS_2;E_2),
\]
and intersections prove equality in (\ref{24}).

\section{The Norm Variation}
It is difficult to get far using variation sets. However if we assume that there is a group norm in $K$ we can define the norm variation and this has easier properties than the variation set. Thus we suppose that there is a function from $K$ to the non-negative real numbers called a \emph{group norm} and written $||u||$, for each $u \in K$, with the properties
\begin{eqnarray}
||u|| =||-u||>0\;\;\mbox{ if }\;\;u \neq z,\;\;\mbox{ while } \;\; ||z||=0; \label{26} &&\vt
||u+v|| \leq ||u||+||v||;\label{27}&& \vt
||u|| \leq a \mbox{ is a closed set.} &&
\label{27A}
\end{eqnarray}
The latter
means that, if every point $ u$ of $ X \subseteq K$  lies in $ ||u|| \leq a$  then every point of $ \bar X$  lies in $ ||u|| \leq a$.

In place of the variation set we define 
\[
V(h;\bS;E) = \sup \left\{ (\D_x)\sum ||h(I,x)||\right\},
\]
the supremum being taken over all divisions $\D$ of an elementary set $E$ from an $\bS\in \A$ that divides $E$. The sup is sometimes $+\infty$.

If $\bS_1 \subseteq \bS_2$, where $\bS_j \in \A$ divide $E$ ($j=1,2$), then every $\D_x$ from $\bS_1$ comes from $\bS_2$, and
\begin{equation}
\label{28}
V(h;\bS_1;E) \leq V(h;bS_2;E).
\end{equation}
Thus we define the \emph{norm variation} of $h$ in $E$ to be
\[
V(h;\A;E) = \inf\left\{V(h;\bS;E) \right\},
\]
the infimum being taken over all $\bS\in \A$ that divide $E$. If the inf is finite we say that $h$ is of \emph{bounded norm variation}. If it is infinite, then $V(h;\bS;E) =+\infty$ for all $\bS \in \A$ that divide $E$. If $V(h;\A;E)=0$, we say that $h$ is of \emph{variation zero}.

Also, we put
\[
V(h;\bS;E;X) = V(h(X;I,x);\bS;E) ,\;\;\;\;\;
V(h;\A;E;X) = V(h(X;I,x);\A;E).
\]
It then follows easily that, for $X_1 \subseteq X \subseteq T$,
\[
V(h;\A;E;X_1) \leq V(h;\A;E;X) \leq V(h;\A;E;T) = V(h;\A;E).
\]
\begin{theorem}
\label{Theorem 15}

\noindent
$\VS(h;\A;E)$ lies in the set of $u$ with $||u|| \leq V(h;\A;E)$.
\begin{equation}
\label{30}
\end{equation}
If $V(h;\A;E)=0$ then $h$ is of variation zero.
\begin{equation}
\label{31}
\end{equation}
If $K$ is the real line $\mathbf{R}$ or the complex plane $\mathbf{Z}$, some point $u \in \VS(h;\A;E)$ has $|u| \geq \frac 14 V(h;\A;E)$ (if finite), or $\vs(h:\bS;E)$ is unbounded for all $\bS \in \A$ dividing $E$ if $V(h;\A;E) = +\infty$.
\begin{equation}
\label{32}
\end{equation}
In (\ref{32}), if $h$ is of variation zero then $V(h;\A;E)=0$.
\begin{equation}
\label{33}
\end{equation}
In (\ref{32}), if $h(I,x)$ is integrable in $P$ to $H(P)$ for each partial set $P$ of $E$, then $|h(I,x) -H(I)|$ is of variation zero.
\begin{equation}
\label{34}
\end{equation}
\end{theorem}
\textbf{Proof:}
Let $\bS\in\A$ divide $E$, and let $\Q$ be a partial division of a division $\D$ of $E$ from $\bS$. Then, by (\ref{26}), \ref{27}),
\[
||(\Q_x)\sum h(I,x)|| \leq (\Q_x) \sum ||h(I,x)|| \leq (\D_x)\sum ||h(I,x)|| \leq V(h;\bS;E).
\]
By (\ref{27A}) we have (\ref{30}), from which (\ref{31}) follows by (\ref{26}). For (\ref{32}) let $\vs(h;\bS;E)$ be a bounded set, for some $\bS\in\A$ dividing $E$. Then for each $\bS_1 \in \A$ dividing $E$ with $\bS_1 \subseteq \bS$, we put 
\[
N(\bS_1) = \sup |u|,
\]
the supremum being taken for all $u \in \vs(h;\bS_1;E)$. Then
\[
|(\Q_x) \sum h(I,x)| \leq N(\bS_1) < \infty.
\]
By the usual argument, $(\D_x) \sum |h(I,x)| \leq 4N(\bS_1)$. Hence
\[
V(h;\A;E) \leq V(h;\bS_1;E) \leq 4N(\bS_1),\;\;\;\;\;\;
N(\bS_1) \geq \frac 14 V(h;\A;E).
\]
Hence $V(h;\A;E)$ is finite, and there is a point $v \in \mbox{Cl}\vs(h;\bS_1;E)$ with $|v| \geq \frac 14 V(h;\A;E)$. Let $B$ be the compact set
\[
\frac 14 V(h;\A;E) \leq |k| \leq N(\bS_1),
\]
an annulus or the boundary of a circle. Then Cl$\vs(h;\bS_1;E)\cap B$ has the finite intersection property for all $bS_1 \in \A$ that divide $E$ (---since $\A$ is directed, a finite number of $\bS_1$'s have an $\bS$ inside them; $|v| \geq \frac 14 V(h;\A;E)$). Hence there is a point $u \in \VS(h;\A;E) \cap B$, proving the first part of (\ref{32}). For the second part, if $V(h;\A;E)=+\infty$, we cannot have a finite $N(\bS)$, so that $\vs(h;\bS;E)$ is an unbounded set, for all $\bS \in \A$ dividing $E$. If, further, $h$ has variation zero, then $N(\bS) \geq 0$ is as small as we please by choice of $\bS \in \A$ dividing $E$. Hence $V(h;\A;E) =0$, giving (\ref{33}). Then (\ref{31}), (\ref{33}) and Theorem \ref{Theorem 12} give (\ref{34}).  \nproof

Another definition of norm variation is obtained as follows. First, for $K$ the real line, a $H(I)$, independent of $x$, is \emph{finitely subadditive} on $\T_1 \subseteq \T$ if, for each $I \in \T_1$ and each division $\D$ of $I$ using intervals $J \in \T_1$ alone, we have $(\D) \sum H(J) \geq H(I)$.

A $H$ is \emph{finitely superadditive} if $-H$ is finitely subadditive. Then we can say that $h(I,x)$ is of \emph{bounded norm variation} in $E$ if there are an $\bS \in \A$ dividing $E$, and a non-negative finitely superadditive function $\chi(I)$ with
\[
||h(I,x)||\leq \chi(I) \;\;\;\mbox{ for all }\;\;(I,x) \in \bS.
\]
If $V(h;\bS;J)$ is finitely superadditive, then for fixed $\bS,J$, it is the least possible $\chi(J)$. Thus for fixed $\bS,J$ we can take the infimum of all such $\chi(J)$, and then we can take the infimum of the result, for fixed $J$ and all $\bS\in \A$ dividing $J$, and we have
\[
\inf \inf \chi(J) = \inf V(h;\bS;J) = V(h;\A;J).
\]
Thus we need the following.
\begin{theorem}
\label{Theorem 16}
Let $(T, \T,\A)$ be a division space.

\noindent
If $\bS \in \A$ divides $E$, then $V(h;\bS;I)$ is a finitely superadditive function of the partial sets of $E$ that are intervals $I$.
\vspace{-10pt}
\begin{equation}
\label{35}
\end{equation}
If $E_1,E_2$ are two disjoint elementary sets and if $\bS_j \in \A$ divides $E_j$ with $I \subseteq E_j$ for all $(I,x) \in \bS_j$ ($j=1,2)$, then
\begin{equation}
\label{36}
V(h;\bS_1;E_1) + V(h;\bS_2;E_2) = V(h;\bS_1\cup\bS_2;E_1 \cup E_2).
\end{equation}
\vspace{-10pt}


If $E_1,E_2$ are two disjoint elementary sets then
\begin{equation}
\label{37}
V(h;\A;E_1) + V(h;\A;E_2) = V(h;\A;E_1 \cup E_2).
\end{equation}

\end{theorem}
\textbf{Proof:} For (\ref{35}) let partial intervals $I_1, \ldots ,I_n$ form a division of $I$ and let $\D_j$ be a division of $I_j$ from $\bS$ ($1\leq j \leq n$). These divisions exist by Theorem \ref{Theorem 1} (a). Then $\D=\bigcup_{j=1}^n \D_j$ is a division of $I$ from $\bS$ and
\[
\sum_{j=1}^n (\D_j) \sum ||h(J,x)|| = (\D)\sum ||h(J,x)|| \leq V(h:\bS;I).
\]
For suitable $\D_j$ we can make the sums $(\D_j)\sum$ on the left tend to the corresponding $V(h;\bS;I_j)$ and (\ref{35}) follows. Similarly in (\ref{36}) we have
\begin{equation}
\label{38}
V(h;\bS_1;E_1) + V(h;\bS_2;E_2) \leq V(h;\bS_1\cup\bS_2;E_1\cup E_2).
\end{equation}
Now if $(I,x) \in \bS_1\cup \bS_2$ then either $(I,x) \in \bS_1$, $I \subseteq E_1$, disjoint from $E_2$, or $(I,x) \in \bS_2$, $I\subseteq E_2$, disjoint from $E_1$, and each division $\D$ of $E_1 \cup E_2$ from $\bS_1 \cup \bS_2$ can be separated into a division $\D_j$ of $E_j$ from $\bS_j$ ($j=1,2$). Thus
\[
(\D)\sum ||h(I,x)|| = (\D_1)\sum + (\D_2)\sum \leq V(h;\bS_1;E) + V(h,\bS_2;E_2).
\]
By taking suitable $\D$ we have the opposite inequality to (\ref{38}), and so equality. Let $\bS \in \A$ divide $E_1 \cup E_2$ and let $\bS_j \in \A$ be a restriction to $E_j$ that divides $E_j$, ($j=1,2$). Then $\bS_1 \cup \bS_2 \subseteq \bS$, and, by (\ref{36}),
\[
\begin{array}{lll}
V(h;\A;E_1)+V(h;\A;E_2) &\leq & V(h;\bS_1;E_1)+V(h;\bS_2;E_2) \vt
&=& V(h;\bS_1\cup \bS_2;E_1 \cup E_2) \leq V(h;\bS;E_1 \cup E_2) .
\end{array}
\]
Taking the infimum as $\bS$ varies,
\begin{equation}
\label{39}
V(h;\A;E_1) + V(h;\A;E_2)  \leq V(h;\A;E_1\cup E_2) .
\end{equation}
Also, by (\ref{36}) and by choice of $\bS_1, \bS_2$ and (by additivity) an $\bS_3\in \A$ dividing $E_1 \cup E_2$ with $\bS_3 \subseteq \bS_1 \cup \bS_2$, we have
\[
\begin{array}{lll}
V(h;\A;E_1 \cup E_2) & \leq & V(h;\bS_3;E_1 \cup E_2)\vt
&\leq & V(h;\bS_1 \cup \bS_2;E_1 \cup E_2) \vt
&=& V(h;\bS_1;E_1 ) + V(h;\bS_2; E_2) \vt
&\leq & V(h;\A;E_1 ) + V(h;\A; E_2) + \varepsilon.
\end{array}
\]
This gives the reverse inequality to (\ref{39}) (for additive $\A$), and so (\ref{37}).   \nproof

\begin{theorem}
\label{Theorem 17}
If $(T,\T,\A)$ is a division space, if $V(h;\A;E) < \infty$, if $\varepsilon>0$, and if $\bS \in \A$ dividing $E$ satisfies
\begin{equation}
\label{40}
V(h;\bS;E) < V(h:\A;E) + \ve,
\end{equation}
then, for all partial sets $P$ of $E$,
\begin{equation}
\label{41}
V(h;\bS;P) < V(h;\A;P) +\ve.
\end{equation}
\end{theorem}
\textbf{Proof:}
As $\bS$ contains the union of its restrictions to $P, E\setminus P$, then Theorem \ref{Theorem 16} (\ref{36}), (\ref{37}) imply that
\[
\begin{array}{rll}
V(h;\bS;P) & \leq & V(h;\bS;E) - V(h;\bS;E \setminus P) \vt
&<& V(h;\A;E) + \ve - V(h;\A;E\setminus P) \vt
&=& V(h;\A;P) + \ve,
\end{array}
\]
giving (\ref{41}).  \nproof

This is a crucial result for monotone convergence theorem.

It is doubtful whether $\vs(h:\bS;E)$ satisfies a similar result in which $\ve$ is replaced by an open set. To replace $V(h;\bS;E)<\infty$ we could use a compact set $C$ and an $\bS \in \A$ dividing $E$ with $\vs(h;\bS;E) \subseteq C$. But the proof with $V$ fails with $\vs$; for let $K=T$ be the real numbers, with
\[
h(u,v) = \left\{
\begin{array}{ll}
1,& u<0<v<1,\vt
v-u, & 1 \leq u<v \leq 2,\vt
0 & \mbox{otherwise.}
\end{array}
\right.
\]
Then $\VS(h;\A;[-1,1])$ consists of the points $0$ and $1$, 
\[\begin{array}{l}
\VS(h;\A;[-1,1])=\{0,1\},\vt \VS(h;\A;[1,2])=[0,1],\vt
\VS(h;\A;[-1,2])=[0,2].
\end{array}
\]
 But from the second and third $\VS$ we cannot reconstruct the first.
 
 \noindent
\textbf{Counterexample by P.~Muldowney:}
$h(u,v):=v-u$ if $2 \leq u<v\leq 3$; $=u$ if $u\leq 1<v<2$; $=0$ otherwise. Let $\bS$ be defined by a gauge $1>\delta(x)>0$ for $x\neq 1$, with $0<\delta(1)=\alpha$. Then
\[
\begin{array}{l}
\vs(h;\bS;[2,3]) =[0,1],\;\;\;\;\;\;\vs(h;\bS;[0,3])=[0,2],
\vt
\VS(h;\A;[2,3]) =[0,1],\;\;\;\;\;\;\VS(h;\A;[0,3])=[0,2],
\vt
\vs(h;\bS;[0,2]) =\{0\}\cup (1-\alpha,1],\vt
\mbox{Cl}\vs(h;\bS;[0,2]) =\{0\}\cup [1-\alpha,1],\vt
\VS(h;\A;[0,2]) =\bigcap \cdots =\{0\}\cup \{1\}.
\end{array}
\]
Take $G=(-\ve,\ve)$ where $0<\ve<\alpha$.  \nproof

\begin{theorem}
\label{Theorem 18}
Let $T,\T,\A)$ be a division space. Let $h(I)$ be real and finitely subadditive over intervals $I$ that are partial sets of an elementary set $E$. If there is an $\bS\in\A$ dividing $E$ such that $\cS(h;\bS;E)$ is bounded above with supremum $s$, then $h$ is integrable in $E$ with integral equal to $s$. If $\cS(h;\bS;E)$ is unbounded above for all $\bS \in \A$ dividing $E$, then, for each integer $n$ there is an $\bS_n\in\A$ dividing $E$ such that 
\[
\cS(h;\bS_n;E) \subseteq [n,\infty).
\]
\end{theorem}
\textbf{Proof:}
By Theorem \ref{Theorem 2} (a), given a division $\D$ of $E$ from an $\bS\in \A$ dividing $E$, there is an $\bS^*\in \A$ dividing $E$ such that $\D' \leq \D$ for each division $\D'$
of $E$ from $\bS^*$. As $h(I)$ is finitely subadditive we obtain
\[
(\D')\sum h(I) \geq (\D)\sum h(I) , = d, \mbox{ say}.
\]
It follows that $\cS(h, \bS^*,E)$ is bounded below by $d$. In the first case we can say that $d>s-\ve$ when $\ve>0$ is given. Then as $\bS^* \subseteq \bS$, $\cS(h;\bS^*;E)$ lies between $s$ and $s-\ve$. Hence $h$ is integrable to $s$. In the second case we can choose $\D=\D_n$ so that $d>n$, and $\bS_n$ can then be taken as the corresponding $S^*$.

\begin{theorem}
\label{Theorem 19}
Let $T,\T.\A)$ be a division space. If $h(I,x)$ is real or complex, integrable with integral $H$ in $E$, and of bounded variation, then $|h(I,x)|$ and $|H(I)|$ are integrable to $V(P)=V(h;\A;P)$, and
\[
V(H;\A;P) = V(P),\;\;\;\;\;\;|H(P) \leq V(P),
\]
for each partial set $P$ of $E$.
\end{theorem}
\textbf{Proof:}
For divisions $\D_x$ of partial sets $P$ of $E$,
\begin{equation}
\label{42}
\left|(\D_x)\sum |h(I,x)|-(\D)\sum |H(I)|\right|
\leq (\D_x)\sum \left||h| -|H| \right| \leq (\D_x)\sum |h-H|.
\end{equation} 
This is as small as we please for all $\D_x $ from $\bS$ by choice of $\bS\in\A$ dividing $E$ by Theorems \ref{Theorem 1} (a), \ref{Theorem 15} (\ref{34}). Hence, from (\ref{42}),
\begin{equation}
\label{43}
V(H;\A;P) = V(h;\A;P) =:V(P).
\end{equation}
By Theorem \ref{Theorem 16}. $H$ is finitely additive, so that $|H|$ is finitely subadditive. Then (\ref{43}) and Theorem \ref{Theorem 18} give the integrability of $|H|$, $|h|$ to $V(P)$, and $|H(P)| \leq V(P)$ follows from the finite subadditivity of $|H|$.

\section{The integrability of functions of interval-point functions}
In this section we take $K$ as the real line, and $r(x_1, \ldots , x_n)$ a function of $n$ real variables.  We proceed to the integrability of $r(h_1(I,x), \ldots , h_n(I,x))$, given the integrability of $h_j(I,x)$ ($1 \leq j \leq n$). 

Naturally we have to assume that $r$ is fairly smooth, and the first condition we have to impose on $r$ is that, for constants $A_j$, ($1 \leq j \leq n$), we have
\begin{equation}
\label{44}
\left|r(y_1, \ldots , y_n)-r(x_1, \ldots , x_n)\right| \leq A_1|y_1-x_1| + \cdots + A_n|y_n-x_n|;
\end{equation}
e.g.~this holds when
\[
\left| \frac{\partial r}{\partial x_j}\right| \leq A_j,\;\;\;\;1\leq j \leq n.
\]
For the second condition, weaker than the first,, we put
\[
r_1(x_1, \ldots , x_n; \ve) := \sup \left|r(y_1, \ldots, y_n) - r(x_1, \ldots ,x_n)\right|,
\]
the supremum being taken over all $(y_1, \ldots ,y_n)$ satisfying $|y_j -x_j| \leq \ve$, ($1 \leq j \leq n$). Here, $\ve >o$; and $(x_1, \ldots ,x_n)$ is fixed. For arbitrarily large $m$, varying $\ve_k>0$ $(1 \leq k \leq m$), and with the following expression fixed:
\begin{equation}
\label{45}
\sum_{k=1}^m x_{jk} \;\;\;\;1 \leq j \leq n,
\end{equation}
we suppose that
\begin{equation}\label{46}
\sum_{k=1}^m r_1(x_{1k}, \ldots , x_{nk}; \ve_k) \rightarrow 0\;\;\mbox{ as }\;\;\sum_{k=1}^m \ve_k \rightarrow 0
\end{equation}
If (\ref{44}) is true, then so is (\ref{46}), since
\[
\begin{array}{rll}
r_1(x_1, \ldots , x_n;\ve) & \leq & (A_1 + \cdots + A_n)\ve, \vt
\sum_{k=1}^m r_1(x_{1k}, \ldots , x_{nk}; \ve) & \leq & \sum_{j=1}^n A_j \sum_{k=1}^m \ve_k \rightarrow 0.
\end{array}
\]
\begin{theorem}
\label{Theorem 20}
Let $\A$ be directed in the sense of divisions, let $\mathbf{R}, \mathbf{R}^+$ be the real line and the line of non-negative numbers, with $\mathbf{R}_n, \mathbf{R}_n^+$ their respective $n$-fold Cartesian products, let $U=\mathbf{R}_n$ or $\mathbf{R}_n^+$, and let
\[
r(x_1, \ldots , x_n): U \mapsto \mathbf{R}.
\]
Let $h_j(I,x)$ be variationally equivalent to $k_j(I,x)$ ($1 \leq j \leq n$). with values in $\mathbf{R}$ or $\mathbf{R}^+$, in an elementary set $E$. Then $r(h_1, \ldots , h_n)-r(k_1, \ldots ,k_n)$ has variation zero in $E$ if \textbf{either} $r$ satisfies (\ref{44}) in $U$ \textbf{or} $r$ satisfies (\ref{46})  in $U$ with the $h_j(I,x)$ ($1\leq j \leq n$) integrable in $E$, and with $(T,\T,\A)$ a division space.
\end{theorem}
\textbf{Proof:}
Given $\ve>0$, there are $\bS_j \in \A$ dividing $E$ such that
\begin{equation}
\label{47}
V(h_j-k_j; \bS_j;E) < \ve,\;\;\;\;1\leq j\leq n.
\end{equation}
As $\A$ is directed in the sense of divisions, there is an $\bS\in \A$ dividing $E$ with
\[
\bS \subseteq \bS_1 \cap \bS_2\cap \cdots \cap\bS_n.
\]
Then as $\bS\subseteq \bS_j$, (\ref{47}) is true with $bS$ replacing $\bS_j$. If (\ref{44}) holds then for each division $\D_x$ of $E$ from $\bS$,
\[
(\D_x) \sum |r(h_1, \ldots , h_n) - r(k_1, \ldots , k_n)| \leq \sum_{j=1}^n A_j (\D_x) \sum |h_j-k_j| \leq \sum_{j=1}^n A_j \ve,
\]
and $r(h_1, \ldots , h_n) - r(k_1, \ldots , k_n)$ has variation zero. If (\ref{46}) holds, with $h_j$ integrable to $H_j$ ($1 \leq j \leq n$) then $k_j$ is also integrable to $H_j$ ($1\leq j\leq n$), since
\[
V(k_j - H_j; \A;E) \leq V(k_j-h_j;\A;E) + V(h_j -H_j;\A;E) =0.
\]
Since 
\[
\begin{array}{rll}
r(h_1,\ldots h_n)-r(k_1,\ldots k_n)&=&
\left\{r(h_1,\ldots h_n)-r(H_1,\ldots H_n)\right\} \;\;- \vt
&& \;\;\;\;-\;\; \left\{r(k_1,\ldots, k_n)-r(H_1,\ldots H_n\right\}),
\end{array}
\]
it is clear that we can replace $k_j$ by $H_j$ ($1\leq j \leq n$), and prove the first $\{ \cdots \}$ of variation zero. The second $\{ \cdots \}$ will be similar. Let $\D_x$ be a division of $E$ from $\bS$, consisting of $(I_l, v_l)$, ($1 \leq l \leq m$). Then the value of (\ref{45}) is $H_j(E)$, fixed, on taking $x_{jl} = H_j(I_l)$. With
\[
y_{jl} = h_j(I_l, u_l),\;\;\;\;
\ve_l=\sum_{l=1}^n V(h_j-H_j;\bS; I_l),
\]
we get
$
\sum_{l=1}^m \left| r(h_1(I_l, u_l) , \ldots , h_n(I_l, u_l)) - r(H_1(I_l) , \ldots , H_n(I_l)) \right| \leq $
\[\leq 
\sum_{l=1}^n  r_1(H_1(I_l) , \ldots , H_n(I_l);\ve_l),
\]
and
\[
\sum_{l=1}^m \ve_l \leq  \sum_{j=1}^n V(h_j - H_j; \bS; E) < n\ve,
\]
using Theorem \ref{Theorem 16} (\ref{35}) for the finite superadditivity of $V$. Hence by choice of $\bS$,
\[
V(r(h_1, \dots , h_n) -r(H_1, \ldots , H_n); \bS;E) \rightarrow 0,
\]
completing the proof.   \nproof 
\begin{theorem}
\label{Theorem 21}
Let $(T,\T,\A)$ be a division space. Let $R$ satisfy (\ref{46}) in $U$ with the $h_j(I,x)$ integrable in $E$ to $H_j$ ($1 \leq n$) and with
\begin{equation}
\label{48}
r(x_1+y_1, \ldots ,x_n+y_n) \leq r(x_1, \ldots ,x_n) 
+ r(y_1, \ldots ,y_n) 
\end{equation}
for all $(x_1, \ldots , x_n), (y_1, \ldots , y_n$ in $U$. Then $r(h_1, \ldots , h_n)$ is integrable with integral equal to the integral of $r(H_1, \ldots , H_n)$ if and only if, for some $\bS \in \A$ dividing $E$, and for some compact set $C$,
\begin{equation}
\label{49}
\cS(r(h_1, \ldots , h_n);\bS;E) \subseteq C.
\end{equation}
\end{theorem}
\textbf{Proof:}
By (\ref{48}) the function $c(P)= r(H_1(P), \ldots H_n(P))$ is a finitely subadditive function of partial sets $P$ of $E$. For if some $\bS \in \A$ divides $E$ then by Theorem \ref{Theorem 1} (a), $\bS$ divides $P$, and if $\D$ is a division of $P$ from $\bS$,
\[
\begin{array}{rll}
(\D) \sum H_j(I) &=& H_j(P) ,\;\;\;\;1 \leq j \leq n,\vt
(\D)\sum r(H_1, \ldots H_n) & \geq & r\left((\D) \sum H_1, \ldots , (\D) \sum H_n\right) \vt
&=& r(H_1(P), \ldots , H_n(P)), \vt
(\D) \sum c(I) & \geq & c(P).
\end{array}
\]
By Theorem \ref{Theorem 18}, $c(I)$ is bounded for some $\bS \in \A$ dividing $E$, and all divisions $\D_1$ of $E$ from $\bS$.
By Theorem \ref{Theorem 20}, $c(I)$ is variationally equivalent to $r(h_1(I,x), \ldots , h_n(I,x))$, which is integrable if and only if $c(I)$ is integrable. Also from this, as $\cS\subseteq \vs$, $\cS(c-r(h_1, \ldots ,h_n);\bS_1;E)$ is bounded for some $\bS_1 \in \A$ divding $E$. Hence $\cS(c;\bS_1 \cap \bS;E)$ is bounded if and only if $\cS(r(h_1, \ldots ,h_n);\bS_1 \cap \bS;E)$  is bounded. Hence $r(h_1, \ldots , h_n)$ is integrable if and only if (\ref{49}) is true.   \nproof

\begin{theorem}
\label{Theorem 22}
Let $h_j(I,x)$ be integrable to $H_j$ in $E$ ($1 \leq j\leq n$). Then \[\max\{h_1(I,x), \ldots ,h_n(I,x)\}\] is integrable to the integrals of $\max\{H_1, \ldots ,H_n\}$ if and only if, for some compact set $C$, and some $\bS \in \A$ dividing $E$,
\begin{equation}
\label{50}
\cS\left(\max(h_1, \ldots , h_n);\bS;E\right) \subseteq C.
\end{equation}
\end{theorem}
\textbf{Proof:}
We show that $r(x_1, \ldots ,x_n) :=\max(x_1, \ldots ,x_n)$ satisfies (\ref{44}) and (\ref{48}). Then Theorem \ref{Theorem 21} completes the proof. For (\ref{44}),
\[
\begin{array}{rll}
x_j &=& (x_j-y_j) + y_j \;\;\leq \;\;|x_j-y_j| +y_j \vt
&\leq & \sum_{j=1}^n |x_j-y_j| + \max(y_1, \ldots ,y_n),
\vt
\max(x_1, \ldots , x_n) &\leq & \sum_{j=1}^n |x_j-y_j| + \max (Y_1, \ldots , y_n),
\end{array}
\]
and interchanging $(x_1, \ldots ,x_n)$ and $(y_1, \ldots ,y_n)$ gives (\ref{44}) with $A_1 = \cdots =A_n =1$. We have (\ref{48}) from 
\[
x_j+y_j \leq \max(x_1, \ldots ,x_n) + \max(y_1, \ldots , y_n).
\]
\textbf{Corollary:}
There is a similar result with $\min$.

Note that (\ref{t50}) is true if
\begin{equation}
\label{51}
(\D_x)\sum h_{j(I,x)} (I,x) \subseteq C
\end{equation}
for every $\D_x$ of $E$ from $\bS$, for all choices of integers $j(I,x)$ in the range from $1$ to $n$. Let $h(I,x)$ be integrable to $H(E)$ in $E$. Then, for some $\bS\in\A$ dividing $E$, we have $\cS(h;\bS;E)$ bounded, and in fact lying in $(H-\ve, H+\ve)$ for some $\ve>0$, and $\bS$ depending on $\ve$. If $h_1 \leq h$ for $1\leq j\leq n$, then 
\[
h_1 \leq \max(h_1, \ldots ,h_n) \leq h,
\]
and so the $\cS$ in (\ref{50}) is in a finite interval from $H_1-\ve$ to $H+\ve$, for suitable $\bS$. If instead we have $h_j \geq h$ for $1 \leq j \leq n$, then $h_j -h \geq 0$, and
\[
\begin{array}{rll}
h_j-h & \leq & \sum_{k=1}^n ( h_k -h) \;\;=\;\; \sum_{k=1}^n h_k - nh, \vt
\max(h_1, \ldots ,h_n) &\leq & \sum_{k=1}^n h_k - (n-1) h ,
\end{array}
\]
and we have (\ref{50}). 
 \nproof
 
 \noindent
 \textbf{Note:} $h$ integrable, $|h|$ not; $\max(h, -h) = |h|$.
 
 \begin{theorem}
 \label{Theorem 50}
 Let $h_j(I,x)\geq 0$ be integrable in $E$ to $H_j$ ($j=1,2$). If $t$ is fixed in $0<t<1$, then $h_1^t h_2^{1-t}$ is integrable in $E$ to the integral of $H_1^t H_2^{1-t}$.
 \end{theorem}
\textbf{Proof:}
Let $r(x_1,x_2) = x_1^tx_2^{1-t}$ ($x_1 \geq0, x_2 \geq 0$). Then (\ref{48}) (subadditivity) is true for $-r$. For let \[
t=\frac 1p,\;\;\;\;1-t - 1- \frac 1p = \frac 1q.
\]
Then H\"{o}lder's inequality gives 
\[
u_1u_2 + v_1v_2 \leq \left(u_1^p + v_1^p \right)^\frac 1p \left(u_2^q + v_2^q\right)^\frac 1q.
\]
With $u_1^p = x_1$, $u_2^q = x_2$, $v_1^p  =y_1$, $v_2^q=y_2$,
\[
x_1^t x_2^{1-t} + y_1^t y_2^{1-t} \leq \left(x_1 + y_1\right)^t\left(x_2 +y_2\right)^{1-t}.
\]
To prove (\ref{46}) we begin with
\begin{equation}
\label{52}
f(x) = x^t + y^t -(x+y)^t \geq 0\;\;\;\;\;\;(x\geq 0, y\geq 0, 0<t<1).
\end{equation}
This follows from $f(0)=0$ and $f'(x) = tx^{t-1} - t(x+y)^{t-1} \geq 0$. (since $t-1<0$); $f(0) \geq 0$. From (\ref{52}) we obtain
\[
(x_1+y_1)^t(x_2+y_2)^{1-t} \leq \left(x_1^t + |y_1|^t\right)\left(x_2^{1-t} + |y_2|^{1-t}\right),
\]
\begin{equation}
\label{53}
(x_1+y_1)^t(x_2+y_2)^{1-t} - x_1^t x_2^{1-t}
 \leq 
x_1^t |y_2|^{1-t} + |y_1|^t x_2^{1-t} + |y_1|^t |y_2|^{1-t},
\end{equation}
($ x_j \geq 0, x_j + y_j \geq 0, j=1,2$.) Next we prove that
\begin{equation}
\label{54}
x_1^t x_2^{1-t} - (x_1+y_1)^t(x_2+y_2)^{1-t} 
 \leq 
x_1^t |y_2|^{1-t} + |y_1|^t x_2^{1-t} + |y_1|^t |y_2|^{1-t}.
\end{equation} 
This is clearly true when $y_j \geq 0$ ($j=1,2$). When 
$y_1<0$, $y_2 \geq O$, we use (\ref{52}):
\[
\begin{array}{rll}
x_1^t &\leq & \left(x_1 - |y_1|\right)^t + |y_1|^t
\;\;\;\;\;\;(x_1 \geq |y_1|),\vt
x_1^t x_2^{1-t} & \leq & \left(x_1 + y_1\right)^t x_2^{1-t} + |y_1|^t x_2^{1-t},
\end{array}
\]
giving (\ref{54}). Similarly if $y_1 \geq 0$, $y_2<0$. If $y_1<0$, $y_2<0$, we have
\[
x_1^t x_2^{1-t} \leq \left((x_1 +y_1)^t + |y_1|^t\right) + \left((x_2 + y_2)^{1-t} +\left|y_2\right|^{1-t}\right)
\]
giving (\ref{54}) again. From (\ref{53}), \ref{54}), and    H\"{o}lder's inequality,
\[
\begin{array}{rll}
|r(x_1+y_1, x_2+y_2) - r(x_1,x_2)| &\leq &
x_1^t|y_2|^{1-t} + |y_1|^t x_2^{1-t} +|y_1|^t|y_2|^{1-t},\vt
\sum_{k=1}^m r_1(x_{1k}, x_{2k};\ve_k) &\leq & \left(\sum_{k=1}^m x_{1k}\right)^t\left(\sum_{k=1}^m \ve_k\right)^{1-t}\;\;
+ \vt
&&
+\left(\sum_{k=1}^m \ve_k\right)^{t}\left(\sum_{k=1}^m x_{2k}\right)^{1-t}
+
\sum_{k=1}^m \ve_k, 
\end{array}\]
which $\rightarrow 0$, if $\sum_{k=1}^m x_{jk}$ ($j=1,2$) are kept fixed, or at least bounded. Hence (\ref{46}) is true for this $r_1$ since $(\D_x)\sum -r(h_1,h_2) \leq 0$.  \nproof
\begin{theorem}
\label{Theorem 24}
Let $h_j(I,x) \geq 0$ and be integrable to $H_j$ ($j=1,2$) in $E$. For some $M \geq 0$, let $h_1(I,x) \leq Mh_2(I,x)$. If $t>1$ is fixed, then $h_1^t/h_2^{t-1}$ is integrable in $E$ to the integral of $H_1^t/H_2^{t-1}$.
\end{theorem}
\textbf{Proof:}
Let
\[
r(x_1,x_2) = \frac{x_1^t}{x_2^{t-1}},\;\;\;\;0\leq x_1\leq Mx_2.
\]
Then
\[
0 \leq \frac{\partial r}{\partial x_1} = \frac{tx_1^{t-1}}{x_2^{t-1}} \leq tM^{t-1},\;\;\;\;
\left|\frac{\partial r}{\partial x_2}\right| = \frac{(t-1)x_1^t}{x_2^t} \leq (t-1) M^t.
\]
Thus (\ref{44}) is satisfied with $A_1 =tM^{t-1}$, $B_1 = (t-1)M^t$. Also, for fixed $b \geq 0$,
\[
1+bc^t - \frac{(1+bc)^t}{(1+b)^{t-1}},\;\;\;\;\;\;(c\geq 0).
\]
(This is $0$ at $c=1$, while its derivative with respect to $c$ is $>0$ in $c>1$, $<0$ in $c<1$.) Putting $x_1=ax_2$, $y_2 = bx_2$, $y_1 =abcx_2$ ($x_j>0$, $y_j>0$, $j=1,2$), we have
\[
\frac{x_1^t}{x_2^{t-1}} + \frac{y_1^t}{y_2^{t-1}} 
- \frac{(x_1+y_1)^t}{(x_2+y_2)^{t-1}}
= a^t x_2 \left(1+bc^t 
- \frac{(1+bc)^t}{(1+b)^{t-1}}\right) \geq 0.
\]
By continuity this is still true as $x_1 \rightarrow 0$ and/or $y_1 \rightarrow 0$, so that (\ref{48}) is true. To show (\ref{49}) true:
\[
(\D_x)\sum \frac{h_1^t(I,x)}{h_2^{t-1}(I,x)} \leq M^{t-1}(\D_x) \sum h_1(I,x) \leq M^{t-1} \left(H_1(E) +1\right)
\]
for all $\D_x$ of $E$ from some suitable $\bS \in \A$ dividing $E$. Hence the theorem follows. \nproof

For fixed $M$ le $h_3$ be real and $|h_3| \leq M$, and let $h_4\geq 0$ be integrable. Then
\[
(\D)\sum |h_3| h_4 \leq M(\D) \sum h_4,
\]
which is bounded for suitable $\bS\in \A$ dividing $E$, and for all $\D$ over $E$ from $\bS$. Hence, if also $h_3h_4$ is integrable, so is $|h_3|h_4$ (Theorem \ref{Theorem 19}), and Theorem \ref{Theorem 24} gives, for $t=2$, the integrability of
\[
h_3^2h_4 = \frac{\left(|h_3|h_4\right)^2}{h_4}.
\]
Another example is, that if $|h_j| \leq Mh_4$ ($j=1,2$), where the $h_1,h_2,h_4$ are real and integrable, then $\frac{h_1h_2}{h_4}$ is integrable. To see this, put $\frac{h_1}{h_4}$ and $h_4$ for $h_3, h_4$ in the previous result, then $\frac{h_1^2}{h_4}$ is integrable, and similarly so are 
\[
\frac{h_2^2}{h_4},\;\;\;\;\;\;\frac{(h_1+h_2)^2}{h_4}.
\]
If $h_5 \geq M>0$, $h_4 \geq 0$, $h_4$ and $h_4h_5$ integrable, then
\[
\frac{h_4}{h_5} = \frac{h_4^2}{h_4h_5}
\]
is integrable. If $M'\geq M>0$, fixed, with $M\leq |h_6| \leq M'$, $h_4 \geq 0$, $h_4$ and $h_4h_6$ real and integrable, then the following are integrable:
\[
h_6^2h_4 = \frac{(h_4h_6)^2}{h_4};\;\;\;\; \;\;\;\;\; 
\frac{h_4}{h_6^2}=\frac{h_4^2}{h_4h_6^2};\;\;\;\;\; \;\;\;\;\; 
\frac{h_4}{h_6} = \frac{\frac{h_4}{h_6^2} \times h_4h_6}{h_4}.
\]
\begin{theorem}
\label{Theorem 25}
Let $r(x)$ be convex with $|r(y)-r(x)|\leq A|y-x|$ for some fixed $A>0$. Let $h \geq 0$ be integrable to $H$, let $f$ be a real valued point function, and let $fh$ be integrable to $H_1$. If $\cS(r(f)h;\bS;E)$ is bounded above, for some $\bS\in \A$ dividing $E$, then $r(f)h$ is integrable. (The function $r$ need only be convex over the range of values of $f$ in $E$.)
\end{theorem}
\textbf{Proof:} 
We can replace $h$ by $H$, disregard those $I \in \D$ with $H(I)=0$, and obtain
\[
\begin{array}{rll}
(\D_x) \sum \left|r(f) H - r\left(\frac{H_1}{H}\right)H\right|
&=& (\D_x)\sum \left|r(f) - r\left(\frac{H_1}{H}\right)\right| H \vt
&\leq& A(\D_x) \sum \left|f-\frac{H_1}{H} \right| 
= A(\D_x)\sum |fH -H_1| \vt
&\rightarrow & 0\;\;\;\mbox{ for suitable }\;\;\bS.
\end{array}
\]
If $I \in \D$ and if $\D'$ is a division of $I$, with $J \in \D'$, then, for $H(I) \neq 0$, we have
\[
\frac{H(J)}{H(I)} \geq 0,\;\;\;\;\;\;(\D') \sum \frac{H(J)}{H(I)} = \frac{H(I)}{H(I)} =1.
\]
Hence, by convexity,
\[
\begin{array}{rll}
r\left(\frac{H_1(I)}{H(I)}\right) H(I) 
&=&
r\left(\frac{(\D') \sum H_1(J)}{H(I)}\right) H(I)  \vt
&=&
r\left((\D') \sum \frac{ H_1(J)}{H(J)} \times \frac{H(J)}{H(I)} \right) H(I) \vt
&\leq &
(\D') \sum r\left(\frac{ H_1(J)}{H(J)} \times \frac{H(J)}{H(I)} \right) H(I) \vt
&=&
(\D') \sum r\left(\frac{ H_1(J)}{H(J)}\right) \times H(J) .
\end{array}
\] 
Hence $r\left(\frac{H_1}{H}\right) H$ is subadditive, and so is integrable by the boundedness condition, and so is $r(f)h$.  

\section{Lebesgue and Denjoy-type theorems: \\ Decomposable division spaces}
The system\footnote{$P$ here is not the partial elementary set of these notes; instead, it is the point (tag-point) of $(I,P)$, in the notation Henstock used in his lecture notes on the Riemann-complete integral. (Note by P.~Muldowney) } using the functions $\delta(P)>0$, and even the system using $\delta (P)$ constant $>0$, are non-additive division spaces, as can easily be checked. But neither is a division space. 

For example, take the dimension $n=1$, and let $a<b<c$. If $\bS_1 \in \A$ divides $[a,b]$, there is a $\delta_1(P) >0$ in $[a,b]$ such that $(I,P)\in \bS_1$ for all $(I,P)$ that are compatible with $\delta_1(P)$, and similarly, if $\bS_2 \in \A$ divides $[b,c]$ there is a corresponding $\delta_2(P)>0$ in $[b,c]$; if $\bS\in \A$ divides $[a,c]$, there is a $\delta(P)>0$ in $[a,c]$ with required properties. 

We can arrange that $\delta(P) \leq \delta_1(P)$ in $a\leq P\leq b$, $\delta(P) \leq \delta_2(P)$ in $b \leq P \leq c$; but whatever the value of $\delta(b)>0$, there are intervals $[u,v]$ with
\[
b-\delta(b) <u<b<v<b+\delta(b)
\]
such that 
\[
([u,v],b) \in \bS\;\; \mbox{ but }\;\; ([u,v],b) \notin \bS_1 \cup \bS_2.
\]
Hence the system cannot be additive, and this is why we need an awkward geometrical theorem.

However, if instead of letting the associated point be anywhere in the closed interval, we suppose that the point is at an end-point of the interval, i.e.~$([u,v],u)$ or $([u,v],v)$, then the system so obtained is additive.

For the intervals connected with $b$ will be either $[u,b]$, when $([u,b],b) \in \bS_1$, or $[b,v]$, when $([b,v],b) \in \bS_2$; and $\bS \subseteq \bS_1 \cup \bS_2$.

Such a system will give the same sums as before. For if $u<x<v$ with $([u,v],x) \in \bS$ (first system), this will contribute a term
\[
f(x)(g(v) - g(u)) = f(x)(g(v) - g(x)) + f(x)(g(x) - g(u))
\]
to a sum over divisions $\D_x$ of $[a,c]$ with $([u,v],x) \in \D_x$.

By repeating over all $[u,v] \in \D$, we get a sum for the second system. Conversely, a sum for the second system is already a sum for the first system.
In $n$ dimensions we use associated points at the vertices of the bricks.

But Lebesgue-type theorems on the limit of a sequence of integrals are false for Riemann integration, so that we need special kinds of division spaces.

First we define
\[
\bS[X] := \{(I,x) \in \bS,\;x \in X\},\;\;\;\;\bS \in \A,\;\;\; X \subseteq T.
\]
Let $E$ be an elementary set, let $(\bS_j)$ be a sequence of members of $\A$, each dividing $E$, and let $(X_j)$ be a sequence of mutually disjoint subsets of $T$. The union of the $X_j$ need not be $T$ or $E$.

If, for all such $(\bS_j), (X_j), E$, there is an $\bS \in \A$, dividing $E$, with
\[
\bS[X_j] \subseteq \bS_j[X_j],\;\;\;\;\;\;j=1,2,3, \ldots,
\]
i.e.~$\bS[X_j] \subseteq \bS_j$, we say that $\A$ has \emph{decomposable families}. If, also, $(T,\T,\A)$ is a division space, we say that is a \emph{decomposable} division space.

Given that, to each $x \in T$, there corresponds an $\bS(x) \in \A$ that divides $E$, i.e.~a mapping from $T$ to $\A$, relative to $E$, if, for each such correspondence, there is an $\bS \in \A$ that divides $E$, with 
\[
\bS[\{x\}] \subseteq \bS(x)[\{x\}],\;\;\;\mbox{ all }\;\;x\in T,
\]
(where $\{x\}$ is the singleton or set containing only the single element $x$) we say that $\A$ has \emph{fully decomposable} families. If, also, $(T,\T,\A)$ is a division space, we say that it is a \emph{fully decomposable} division space.

The system in Euclidean space $T$ of $n$ dimensions that uses functions $\delta(x)>0$ at each point $x$, with associated points at corner points (vertices) of bricks, is a decomposable\footnote{Fully decomposable (P.~Muldowney).} division space.

First we have an important result for the variation.

\begin{theorem}
\label{Theorem 26}
Let $(T,\T,\A)$ be a decomposable division space, and let $(X_j)$ be a monotone increasing sequence of sets in $T$ with union $X$. Then
\[
\lim_{j \rightarrow \infty} V(h;\A;E;X_j) = V(h;\A;E;X).
\]
\end{theorem}
\textbf{Proof:}
As $X \supseteq X_j$, it follows that $V(h;\A;E;X) \geq V(h;\A;E;X)$,
\begin{equation}
\label{55}
V(h;\A;E;X) \geq \lim_{j\rightarrow \infty} V(h;\A;E;X_j).
\end{equation}
Thus if the limit is conventional $+\infty$ the result is true. So we can now assume the limit is finite.
Let $\bS_j \in \A$, dividing $E$, satisfy
\begin{equation}
\label{56}
V(h;\bS_j;E;X_j) < V(h;\A;E;X_j) + \frac \ve{2^j}
\end{equation}
for $j \geq 1$. As $\A$ has decomposable families there is an $\bS \in \A$ dividing $E$ with
\[
\bS[X_j \setminus X_{j-1}] \subseteq \bS_j[X_j \setminus X_{j-1}],\;\;\;\;\;\;\;\;\; j=1,2,3,\ldots, \;\;\;\;X_0 \;\mbox{ empty }.
\]
If $\D_x$ is a division of $E$ from $\bS$, and if $\Q_x, \Q_{jx}$ are the partial divisions from $\D_x$ with $x \in X$, $x \in X_j\setminus X_{j-1}$, respectively ($j=1,2$), there is a greatest integer $m$ (depending on $\D_x$) such that $\Q_{mx}$ is not empty. Let $P_j$ be the partial set from $\Q_{jx}$ ($1 \leq j \leq m$). Then, from Theorems \ref{Theorem 16}, \ref{Theorem 17}, and (\ref{55}), (\ref{56}), (and with $\chi(X,x)$ denoting the indicator or characteristic function of $X$),

\begin{eqnarray*}
(\D_x)\sum ||h(I,x)|| \chi(X,x)
&=& \sum_{j=1}^m ||h(I,x)|| \;\;=\;\;
\sum_{j=1}^m (\Q_{jx}\sum ||h(I,x)|| \vt
&\leq & \sum_{j=1}^m V(h;\bS_j;P_j;X_j)
\vt
&<& \sum_{j=1}^m \left(V(h;\A_j;P_j;X_j) + \frac \ve{2^j} \right) \vt
&<& \sum_{j=1}^m V(h;\A;P_j;X_m) + \ve \vt
& \leq &
 V(h;\A;E;X_m) + \ve \vt
& \leq & \lim_{j \rightarrow \infty} V(h;\A;E;X_j) + \ve.
\end{eqnarray*}
Hence
\[
V(h;\A;E;X) \leq
V(h;\bS;E;X) \leq \lim_{j\rightarrow \infty} V(h;\A;E;X_j) + \ve,
\]
giving the opposite inequality to (\ref{55}), and so the result.  \nproof

Theorem \ref{Theorem 26} is false for Riemann integration, for let $X_j$ be the set of the first $j$ rationals in $0<x<1$, let $E=[0,1)$, and let $h([u,v)) = v-u$. Then
\[
V(h;\A;E;X_j) = 0,\;\;\;\;\;\;V(h;\A;E;X) = 1.
\]

\begin{theorem}
\label{Theorem 27}
Let $(t,\T,\A)$ be a decomposable division space, and let $(X_j)$ be a sequence of subsets of $T$ with union $X$. Then
\[
V(h;\A;E;X) \leq \sum_{j=1}^\infty
V(h;\A;E;X_j). 
\]
\end{theorem}
\textbf{Proof:}
$V(h;\A;E;Y_1)\leq V(h;\A;E;Y_2)$ if $Y_1 \subseteq Y_2$. Hence we can assume the $X_j$ disjoint. Also we can take the right hand side finite (or else there is nothing to prove). We begin the proof as in the proof of Theorem \ref{Theorem 26}, replacing $X_j\setminus X_{j-1}$ by $X_j$ and obtaining

\begin{eqnarray*}
(\D_x)\sum ||h(I,x)|| \chi(X,x)
&\leq & \sum_{j=1}^m V(h;\bS_j;P_j;X_j)
\vt
&\leq & \sum_{j=1}^m V(h;\bS_j;E;X_j)\vt
&<& \sum_{j=1}^m \left(V(h;\A;E;X_j) + \frac \ve{2^j} \right), \vt
 V(h;\A;E;X)  
& \leq & 
V(h;\bS;E;X) \vt
&\leq & \sum_{j=1}^\infty V(h;\A;E;X_j) + \ve.
\end{eqnarray*}
Hence the result.   \nproof

An analogue of Theorem \ref{Theorem 27} can easily be proved.

\begin{theorem}
\label{Theorem 28}
Let $(T, \T, \A)$ be a decomposable division space, and let $C \subseteq K$ be a compact set. IF $\vs(h;\bS;E;X) \subseteq C$ for some $\bS \in \A$ dividing $E$, then
\[
\VS(h;\A;E;X) \subseteq \mbox{Cl} \bigcup_{m=1}^\infty
\left\{\VS(h;\A;E;X_1) + \cdots +
\VS(h;\A;E;X_m)\right\}.
\]
\end{theorem}
\textbf{Proof:} (Analogous.)

\begin{theorem}
\label{Theorem 29}
Let $(T,\T,\A)$ be a decomposable division space and let $K$ have scalar multiplication with $||bk||\leq |b|.||k||.$ ($k\in K$, $b$ real or complex).
For a scalar point function $f(x)$, if 
\begin{equation}
\label{57}
V(h;\A;E;X)=0\;\;\;\mbox{ then }\;\;\;V(fh;\A;E;X)=0.
\end{equation}
For a scalar point function $f(x)\neq 0$ in a set $X_1 \subset X$, 
\begin{equation}
\label{58}
\mbox{if } \;\;\;V(fh;\A;E;X) =0\;\;\;\mbox{ then }\;\;\;V(h;\A;E;X_1)=0.
\end{equation}
Let $K$ be the real line or complex plane. If $h_1,h_2$ are variationally equivalent, if $f$ is a point function, and 
\begin{equation}
\label{59}
\mbox{if } \;\;\;\int_Efdh_1\;\;\;\mbox{ exists then so does }\;\;\;\int_Efdh_2
\end{equation}
and they are equal. Conversely if both integrals exist and are equal over each partial set $E$ of an elementary set $E_1$, then $f=0$ except in a set over which $h_1-h_2$ has variation zero. 
\end{theorem}
\textbf{Proof:}
Let $X_j$ be the subset of $X$ with $|f|\leq j$ ($j=2,3, \ldots$). By Theorem \ref{Theorem 27}, (\ref{57}) follows from
\[
V(fh;\A;E;X) \leq \sum_{j=1}^\infty V(fh;\A;E;X_j) \leq 
\sum_{j=2}^\infty jV(h;\A;E;X) =0.
\]
Then (\ref{58}) follows on replacing $X$ and $f$ by $X_1$ and $1/f$ in (\ref{57}). Then for (\ref{59}) we use (\ref{57}) and Theorem \ref{Theorem 15} (\ref{34}) one way; and the linearity of the integral and Theorem \ref{Theorem 15} (\ref{34}) and (\ref{58}) for the other way. \nproof

\section{Limits of Integrals}
When the linear space $K$ has a norm $|| \cdot ||$, the property
\begin{equation}
\label{60}
||\int_E f_j(x)d\mu - \int_E f(x)d\mu|| \rightarrow 0\;\;\mbox{ as }\;\; j \rightarrow \infty
\end{equation}
holds if
\begin{equation}
\label{61}
\mu \geq 0;
\end{equation}
\begin{equation}
\label{62}
|| f_j(x)- f(x)|| \rightarrow 0\;\;\mbox{ almost everywhere on }\;\;E;
\end{equation}
\begin{equation}
\label{63}\mbox{each }  f_j(x) \mbox{ is integrable wrt } \mu \mbox{ on }E;
\end{equation}
\begin{equation}
\label{64}
||f_j(x)|| \leq F(x)
\end{equation}
where $F(x)$ is some real-valued point-function integrable wrt $\mu$ on $E$, as we shall see.

But not all our topological groups $K$ have even a group norm, so that to obtain similar results we have to consider different sufficiency conditions.

As usual, we generalize $f_j(x) \mu(I)$ to become $h_j(I,x)$, and we have to modify (\ref{62}).

To define the convergence of $h_j(I,x)$ to $h(I,x)$ we suppose that to each $I \in \T$ there corresponds a sequence $(Y^k(I))$ of sets in $K$ that contain the zero $z$, with the property that if $G$ satisfies $z \in G\in \G$, there are a positive integer $k$ and an $\bS \in \A$ dividing $E$, such that for all divisions $\D$ over $E$ from $\bS$,
\begin{equation}
\label{65}
(\D)\sum Y^k(I) := \left\{(\D) \sum y^k(I): y^k(I)\in Y^k(I) \right\} \subseteq G.
\end{equation}
Then we suppose that there are a set $X \subseteq T$, an $\bS_1 \in \A$ dividing $E$, and integers $l=l(k,x)$ ($k \geq 1$) such that
\begin{equation}
\label{66}
h \mbox{ and } h_j \;\;\;(j \geq 1)\;\;\;\mbox{ are of variation zero in } X \mbox{ relative to } E,\A;
\end{equation}
\begin{equation}
\label{67}
h_jI,x) - h(I,x) \in Y^k(I),\;\;\;\;\;((I,x) \in \bS_1,\;\; x \notin X,\;\; j \geq l(k,x)).
\end{equation}
\begin{theorem}
\label{Theorem 30}
Let $(T,\T,\A)$ be a division space, and let $K$ be complete. If $(h_j)$ is a sequence of interval-point functions $h_j(I,x)$ integrable to $H_(P)$ over each partial set $P$ of $E$, and if $h(I,x)$ is an interval-point function such that (\ref{65}), (\ref{66}), (\ref{67}) hold with $l(k,x) = l(k)$, independent of $x$, then $H_j(P)$ tends to a limit $H(P)$, and $h$ is integrable to $H(P)$ over each partial set $P$ of $E$.
\end{theorem}
\textbf{Proof:}
As $z \in Y^k(I)$ for each $I \subseteq E \setminus P$, the convergence condition\footnote{\textbf{Note by P.~Muldowney} We want condition (59)': Given $z\in G\in \G$, $\exists k$, $\bS \in \A$ dividing $P$, such that for all $\Q$ from $\bS$,
\[
(\Q)\sum Y^k(I) = \left\{(\Q)\sum y^k(I): y^k(I) \in Y^k(I) \right\} \subseteq G;
\]
so, in (\ref{65}), if $(I,x) \in \D$, $I \subseteq E \setminus P$, we take $y^k(I) = z$, and then $(\D)\sum Y^k(I)$ reduces to $(\Q)\sum Y^k(I)$.}
condition holds for $P$ if it holds for $E$. Hence we can take $P=E$. Also by (\ref{66}) and Theorem \ref{Theorem 11} [14] we can assume $X$ empty. Let $z \in G\in \G$, let $s>t \geq l(k)$ and let $\bS_{st} \in \A$, dividing $E$, be such that $\bS_{st}\subseteq \bS \cap \bS_1$, and such that for all divisions $\D_x$ of $E$ from $\bS_{st}$,
\[
(\D_x) \sum h_s(I,x) - H_s(E) \in G,\;\;\;\;\;\;
(\D_x) \sum h_t(I,x) - H_t(E) \in G.
\]
Using the same $\D_x$ with (\ref{65}), (\ref{67}),
\begin{eqnarray*}
H_s(E) - H_t(E) &=& \left((\D_x) \sum h_t(I,x) - H_t(E) \right) - \vt
&&-\left((\D_x) \sum h_s(I,x) - H_s(E) \right) + \vt
&& + (\D_x) \sum \left(h_s(I,x) - h(I,x) \right) -\vt
&&-(\D_x) \sum \left(h_t(I,x) - h(I,x) \right) \vt
&\in & G-G + (\D) \sum Y^k(I) - (\D) \sum Y^k(I)  \vt
&\subseteq & G-G+G-G.
\end{eqnarray*}
As $K$ is a topological group, given $G^*$ in $z \in G^* \in \G$, by choice of $G$ we have 
\[
G-G+G-G \subseteq G^*.
\]
But $H_s(E), H_t(E)$ are independent of $\bS_{st}$ so that there is a least value $l^*$ of $l(k)$, depending only on $G^*$, for which
\[
H_s(E) - H_t(E) \in G^*\;\;\;\;\;(s>t \geq l^*).
\]
Hence $(H_s(E))$ is fundamental (or Cauchy) and so convergent as $K$ is complete. If the limit is $H(E)$, then for all divisions $\D_x$ of $E$ from $\bS_{t+1,t}$,
\begin{eqnarray*}
(\D_x)\sum h(I,x) - H(E) &=& \left((\D_x) \sum h_t(I,x) - H_t(E) \right) + \vt
&&+\left(H_t(E) - H(E) \right) + \vt
&& +
(\D_x) \sum \left(h(I,x) - h_t(I,x) \right) \vt
&\in & G+H_t(E) -H(E) -(\D) \sum Y^k(I) \vt
&\subseteq & G-G + H_t(E) -H(E).
\end{eqnarray*}
By choice of $G$ and $t$ and so of $\bS_{t+1,t}$, this last set is contained in an arbitrary open nbd of $z$, and the theorem is proved\footnote{\textbf{Note:} Decomposability not used.}.  \nproof

It may not always be possible or convenient to prove that (\ref{67}) holds for an $l(k,x)$ independent of $x$ so that we need another approach.

We apply a real continuous linear functional $\mathcal{F}$ to obtain real-valued interval-point functions $\mathcal{F}h_s(I,x)$. We need extra conditions, the most useful being those invariant for $\mathcal{F}$.

For example the sets $\mathcal{F}Y^k(I)$ will be shown to have the same properties as the sets $Y^k(I)$. Also compacr sets (?) are invariant for $\mathcal{F}$ as $\mathcal{F}$ is continuous. But open sets are not invariant and this is one reason why we do not assume the $Y^k(I) $ open. Another reason is that we do not need this assumption even though it is easier to construct open $Y^k(I)$.

We first have an approximation theorem.

\begin{theorem}
\label{Theorem 31}
Let $(T,\T,\A)$ be a decomposable division space, $C$ a compact set, and $h,h_j$ ($j \geq 1$) interval-point functions such that (\ref{65}), (\ref{66}), (\ref{67}) hold with
\begin{equation}
\label{68}
(\D_x) \sum h_{j(I,x)} (I,x) \in C
\end{equation}
for all divisions $\D_x$ of $E$ from an $\bS_2 \in \A$, dividing $E$, and all choices $j(I,x)$ of integers $1,2,3, \ldots$ for $(I,x) \in \D_x$. Then the non-empty
\begin{equation}
\label{69}
\bar{\cS}(h;\A;E) \subseteq C
\end{equation}
and Theorem \ref{Theorem 9} (\ref{7}) holds.
If also the $h_j$ are integrable to $H_j$ on the partial sets of $E$, for each $j \geq 1$, then, given an arbitrary integer $N$, there is an $\bS_3 \in \A$, dividing $E$, such that all divisions $\D_x$ of $E$ from $\bS_3$ satisfy
\begin{equation}
\label{70}
(\D_x) \sum h(I,x) - \sum_{j=N}^q H_j(P_j) \in G-G
\end{equation}
for certain disjoint partial sets $P_N, \ldots , P_q$, with union $E$, formed from the intervals of $\D$. Further if $K$ is the real line with $h_j(I,x)$ monotone increasing in $j$ for each fixed $(I,x)$, then $H_j(E)$ tends to a limit $H(E)$ as $ j \rightarrow \infty$, and $h$ is integrable over $E$ to $H(E)$.
\end{theorem}

\noindent
\textbf{Proof:} By (\ref{66}) and Theorem \ref{Theorem 11} [14] we can replace $h_j(I,x)$ and $h(I,x)$ by $z$ in $X$ and so take $X$ empty. Let $l_0(\D_x)$ be the maximum of  $l(k,x)$ for $(I,x) \in \D_x$. Then, from (\ref{65}), (\ref{67}), 
\[
(\D_x)\sum h_j(I,x) - (\D_x)\sum h(I,x)  \in (\D)\sum Y^k(I) \subseteq G\;\;\;\;\;\;(j \geq l_0(\D_x)).
\]
Now the compact $C$ is closed as $\G$ is Hausdorff. As $G$ is an arbitrary open nbd of $z$ (\ref{68}) shows that
\[
(\D_x) \sum h(I,x) = \lim_{j \rightarrow \infty}(\D_x)\sum h_j(I,x) \in \bar C = C,\;\;\;\;\cS(h;\bS;E) \subseteq C.
\]
Then the conditions of Theorem \ref{Theorem 9} , and we have (\ref{7}). That is, $
\exists \bS_1$  such that \[ \cS(h;\bS_1,E) \subseteq \left\{\bar{\cS}(h;\A;E) + G\right\}\cap C.
\]
When the $h_j$ are integrable to $H_j$, let $G$ be given in $z \in G\in \G$. Then we can find a sequence $(G_j) \subseteq \G$ that satisfies\footnote{See Henstock, \emph{Theory of Integration}, Page 129.}
$z \in G_j \in \G$ (all $j$), $G_0=G,$
\begin{equation}
\label{71}
\sum_{N<j \leq q}G_j \subseteq G_N\;\;\;(q>N),\;\;\;G_j-G_j \subseteq G_{j-1}\;\;\;(j\geq 1).
\end{equation}
By Theorem \ref{Theorem 10} (\ref{9}) there is an $\bS_j \in \A$, dividing $E$, such that if a division $\D_x$ over $E$ is from $\bS_j$, with a partial division $\Q_x$ of $\D_x$ forming a partial set $P$, then
\begin{equation}
\label{72}
(\D_x) \sum h_j(I,x) - H_j(E) \in G_j,\;\;\;
(\Q_x) \sum h_j(I,x) - H_j(P)\in G_j-G_j \subseteq G_{j-1}.
\end{equation}
For a given integer $N$ let $X_{kj}$ be the set of all $x$ with $l(k,x)=j$ ($j>N$), and let $X_{kN}$ be the set of all $x$ with $l(k,x)\leq N$. As $\A$ has decomposable families, there is an $\bS^* \in \A$, dividing $E$, with
\[
\bS^*[X_{kj}] \subseteq \bS_j[X_{kj}],\;\;\;\;\;(j\geq N).
\]
If $\D_x$ over $E$ is from $\bS^*$ then, by (\ref{67}),
\begin{equation}
\label{73}
(\D_x)\sum \left(h_{l(k,x)}(I,x) - h(I,x)\right) \in (\D) \sum Y^k(I) \subseteq G,
\end{equation}
where we replace $l(k,x)$ by $N$ if $l(k,x) <N$. We group the $(I,x)$ in $\D_x$ into partial divisions $\Q_{N,x}, \ldots , \Q_{q,x}$, where $I \in \Q_j$ and where $q=\max(j)$ for the particular $\D_x$. Let $P_j$ be the union of the intervals of $Q_j$. By (\ref{71}), \ref{72}), (\ref{73}),
\[
(\D_x) \sum h_{l(k,x)}(I,x) - \sum_{j=N}^q H_j(P_j) \in \sum_{j=N}^q G_{j-1} \subseteq G_{N-2} \subseteq G, \;\;\;(N \geq 2).
\]
This gives (\ref{70}). For the last part of the theorem we take $G=(-\ve/2, \ve/2)$ in (\ref{70}), for arbitrary $\ve>0$. Then
\begin{equation}
\label{74}
\left|(\D_x)\sum h(I,x)-\sum_{j=N}^q H_j(P_j)\right| < \ve.
\end{equation}
By monotonicity of the $h_j(I,x)$ in $j$, $H_j(P)$ is monotone increasing in $j$. It is finitely additive for partial sets $P$ of $E$. In (\ref{68}) taking $j(I,x)=j$, constant, for all $(I,x)$, we have $H_j(E) \in \bar C = C$, and $H(E) = \lim_{j\rightarrow \infty} H_j(E)$ exists. Hence, by (\ref{74}), and for $N$ large enough, 
\[
\begin{array}{rll}
H_N(E) &=& \sum_{j=N}^q H_N(P_j) \;\;\leq \;\; \sum_{j=N}^q H_j(P_j) \vt
&\leq & \sum_{j=N}^q H_q(P_j) \;\;=\;\;H_q(E)\;\;\leq \;\;H(E),
\end{array}
\]
$0 \leq H(E) - H_N(E) < \ve, $
giving
$H(E) - 2 \ve < (\D_x)\sum h(I,x) < H(E) + \ve$, so $h$ is integrable, and $H$ is its integral.  \nproof

To show that it is not enough to assume the monotone convergence of $h_j(I,x)$ in $j$ for each fixed $(I,x)$,  omitting (\ref{67}), we take $T=[0,1)$, $K$ the real line, and $\T$ the set of half-closed intervals $[a,b)$ in $T$, and we replace the sequences by series.

\noindent
\textbf{Example:}
For $j>1$ let
\[
\begin{array}{rll}
h_j([u,v),x)& := &\left\{
\begin{array}{rl}
v-u & (\frac{j+1}{j} u<v\leq \frac j{j-1} u,\vt
0& \mbox{otherwise};

\end{array} \right.
\vspace{10pt}
 \\
\vspace{10pt}
h([u,v),x) &:=& \sum_{j=2}^\infty h_j([u,v),x) \vt
&=& \left\{
\begin{array}{rl} v-u & (u<v \leq 2u,\;\;\;x=u,v, \vt
0 & \mbox{otherwise}.
\end{array}\right.
\end{array}
\]
(\textbf{Note by P.~Muldowney:} What if $u=0$ in def of $h_j$ above?)

\noindent
Then $\int_T h=1$. But if we take $\delta_j(0)=1$, $\delta_j(u)=u/j$ ($u>0$), then $h_j([u,v), x) =0$ for all $[u,v)$ compatible with $\delta_j(u)$, and hence $\int_T h_j =0$ for all $j$. 

Also it can be proved that $\sum_{j=1}^\infty h_{2j}([u,v),x)$ is not integrable, while
\[
\sum_{j=2}^\infty jh_j ([u,v),x)
\] gives unbounded sums. This last shows that, from the boundedness of the $H_j(P)$ of Theorem \ref{Theorem 31}, last part, we cannot deduce (\ref{67}) or (\ref{68}) even when the $h_j$ are monotone increasing in $j$.

\noindent
\textbf{Notes on Example, by P.~Muldowney:}

\noindent
$\delta_j(u) = u/j$:

If $x=u$, should have $v<u+u/j = \frac{j+1}j u$, contradicts $\frac{j+1}j u<v$. Therefore $h_j([u,v), u) =0$.

?If $x=v$, $u>v-v/j = \frac{j-1}j v$.

\noindent
If $u, j$ given, $u<v$ and $v \leq 2u$ then [DIAGRAM LOCATING $v$]. There exists $j \geq 2$, such that $\frac{j+1}j u<v \leq \frac j{j-1} u$, and this is the only such $j$. On the other hand,
\[
v>2u \Longrightarrow v>\frac j{j-1}u\;\;\; \forall j\;\;\; ; u\geq v \Longrightarrow \frac{j+1}j u \geq v\;\;\; \forall j.
\]

*********************

However, better results are true for special $h_j$.
\begin{theorem}
\label{Theorem 32}
Let $(T,\T,\A)$ be a decomposable division space, let $h(I) \geq 0$ be an interval function, and, for each $x \in T$ let $f(x,y)$ be monotone increasing in $y \geq 1$. If, for each $y \geq 1$, $f(x,y$ is integrable with respect to $h$ in an elementary set $E$, with indefinite integral $H(\cdot, y)$, and f $H(E,y)$ is bounded above as $y \rightarrow \infty$, then 
\[
f(x) := \lim_{y \rightarrow \infty} f(x,y)
\]
exists (finitely) for all $x$ save a set $X$ with $V(h;\A;E;X) =0$. Putting $f(x,y)=0=f(x)$ in $X$, then $f(x)$ is integrable with respect to $h$ in $E$, with indefinite integral $H$,
\[
H=\lim_{y \rightarrow \infty} H(\cdot , y).
\]
\end{theorem}
\textbf{Proof:} (See Henstock, \emph{Linear Analysis}, page 238.)

We now consider more general spaces $K$ of values. Assume that real continuous linear functionals $R$ exist to separate all points of $C$, i.e.~if $a,b \in C$, $a \neq b$, there is such an $R$ with $R(a) \neq R(b)$.

As $R$ is continuous, the set $G \subseteq K$ where, for fixed $a$, 
\[
|R(x)-R(a)|< \rho = |R(b)-R(a)!,
\]
is open. If $x_1,x_2 \in G$, and $x$, are such that $x+x=x_1+x_2$, then
\[
\begin{array}{rll}
2|R(x)-R(a)| &=& |R(x+x)- 2R(a)| = |R(x_1) +R(x_2) - 2R(a)| \vt &\leq & (R(x_1)-R(a)| + |R(x_2)-R(a)| < 2\rho,
\end{array}
\]
so that $x\in G$. Thus $G$ is an open convex nbd of $a$ that does not contain $b$. Being true for all pairs of distinct points of $C$, we can say that $C$ is \emph{locally convex} in $K$. The converse may also hold.

\begin{theorem}
\label{Theorem 33}
Let $(T,\T,\A)$ be a decomposable division space, let $K$ be a real linear space, and let $C \subseteq K$ be a compact set such that real continuous linear functionals exist to separate all points of $C$. Let $h \geq 0$ be an interval-point function and let $f_j$ ($j\geq 1$) and $f$ e point functions with values in $K$, such that $f_jh$ and $fh$ satisfy all conditions on $h_j, h$ respectively in Theorem \ref{Theorem 31}, second part (i.e.~(\ref{65}), (\ref{66}), (\ref{67}), (\ref{68}), and the integrability of $f_jh$ to $K_j$ in $E$). [NOTE: not assuming $l(k,x) = l(k)$.] Then the limit of $K_j(E) = \int_E f_jdh$ exists as $j \rightarrow \infty$, say as $K(E)$, and $fh$ is integrable over $E$ to $K(E)$.
\end{theorem}
\textbf{Proof:}
By Theorem \ref{Theorem 31} (\ref{69}), $\bar \cS(fh;\A;E) \subseteq C$ and is not empty. Let it contain two points $a \neq b$, so that there is a real continuous linear functional $R$ with $R(a) \neq R(b)$. By linearity, (\ref{68}), and continuity,
\[
(\D_x)\sum R\left(f_{j(I,x)}(x)\right) h(I,x) = R\left((\D_x)\sum f_{j(I,x)} (x)h(I,x) \right) \in R(C),
\]
a compact set. Hence $R(f_j)h$ satisfies (\ref{68}). Also by linearity,
\[
\begin{array}{rl}
&(\D_x)\sum R\left(f_{j(I,x)}(x)\right) h(I,x) - R\left(K_j(E)\right) 
\vt
=& R\left((\D_x)\sum f_{j(I,x)} (x)h(I,x) -K_j(E)\right) .
\end{array}.\]
As $R(z)=0$, and by continuity, if $\ve>0$, there is an open nbd $G$ of $z$ with $R(G) \subseteq (-\ve, \ve)$. Hence $R(f_j(x))h(I,x)$ is integrable over $E$ to $R(K_j(E))$. Further, from linearity and (\ref{67}), for all $j \geq l(k,x)$, 
\[
\begin{array}{rl}
&R\left(f_j(x)\right)h(I,x) - R\left(f(x)\right) h(I,x) \vt
=& R\left(f_j(x)-f(x) h(I,x) \right) \in R\left(Y^k(I)\right), = Z^k(I) \;\;\mbox{ say},
\end{array}
\]
As $R$ is additive and continuous, then,given $\ve>0$, and a $G$ depending on $\ve$, with $z \in G\in \G$, we have
\[
\begin{array}{rll}
(\D)\sum Z^k(I) &=& \left\{ (\D) \sum R\left(y^k(I)\right): y^k(I) \in Y^k(I)\right\} \vt
&=& 
\left\{ R\left((\D) \sum y^k(I)\right): y^k(I) \in Y^k(I)\right\} \vt
& \subseteq & R(G) \;\;\;\;\subseteq \;\;\;\;(-\ve, \ve).
\end{array}
\]
As $R(z)=0$, the $Z^k(I)$ have the same properties as the $Y^k(I)$. Thus we have reduced the problem to the case $K=$ the real line, and temporarily we can omit the $R$. By Theorem \ref{Theorem 31} (\ref{68}) and Teorem \ref{Theorem 22}, (\ref{51}) implying (\ref{50}), we have the integrability of
the functions
\[
\max\left(f_N(x), \ldots , f_q(x) \right) h(I,x),\;\;\;\;\;\;
\min\left(f_N(x), \ldots , f_q(x) \right) h(I,x),
\]
for each $q,N$ in $q>N \geq 0$. By Theorem \ref{Theorem 32}, the following exist except in a set of $X_1$ of $h$-variation zero, and are integrable in $E$:
\begin{eqnarray*}
\inf_{j\geq N} f_j(x) h(I,x) &=& \lim_{q \rightarrow \infty} \min_{N \leq j \leq q} f_j(x) h(I,x), \vt
\sup_{j\geq N} f_j(x) h(I,x) &=& \lim_{q \rightarrow \infty} \max_{N \leq j \leq q} f_j(x) h(I,x).
\end{eqnarray*}
Let $X_2(\ve)$ be the set of $x \notin X_1$ where, for an infinity of $j$, $|f_j(x)-f(x)| \geq \ve$, ($\ve>0$), and let $\delta >0$. Then by (\ref{65}), (\ref{67}), for $G= (-\ve \delta, \ve \delta)$,
\[
\ve(\D_x)\sum h(I,x) \chi (X_2(\ve;x)) < \ve \delta,\;\;\;\;\;\;V(h;\A;E;X_2(\ve))=0.
\]
Taking $\ve = \frac 1r$ ($r=1,2,3, \ldots$), with Theorem \ref{Theorem 27}, we have
 $\bigcup_{r=1}^\infty X_2(\frac 1r)$ of $h$-variation zero.  Hence
 \[
 fh=\lim_{j \rightarrow \infty} f_jh = \lim_{N \rightarrow} \inf_{j \geq N} f_jh
 \]
 exists except in a set $X_3$ f $h$-variation zero. Further $fh$ is integrable in $E$ since $\inf_{j \geq N}f_j h$ is monotone increasing with $N$. Also,
\begin{eqnarray*}
\int_Efdh &=& \int_E \lim_{N \rightarrow \infty} \inf_{j \geq N} f_j(x) dh \vt
&=& \lim_{N \rightarrow \infty} \lim_{q \rightarrow \infty} \int_E \min_{N\leq j \leq q} f_j(x) dh \vt
&\leq & \lim_{N \rightarrow \infty} \lim_{q \rightarrow \infty}  \min_{N\leq j \leq q} \int_E f_j(x) dh \vt
&=& \liminf_{j \rightarrow \infty} \int_E f_j(x) dh \vt
& \leq & \limsup_{j \rightarrow \infty} \int_E f_j(x) dh,
\end{eqnarray*}
and, by a similar argument this is $\leq \int_E fdh$. Hence the result when $K$ is the real line. In the general case this proof contradicts $R(a) \neq R(b)$, Hence $\bar \cS9fh;\A;E)$ contains just one point $a$. For all the given $R$,
\begin{equation}
\label{75}
R(a) = \lim_{j \rightarrow \infty} \int_E R\left(f_j(x)\right) dh.
\end{equation}
Also we have
\begin{equation}
\label{76}
K_j(E):= \int_E f_j(x) dh \in C.
\end{equation}
By Theorem \ref{Theorem 6}, (\ref{76}) shows that $K_j(E)$ has at least one limit point, $b$ say, and for each such limit point, the sequence $\left(R(K_j(E))\right)$ has a limit point $R(b)$ by continuity of $R$. From (\ref{75}) $R(b) = R(a)$, so that, by choice of $R$ we have 
\[
b=a,\;\;\;\;\;\lim_{j\rightarrow \infty} K_j(E) = a,
\]
and the theorem is proved.   \nproof

Theorem \ref{Theorem 30} covers the case of uniform convergence and a normed linear space $K$. If $\mu\geq 0$ is finitely additive, if $f_j\mu$ and $f\mu$ are integrable over $E$, if $\mu(E)$ is finite, and if $||f_j(x)-f(x)|| \rightarrow 0$ uniformly in $x$ as $j \rightarrow \infty$, then $Y^k(I)$ can be the open sphere with centre $z$ and radius $\mu(I)/k$, say, $OS(\mu(I)/k)$, and the sum of the
$Y^k(I)$ would lie in $OS(\mu(E)/k)$, small for large $k$. For the we have
\[
\begin{array}{rll}
||f_j(x) - f(x)|| & < &\frac 1k\;\;\;(j \geq j_0(k)), \vspace{15pt}\vt\left(f_j(x) - f(x)\right) \mu(I)|| &<& \frac{\mu(I)}k, \vspace{15pt} \vt \left(f_j(x) - f(x)\right) \mu(I) &\in &Y^k(I).
\end{array}
\]
A result of D.~Przeworska-Rolewicz and S.~Rolewicz, ``\emph{On integrals of functions with values in a complete linear metric space}'', shows that the uniform convergence result may fail if the distance from $z$ is not linear. For fixed $p$ in $0<p<1$, let $K=L^p[0,1]$, and let 
\[
f_j(x,y) = j . \chi\left(\left[ \frac{k-1}j, \frac k j\right),y\right), \;\;\;\;(k-1 \leq jx<k,\;\;\;k=1, 2, \ldots ,j).
\]
Then
\[
||f_j(x, \cdot)||_p = \int_{[0,1]}|f_j(x,y)|^p dy = j^{p-1} \rightarrow 0
\]
uniformly in $x$ as $j \rightarrow \infty$. But
\[
\int_{[0,1]}f_j(x,y) dy = \sum_{k=1}^\infty \frac{j . \chi\left(\left[ \frac{k-1}j, \frac k j\right),y\right)}j
= \chi\left([0,1],y\right),
\]
and the latter does not tend to zero. The failure of the conclusion of Theorem \ref{Theorem 30} in this case is probably due to the following. If $m := \max \{\mu(I): I \in \D\}$ then
\begin{eqnarray*}
||f_j(x) \mu(I) || &=& j^{p-1} \mu^p(I),\vt
(\D) \sum \mu^p(I) &=& (\D) \sum \frac{\mu(I)}{\mu^{1-p}(I)} \vt
&\geq & \frac{(\D) \sum \mu(I)}{m^{1-p}} \vt
&=& \mu(E) m^{p-1}\vt & \rightarrow & \infty
\end{eqnarray*}
as $m \rightarrow \infty$. Here uniform convergence of $f_j$ in the ordinary sense does not induce uniform convergence of $f_jh$ in the sense of using $\left(Y^k(I)\right)$. 

On the other hand, when $K$ is the real line or complex plane, this latter convergence lies deeper than ordinary uniform convergence.
\begin{theorem}
\label{Theorem 34} (Egoroff) Let $(T, \T, \A)$ be a decomposable division space and let $f_j(x) \geq 0$ and tend to zero as $j \rightarrow \infty$, except in a set of $h$-variation zero, where $h(I,x) \geq 0$. If $h$ is integrable in $E$, and if each $f_j$ is $h$-measurable (i.e.~the characteristic function of the set $X(j, \ve)$, where $f_j(x)<\ve$, is integrable with respect to $h$ for each $\ve>0$) then, given $\ve>0$, there is a set $X(\ve)$ with $V(h:\A;E;X(\ve))<\ve$, such that $\left(f_j(x)\right)$ is uniformly convergent to $0$ in $\setminus X(\ve)$.
\end{theorem}
\textbf{Proof:}
For the integrability of $h$, and since every characteristic function is bounded, Theorems \ref{Theorem 22} and \ref{31} (last part) show that
\begin{eqnarray*}
\min_{N\leq j \leq q} \chi \left(X(j, \ve),x\right) &=&
\chi\left(\bigcap_{j=N}^q X(j, \ve),x\right), \vt
\inf_{ j \geq N} \chi \left(X(j, \ve),x\right) &=&
\chi\left(\bigcap_{j=N}^\infty X(j, \ve),x\right),
\end{eqnarray*}
are integrable. The set
$Y(N, \ve) := \bigcap_{j=N}^\infty X(j, \ve)$ satisfies $f_j(x)< \ve$ (all $j>N$), and 
\[
Y(M, \ve) \subseteq Y(N, \ve) \;\;\;\;\;\;(M<N).
\]
Also in $\setminus \bigcup_{N=1}^\infty Y(N, \ve)$ we have $f_j(x) \geq \ve$ for an infinity of $j$. As $f_j(x) \rightarrow 0$ except in a set of $h$-variation zero, we have
\begin{eqnarray*}
V\left(h;\A;E; \setminus \bigcup_{N=1}^\infty Y(N,\ve)\right) &=&0,\vt
V\left(h;\A;E;  \bigcup_{N=1}^\infty Y(N,\ve)\right) &=& V(h;\A;E).
\end{eqnarray*}
By Theorem \ref{Theorem 26}, given $\delta>0$, there is an integer $N$ depending on $\ve>0$, $\delta>0$, such that
\begin{eqnarray*}
V\left(h;\A;E; Y(N,\ve)\right) &>& V(h;\A;E) - \delta,\vt
\mbox{i.e. }\; \int_E \chi\left(Y(N, \ve),x\right) dh &>& \int_E h - \delta, \vt
\int_E \chi\left(\setminus Y(N, \ve),x\right) dh &<&  \delta, \vt
V\left(h;\A;E;  \setminus Y(N,\ve)\right) &<& \delta.
\end{eqnarray*}
For $\ve>0$, $\delta >0$ we now substitute $1/k$ and $\delta 2^{-k}$, so that $N = N(k, \delta)$.
Putting
\[
X(\delta) = \bigcup_{k=1}^\infty \setminus Y\left(N(k, \delta), \frac 1k\right),
\]
Theorem \ref{Theorem 27} gives $V(h;\A;E;X(\delta))<\delta$, while in 
\[
\setminus X(\delta) = \bigcap_{k=1}^\infty Y\left(N(k,\delta), \frac 1k\right)
\]
we have
\[
\sup_{j\geq N(k,\delta)}f_j(x) \leq \frac 1k,
\]
i.e.~$f_j \rightarrow 0$ uniformly in $\setminus X(\delta)$.

\begin{theorem}
\label{Theorem 35}
Let $(T, \T,\A)$ be a decomposable division space, let $\mu(I)\geq 0$ be integrable in $E$, let $F(x)\geq 0$, let $f_j(x)$ and $f(x)$ have values in $K$, and let $F$, $||f_j-f||$ be integrable in $E$ with respect to $\mu$ ($j \geq 1$), with $||f_j(x) -f(x)|| \rightarrow 0$ as $j \rightarrow \infty$, except in a set of $\mu$-variation zero, and $||f_j(x)|| \leq F(x)$ (all $j,x$). If Vitali's theorem holds in $T$, then (\ref{65}), (\ref{66}), (\ref{67}) are true for suitable $Y^k$.
\end{theorem}
\textbf{Proof:}
$||f_j-f||$ is integrable and so measurable (see M313 proof --- Henstock Lecture Notes?). Hence, by Theorem \ref{Theorem 34} there is a set $X(\ve)$ with $V(\mu;\A;E;X(\ve))<\ve$, such that $\left(||f_j(x) - f(x)\right)_j$ is uniformly convergent to $0$ in $\setminus X(\ve)$. By linearity of the norm
\[
||f_j(x)\mu(I) - f(x) \mu(I)|| = ||f_j(x) -f(x)|| \mu(I) \leq 2F(x)\mu(I).
\]
Also $F(x) \chi\left(X(1/k,x)\right)$ is integrable. 

[GAP IN PROOF]

\noindent
Taking
\[
Y^k(I) = OS\left(\frac{\mu(I)}{k} + 2 \int_I F(x) \chi\left(X\left(\frac 1k\right), x\right)d\mu\right)
\]
and using Vitali's theorem to show that the latter integral is differentiable except in a set of variation zero, we have (\ref{66}), (\ref{67}), while (\ref{65}) follows from the absolute continuity of the integral (see M313) on taking $k$ large enough. If $K$ is a normed linear space and if
(\ref{65}), (\ref{66}), (\ref{67}) hold with $h_j = f_j \mu$, $h=f\mu$, then $f_j \rightarrow f$ almost everywhere. For, as in the proof of Theorem \ref{Theorem 34}, if
\[
X(\ve) = \left\{ x \in T: ||f_j(x)-f(x)|| \geq \ve \mbox{ for an infinity of } j\right\}
\]
then $V(\mu;\A;E;X(\ve))=0$ ($ \ve>0$).  \nproof

We may wish to extend the theory to convergence in measure.

If $||f_j|| \leq F$ we can take a sequence $(X_j)$ of sets such that (\ref{67}) for $h_j = f_j\mu$ need not hold while $x \in X_j$, and such that 
\[
V(\mu;\A;E;X_j) \rightarrow 0\;\;\mbox{ as }\;\;j\rightarrow \infty.
\]
The absolute continuity of the integral of $F$ seems to lead to Theorems \ref{Theorem 30}, \ref{Theorem 33}, \ref{Theorem 35} for this case.

\section{The Denjoy Extension}
In one dimension, an open set $G$ is the union of a sequence of disjoint open intervals $(a_j, b_j)$. When the integrals over the $(a_j,b_j)$ exist, possibly by the Cauchy extension, and when the Lebesgue integral over $[a,b]\setminus G$ (closed) exists, we can sometimes define the special Denjoy integral over $[a,b]$ to be the integral over $[a,b] \setminus G$, together with the infinite sum of the integrals over the separate $[a,b] \cap [a_j,b_j]$. This idea is applied transfinitely to give the Denjoy integral.

We turn to the corresponding theorem in generalized Riemann integration. Does not know whether the result holds in higher dimensions. The following example of open sets could not be dealt with by an analogue of the proof for one dimension: [CIRCLE: half, quarter, eighth, ... (halving each time). In the final quarter the proof fails.]

\begin{theorem}
\label{Theorem 36}
In Euclidean space of one dimension, we use intervals $[u,v)$ with associated point either $u$ or $v$. If $G = \bigcup_{j=1}^\infty (u_j, v_j)$ is an open set in a finite interval $(a,b)$, where the $(u_j,v_j)$ are disjont, if $h(I,x) \chi (\setminus G,x)$ is integrable over each $[u_j, v_j)$ ($j=1,2,\ldots $), and if, given $\ve>0$, there is an integer $J$ such that for every finite collection $\Q$ of intevals $[u,v)$, each contained in a $[u_j, v_j)$, for some $j \geq J$, no two intervals $[u,v)$ lying in the same $[u_j,v_j)$, we have
\begin{equation}
\label{77}
\left| (\Q)\sum \int_u^v \chi (G,x) dh\right| < \ve,
\end{equation}
then there exists
\begin{equation}
\label{78}
\int_a^b h = \int_a^b \chi(\setminus G,x)dh + \sum_{j=1}^\infty \int_{u_j}^{v_j} \chi (G,x) dh.
\end{equation}
\end{theorem}
\textbf{Proof:}
Subtracting the first integral on the right of (\ref{78}), from the left hand side, we have to prove that there exists
\begin{equation}
\label{79}
\int_a^b \chi(G,x) dh = \sum_{j=1}^\infty \int_{u_j}^{v_j} \chi(G,x) dh.
\end{equation}
We therefore define
\begin{equation}
\label{80}
H_1(E):=\sum_{j=1}^\infty \int_{[u_j,v_j) \cap E} \chi(G,x) dh\;\;\;\;\;\;(E \subseteq [a,b)\;)
\end{equation}
and our first task is to show that the series is convergent. As $E$ is an elementary set, $[u_j, v_j) \cap E$ is either empty or $[u_j, v_j)$, for all but a finite number of $j$ for which $(u_j,v_j)$ contains a boundary point of $E$. (This property fails in two dimensions, but the failure is not crucial.) Thus from (\ref{77}), the sequence of partial sums for (\ref{80}) is fundamental (Cauchy) and so convergent, and $H_1(E)$ exists in $[a,b)$. Next, there is an $\bS_j \in \A$, defined by $\delta_j(x)>0$, and dividing $[u_j,v_j)$, for which
$
\cS(\chi(G,x)h(I,x); \bS_j; [u_j, v_j))$
lies in the circle with centre $H_1(u_j,v_j)$ and radius $\ve 2^{-j}$, ($j=1,2, \ldots $). Hence $\vs(\chi(G,x)h(I,x) - H_1(I);\A;[u_j, v_j))$ lies in the circle with centre the origin and radius $\ve 2^{-j}$. From the separate $\bS_j$ we construct an $\bS \in \A$ defined by a function $\delta(x)>0$, and dividing $[a,b)$, in the following way. Let $J$ be an integer satisfying (\ref{77}). If $x \in G$ we have $u_j < x< v_j$ for some $j$, and we take the largest $\delta(x)>0$ such that
\begin{equation}
\label{81}
(x-\delta(x), x+\delta(x)) \subseteq (u_j, v_j),\;\;\;\;\;\delta(x) \leq \delta_j(x).
\end{equation}
On the other hand, if $x \in \setminus G$ we take the largest $\delta(x)>0$ such that
\begin{equation}
\label{82} (x-\delta(x),x) \subseteq (u_j, v_j),\;\;\;\delta(x) \leq \delta_j(x)\;\;\mbox{ when } x = v_j\;\;\mbox{ for some } j; \mbox{ \textbf{and}} \end{equation}
\begin{equation}
\label{83}
(x-\delta(x),x) \cap (u_j, v_j),\;\;\mbox{ is empty } \;(1\leq j \leq J-1) \;\; \mbox{ when }x \neq v_k \mbox{ for any } k, \end{equation} \begin{equation}
\label{84}
(x-\delta(x),x) \subseteq (u_j, v_j),\;\;\delta(x) \leq \delta_j(x) 
\mbox{ when } x=u_j \mbox{ for some } j
\end{equation} \begin{equation}
\label{85}
(x, x+\delta(x)) \cap (u_j, v_j),\;\;\mbox{ is empty } \;
(1\leq j \leq J-1) \;\; \mbox{ when }x \neq u_k \mbox{ for any } k. \end{equation}
Let $\D$ be a division of $[a,b0$ that is compatible with this $\delta(x)$. If $x$ is the associated point of $[u,v) \in \D$ then $x=u$ or $x=v$. By the construction (\ref{82}), if $x=v=v_j$ then $u_j \leq v <v_j$. If $x=v \neq v_j$, $x \in \setminus G$, then by (\ref{83}), $[u,v)$ can only overlap with $(u_j,v_j)$ when $j \geq J$. If $x=v$, $x \in G$, then $[u,v) \subseteq [u_j,v_j)$ for some $j$.
Similarly when $x=u$

\noindent
\textbf{[A]} Using (\ref{83}), (\ref{85}), the sum of $H_1(u,v)$, for $[u,v)$ with $x \in \setminus G$, $x=u \neq v_j$ or $x=v \neq v_j$ (all $j$), is the limit of a sequence of sums over various $\Q$ satisfying (\ref{77}), and the modulus of the sum is $\leq \ve$. Here $x\notin G$, so that
\[
H_1(I) - h(I,x) \chi(G,x) =H_(I),\;\;\;\;\;I=[u,v).
\]
\textbf{[B]} If $x=u=u_j$ or $x=v=v_j$ (when $x\in \setminus G$), or $x \in u_j,v_j)$, then by (\ref{81}), (\ref{82}), (\ref{84}),
\[
([u,v),x) \in \bS_j.
\]
Thus we see that
\[
[u_1,v_1), \ldots , [u_{J-1}, v_{J-1})
\]
are each divided by partial divisions of $\D$, using these $[u,v)$; and partial sets of $[u_J,v_J)$, $[u_{J+1}, v_{J+1})$, $\ldots$ are divided by such $[u,v) \in \D$. Using (\ref{80}), the results for 
\[
\cS\left(\chi(G,x)h(I,x);\bS_j;[u_j,v_j)\right),\;\;\;\;\;
\vs\left(\chi(G,x)h(I,x)-H_1(I);\A;[u_j,v_j)\right),
\]
and the preceding remarks \textbf{[A]} and \textbf{[B]},
$\left|(\D_x) \sum \chi(G,x)h(I,x) -H_1(a,b)\right|$=
\[
=\left|
(\D_x) \sum \left(\chi(G,x)h(I,x) -H_1(I)\right)\right|
\leq \ve + \sum_{j=1}^{J-1} \frac \ve{2^j}
+ \sum_{j=J}^{\infty} \frac \ve{2^j} < 3\ve,
\]
proving the result.  \nproof

An interval-point functin $h(I,x)$ is of \emph{bounded variation} (VB*) in $X$ relative to $E$ if $V(h(I,x);\A;E;X)$ is finite. We say that $h$ is of \emph{generalized bounded variation} (VBG*) in $X$, relative to $E$, if $X$ is the union of a sequence $(X_j)$ of sets with $V(h(I,x);\A;E;X_j)$ finite for all $j$.
\begin{theorem}
\label{Theorem 37}
If there is a positive function $k(x)>0$ for which \[V(k(x)h(I,x);\A;E;X)\] is finite, then $h$ is VBG* in $X$ relative to $E$, $(T,\T,\A)$ being a division space. If $(T,\T,\A)$ is a decomposable division space, the converse holds.
\end{theorem}
\textbf{Proof:}
If $V(k(x)h(I,x);\A;E;X)$ is finite, let $X_j$ be the set where
\[
\frac 1j \leq k(x) < \frac 1{j-1},\;\;\; j=2,3, \ldots .
\]
Then in $X$ we have
\begin{eqnarray*}
h(I,x) &=& \frac{k(x)h(I,x)}{k(x)} \;\;\leq\;\; jk(x)h(I,x), \vt
V(h;\A;E;X_j )&\leq & jV(kh;\A;E;X_j) \;\;\leq \;\; jV(kh;\A;E;X_j)
\end{eqnarray*}
which is finite.
Conversely, if $(T,\T,\A)$ is a decomposable division space, and if $V(h;\A;E;X_j) = a_j$ is finite, then give $\ve>0$ there is an $\bS_j \in \A$ that divides $E$ such that
\[
V(h;\bS_j;E;X_j) \leq a_j + \ve.
\]
We can assume that the sets $X_j$ are disjoint, for otherwise we can replace them by
\[
X_1, X-2 \setminus X_1, X_3 \setminus \left(X_1 \cup X_2\right), \ldots
\]
in each of which $h$ is still VB*. By decomposability there is an $\bS \in \A$ dividing $E$ such that 
\[
\bS[X_j] \subseteq \bS_j[X_j],\;\;\;\;j=1,2,3, \ldots.
\]
Also let us define
\[
k(x) = \left\{
\begin{array}{ll}
\left(a_j + \ve\right) 2^{-j}
& (x \in X_j,\;\;j=1,2,3, \ldots) , \vt
1 & (x \notin X).
\end{array}
\right.
\]
Then
\[
\begin{array}{rll}
V(kh;\A;E;X) &\leq & V(kh;\bS;E;X) \vt
&\leq & 
 \sum_{j=1}^\infty V(kh;\bS;E;X_j) \vt
 & \leq & \sum_{j=1}^\infty  \left(a_j + \ve\right) 2^{-j}   V(kh;\bS_j;E;X_j) \vt
 &\leq & \sum_{j=1}^\infty 2^{-j} \;\;\;=\;\;\;1,
 \end{array}
 \]
 giving the result.  \nproof
 
 \section{Cartesian Product of Two Spaces}
 From two spaces $T^x, T^y$ of points we construct the product space of points $z=(x,y)$:
 \[
 T^z = T^x \times T^y = \left\{(x,y) : x \in T^x, y\in T^y \right\}.
 \]
 We put the extra index $t$ on objects connected with $T^t$ ($t=x,y,z$). To obtain a division space for $T^z$, let $\T^t$, $\A^t$ be the respective families of all $t$-intervals $I^t \subseteq T^t$, and of $(I^t,t)$, ($t=x,y,z$), where $z$-intervals are the products $I^z=I^x\times I^y$. We can write
 \[
 \left(I^z,z\right) =  \left(I^x,x\right) \times  \left(I^y,y\right).
 \]
 \textbf{Stable Families:} $\A$ is \emph{stable} relative to the associated points if, for each elementary set $E$, and all $\bS \in \A$ that divide $E$,
$\left\{x:(I,x) \in \bS, I \subseteq E\right\}$ is independent of $\bS$, depends only on $E$, and is the set
\[
E^* = \left\{ x: (I,x) \in \T^x, I \subseteq E \right\}.
\]
($E^8$ will usually be the closure of $E$.) We use stable families $\A^t$ of $\bS^t$, these $\bS^t$ being sets of $(I^t,t)$, ($t=x,y,z$).

$\A^x, \A^y, \A^z$ have the \emph{Fubini property} in common if two properties\footnote{\textbf{Note by P.~Muldowney:} The notation here is a bit confusing. Subscripts $_x$, superscripts $^x$, and arguments $(x)$ have different roles in the definition, and the reader is expected to assign the appropriate role in each case. (Likewise, of course, in notation $\D_x = \{(I,x)\}$, the first $x$ indicates that the partition $\D$ has linked tag points, while the second $x$ indicates a member of the set of tag points; and a similar mental distinction has to be applied by the reader in that case also.)} hold:

First, let $E^x, E^y$ be arbitrary elementary sets and let $\bS^z$ be an arbitrary member of $\A^z$ that divides $E^z = E^x \times E^y$. Then, to each $x \in E^{*x}$ there is an $\bS^y(x) \in \A^y$ that divides $E^y$.  Further, to each collection of divisions $\D_y^y(x)$ of $E^y$ from $\bS^y(x)$, one division for each such $x$, there corresponds an $\bS^x \in \A^x$ that divides $E^x$, such that if
\begin{equation}
\label{86}
(I^x,x) \in \bS^x,\;\;\;(I^y,y) \in \D_y^y(x),\;\;\;\mbox{ then }\;\;\;(I^x,x) \times (I^y,y) \in \bS^z.
\end{equation}
As this property is unsymmetrical in $x$ and $y$, we assume also the property in which $x$ and $y$ are interchanged, keeping the product space as $T^x \times T^y$

If $\A^x, \A^y, \A^z$ have the Fubini property in common, and if the $(T^t, \T^t, \A^t)$ ($t=x,y,z$) are division spaces, we call $(T^z, \T^z, \A^z)$ a \emph{Fubini division space}.
\begin{theorem}
\label{Theorem 38}
Let $(T^z, \T^z, \A^z)$ be a Fubini division space\footnote{and let $(T^x,\T^x, \A^x)$ be decomposable.}.
Let the real- or complex-valued $h^x(I^x,x), h^y(x;I^y,y)$
be defined for $(I^t,t) \in \T^t$, ($t=x,y$), respectively, and for all $x \in E^{*x}$. If 
\begin{equation}
\label{87}
h^z(I^z,z) = h^x(I^x,x)  h^y(x;I^y,y), \;\;\mbox{ for }\;\; 
(I^z,z) \times (I^y,y) = (I^x,x) \times (I^y,y),
\end{equation}
has a generalized Riemann integral $H(E^z)$ in $E^z = E^x \times E^y$, then the function $J(x) = \int_{E^y} h^y(x;I^y, y)$ exists, except for the set of $x \in X^x$ with 
$V(h^x;\A^x;E^x;X^x) =0$, and if we put $J(x)=0$ in $X^x$, the generalized Riemann integral of $J(x) h^x(I^x,x)$ in $E^x$ is 
\[
H(E^z) = \int_{E^z} h^xh^y(x; \cdot) = \int_{E^x} \left( \int_{E^y} h^y(x; \cdot) \right) dh^x.
\]
\end{theorem}
\textbf{Proof:} (\emph{In earlier proofs it was assumed that $h^x$ is VBG*. Saks (\emph{Theory of the Integral}, second edition, pp.~87--88) gives three examples to show that, for Lebesgue integrals, the measure has to be VBG*. However, in the first two example the measure of an interval is $+\infty$, while the third example depends on a peculiarity of measure spaces not shared by division spaces.}) Denoting the integral of $h^z$ over $E_1^z$ by $H(E_1^z)$, given $\ve>0$, let $\bS^z \in \A^z$ and dividing $E^z$ be such that for all divisions $\D_z^z$ of $E^z$ from $\bS^z$,
\begin{equation}
\label{88}
(\D_z^z) \sum \left| h^z(I6z,z) - H(I^z)\right| < \ve.
\end{equation}
By the Fubini property nthere are suitable $\bS^y(x) \in \A^y$ dividing $E^y$, divisions $\D_y^y(x)$ of $E^y$ from $\bS^y(x)$, and an $\bS^x \in \A^x$, dividing $E^x$, that satisfy (\ref{86}). Let $\D_x^x$ be a division of $E^x$ from $\bS^x$. By (\ref{86}), (\ref{88}) and the finite additivity of $H$,
\begin{equation}
\label{89}
(\D_x^x) \sum \left| (\D_y^y)\sum h^x(I^x,x) h^y(x;I^y,y) - H(I^x \times E^y) \right| < \ve.
\end{equation}
Let $X_j^x$ be the set of $x$ where there are at least two divisions $\D_{y1}^y(x), \D_{y2}^y(x)$ of $E^y$ from each $\bS^y(x)$, for which
\begin{equation}
\label{90}
\left|(\D_{y1}^y(x))\sum h^y(x;I^y,y) -
(\D_{y2}^y(x))\sum h^y(x;I^y,y)\right| > \frac 1j.
\end{equation}
Then, from (\ref{89}), (\ref{90}), we have
\[
(\D_x^x) \sum  j^{-1} \left|h^x(I^x,x)\right| \chi(X_j^x,x) < 2\ve,\;\;\;\;\;\;V(h^x;\bS^x;E^x;X_j^x) < 2j\ve,
\]
\begin{equation}
\label{91}
V(h^x;\A^x;E^x;X_j^x)=0,\;\;\;\;\;\;V(h^x;\A^x;E^x;X^x)=0
\end{equation}
since $X^x$ is the countable union of sets $X_j^x$.
In $\setminus X^x$, $J(x)$ exists as an integral, and we make the following construction. Let $\bS_j^z$ be an $\bS^z$ for which $\ve = j^{-1}$ in (\ref{88}), and let $\bS_j^y(x)$ be an $\bS^y(x)$ from $\bS_j^z$, by using the Fubini property. As $J(x)$ exists, there are $\bS_{j1}^y(x) \in \A^y$, dividing $E^y$, such that
\begin{equation}
\label{92}
|f_j(x) - J(x)|<\frac 1j,\;\;\;\;\;\;f_j(x) := (\D_{yj}^y(x))\sum h^y(x;I^y,y),
\end{equation}
where $\D_{yj}^y(x)$ is an arbitrary division of $E^y$ from $\bS_{j1}^y(x)$. By taking intersections of $\bS_k^y(x) \cap \bS_{k1}^y(x)$ for $1\leq k \leq j$, if necessary, we can assume that the sequence $(\bS_j^y(x))_j$ is a monotone decreasing sequence of families with respect to $j$ that satisfy (\ref{92}). Then $\D_{yk}^y(x)$, chosen from $\bS_k^y(x)$, is also from $\bS_j^y(x)$ when $k \geq j$.
By (\ref{89}) again, for $\ve=j^{-1}$ and $\D_y^y(x) = \D_{yj}^y(x)$, we have
\[\begin{array}{c}
(\D_x^x)\sum \left|h^x(I^x,x)\right|\left|f_{j+1}(x) - f_j(x) \right| < \frac 2j,\vspace{10pt} \vt 
V \left(\left|f_{j+1}-f_j\right| h^x;\A^x;E^x; \setminus X^x\right) \leq \frac 2j.
\end{array}
\]
As $|f_{j+1}-f_j|$ is a $k(x)$ as for VBG* function, we see that $h^x$ is VBG*, relative to $E^x$, in the set $X_1^x$ contained in $\setminus X^x$, where $f_{j+1} \neq f_j$ for some $j$. In $\setminus (X^x \cup X_1^x)$,
\begin{equation}
\label{93}
f_j(x) = f_1(x),\;\;\;f_j(x) = J(x)\;\;\;\;\;\mbox{ all }\;\;j.
\end{equation}
We now split $T^x$ into disjoint sets $X^x$, $X_1^x$, $ \setminus(X^x \cup X_1^x)$.
By (\ref{91}) we can ignore $X^x$. For $x$ in the last set, we can in (\ref{89}) replace
\begin{equation}
\label{94}
(\D_y^y(x)) \sum h_x^x(I^x,x) h^y(x;I^y,y)\;\;\; \mbox{ by }\;\;\; h_x^x(I^x,x) J(x).
\end{equation}
The rest of the proof shows that we can do the same for $x\in X_1^x$, on changing $\ve = j^{-1}$ to $2.j^{-1}$. Now, $h^x$ is VBG* relative to $E^*$ in $X_1^x$. Thus, by Theorem \ref{Theorem 37} (\textbf{which requires $(T^x, \T^x,\A^x)$ decomposable)}, there are a $k_1(x)>0$ in $X_1^x$, and an $\bS_1^x \in \A^x$, dividing $E^x$, for which
\begin{equation}
\label{95}
V\left(k_1(x)h^x;\A^x;E^x;X_1^x\right)<1,\;\;\;\;\;\;
V\left(k_1(x)h^x;\bS_1^x;E^x;X_1^x\right)\leq 1.
\end{equation}
By (\ref{92}) we can choose $r=r(x)$ so high tat
\begin{equation}
\label{96}
\left|f_r(x) -J(x)\right| < \frac 1{jk_1(x)}.
\end{equation}
By the Fubini property there is an $\bS^X$ to satisfy (\ref{86}) with
\[
\D_y^y(x) = \left\{
\begin{array}{rl}
\D_{y\,j(x)}^y(x) & (x \in X_1^x),\vt
\mbox{arbitrary}&(x\in X^x,\vt
\D_{y1}^y(x) & (x\notin X^x \cup X_1^x).
\end{array}
\right.
\]
By (\ref{91}), (\ref{95}), (\ref{96}) we can make the replacement (\ref{94}) in (\ref{89}) and obtain
\[
(\D_x^x)\sum \left|J(x)h^x(I^x,x) - H(I^x \times E^x)\right| < \frac 2j,
\]
completing the proof.  \nproof
\begin{itemize}
\item
The original Fubini is got from:
$h^z(I^z,z) = f(x,y) h^x(I^x,x) h^y(I^y,y)$.
\item
A result of Cameron \& Martin \& Robbins \& Rodgers uses
\[
h^z(I^z,z) = f(y)h^x(I^x,x) F(x;I^y),
\]
$F$ being finitely additive for $I^y \subseteq E^y$ and each
fixe $x$, and $h^x(I^x,x) F(x;I^y)$ being integrable over $E^x$ for each $I^y \subseteq E^y$.
\item
Lee gives a Fubini theorem on more general $h^z(I^z,z)$.
\end{itemize}

\begin{theorem}
\label{Theorem 39}
et $(T^z, \T^z,\A^z)$ be a Fubini division space. Let the real- or complex-valued $h^z(I^z,z)$ be defined in $E^x \times E^y = E^z$. Let
\begin{equation}
\label{97}
H=\int_{E^z} h^z,\;\;\;\;\;\;G(I^x) = \int_{E^y} h^z\left(I^x\times I^y, (x,y)\right)
\end{equation}
exist for each $I^x \subseteq E^x$. Let $\mu(I^x)$ be of bounded variation in $E^x$. If, given $\ve>0$ there is an $\bS^y \in \A^y$, dividing $E^y$, such that (\textbf{Note:} ``\emph{uniformity}'')
\begin{equation}
\label{98}
\left|G(I^x) - s(I^x)\right| \leq \ve \mu(I^x),\;\;\;\;\;
s(I^x) := (\D_y^y)\sum h^z(I^x\times I^y, (x,y)),
\end{equation}
for all divisions $\D_y^y$ of $E^y$ from $\bS^y$, with fixed $I^x$, then
\[
H=\int_{E^x}\left(\int_{E^y} h^z\right) = \int_{E^x} G(I^x).
\]
\end{theorem}
\textbf{Proof:} Apply (\ref{98}) to  (\ref{89}).  \nproof

From families  $\A^t$ ($t=x,y$), we now construct $\bS^z$, and a suitable $\A^z$. Given an elementary set $E^t$, for each $z \in T^z$ let there be an $\bS_1^t(z) \in \A^t$ dividing $E^t$ ($t=x,y$). Then $\A^z$ is the family of all finite unions $\bS^z$ of families
\[
\left\{\left(I^x\times I^y,(x,y)\right): I^t \subseteq E^t, (I^t,t) \in \bS_1\left((x,y)\right), t=x,y\right\},
\]
the finite unions being taken over disjoint products $E^x \times E^y$. Such a space $(T^z,\T^z,\A^z)$ is called the \emph{product division space} of the $(T^t, \T^t, \A^t)$, ($t=x,y$).
\begin{theorem}
\label{Theorem 40}
The product division space of the fully decomposable division space $(T^t,\T^t, \A^t)$ with stable $\A^t$, ($t=x,y$), is a fully decomposable Fubini division space with stable $\A^z$.
\end{theorem}
\textbf{Proof:}
Let $\bS_1^t(z) \in \A^t$, dividing $E^t$ ($t=x,y$), give rise to $\bS^z$, connected with $E^x \times E^y$. For fixed $x \in E^{*x}$ and all $y \in E^{*y}$ (\emph{stability}) we take $\bS_1^y\left((x,y)\right)$. As $(T^y, \T^y, \A^y)$ is fully decomposable and $\A^y$ stable, there is an $\bS_2y(x)\in \A^y$, dividing $E^y$, such that
\[
(I^y,y) \in \bS_2^y(x),\;\;\;(I^x,x) \in \bS_1^x((x,y)) \;\;\;\Longrightarrow \;\;\;(I^x,x) \times (I^y,y) \in \bS^z
\]
since $\bS_2^y(x)[\{y\}] \subseteq  \bS_1^y((x,y))[\{y\}]$.
(The set $\{y\}$ is the singleton of $y$.) For each $x \in E^{*x}$ let $\D_y^y(x)$ be a division of $E^y$ from $\bS_2^y(x)$, and let $y_1, \ldots ,y_n$ be the associated points in $\D_y^y(x)$. As $\A^x$ is directed in the sense of  divisions there is an $\bS_2^x(x) \in \A^x$, dividing $E^x$, that lies in $\bS_1^x((x,y_j))$, ($1 \leq j \leq n$). As $\A^x$ has fully decomposable families there is an $\bS_3^x \in \A^x$, dividing $E^x$, that has the same property. Hence the Fubini property is true with $\bS^y(x), \bS^x$, respectively, replaced by $\bS_2^y(x)$ and $\bS_3^x$; and since we can interchange $x$ and $y$ in the proof. To show that $\bS^z$ divides $E^x \times E^y$, we have a division $\D_x$ of $E^x$ from $\bS_3^x$, and for each associated point $x$ in that division, a division $\D_y^y(x) $ of $E^y$ from $\bS^y(x)$. The products of the
\[
(I^x,x) \in \D_x^x,\;\;\;\;\;(I^y,y) \in \D_y^y(x),
\]
lie in $\bS^z$, and together form a division of $E^x \times E^y$. Similarly for finite union of disjoint products of elementary sets, so that $\A^z$ divides all elementary sets. Finally, $\A^z$ is directed in the sense of divisions, is stable, additive, and has the restriction property, as $\A^x, \A^y$ have these properties; while $\A^z$ is fully decomposable by construction, and
\begin{equation}
\label{99}
\left(E^x \times E^y\right)^* = E^{*x}\times E^{*y}.
\end{equation}
\section{Division Spaces for \\ Cartesian Product Spaces}
Let a space $T(b)$ correspond to each $b$ of an index set $B$, and let $T_B$ be the union of the $T(b)$ for all $b \in B$. Then the \emph{Cartesian product space} $T = \mbX_B T(b)$ is the set of all functions $f:B\mapsto T_B$ with $f(b) \in T(b)$, $b\in B$, and
\[
\prod\left(X(b):b\in B\right) :=\mbX_B X(b)
\]
denotes the Cartesian product space contained in $T$ with
\[
X(b) \subseteq T(b)\;\;\;(b \in B_1),\;\;\;\;\;X(b) = T(b)\;\;\;(b \notin B_1),
\]
where $B_1 \subseteq B$. Further, we write
\[
T(B_1) := \mbX_{B_1} T(b),\;\;\;\;\;X(B_1) :=\mbX_{B_1}X(b),\;\;\;(B_1 \subseteq B).
\]
For each $b \in B$, let $(T(b), \T(b), \A(b))$ be a division space in which
\begin{equation}
\label{100}
T(b) \in \T(b),
\end{equation}
\begin{equation}
\label{101}
\mbox{each } I(b) \in \T(b) \mbox{ is a partial division of } T(b),
\end{equation}
\begin{equation}
\label{102}
\A(b) \mbox{ is stable}.
\end{equation}
When $B_1$ is finite, $B_1 \subseteq B$, we can obtain analogues of the preceding results for products of two spaces; in particular, we can construct families $\bS(B_1), \A(B_1)$ like $\bS^z, \A^z$, respectively.

If $E(b)$ is an elementary set for each $b \in B_1$, (\ref{99}) and (\ref{102}) give
\begin{equation}
\label{103}
E(B_1)^* = E^*(B_1).
\end{equation}
For $T$ we consider the family $\T$ of \emph{generalized) intervals} $I=\prod\left(I(b):b\in B_1\right)$ for all finite sets $B_1 \subseteq B$ and all $I(b) \in \T(b)$, ($b \in B_1$).

To construct $\bS$ dividing an elementary set $E$ of $T$ we associate with each $f \in T$ a finite set $B_1(\bS,f) \subset B$ and, for each finite set $B_2 \subseteq B$, an $\bS_1(B_2,\bS,f) \in \A(B_2)$, dividing $T(B_2)$.

We also suppose given a sequence $B_0=(b_j) \subseteq B$ that depends on $\bS$, with $B_1(\bS,f) \subseteq B_0$, (all $f$). 

Then $\bS$ is the family of all $(I,f)$ with
\begin{equation}
\label{104}
I=\prod(I(b): b\in B_2) \subseteq E,\;\;\;\;\;
B_2\mbox{ finite, }\;\;\;B \supseteq B_2\supseteq B_1(\bS,f),
\end{equation}
\begin{equation}
\label{105}
\mbX_{B_2}\left(I(b), f(b)\right) \in \bS_1(B_2,\bS,f),
\end{equation}
\begin{equation}
\label{106}
f(b) \in T(b)^*,\;\;\;(b \in B \setminus B_2).
\end{equation}
If $\A$ is the family of such $\bS$, for all sequences $B_0 \subseteq B$, all $B_1(\bS,f) \subseteq B_0$, all $\bS_1(B_2,\bS,f) \in \A(B_2)$ dividing $T(B_2)$ (all finite $B_2$ in $B \supseteq B_2 \supseteq B_1(\bS,f)$, then $(T,\T,\A)$ is the product division space of the division spaces $(T(b), \T(b), \A(b))$, ($b \in B$).

By construction, $\A$ has decomposable families that are not fully decomposable only because of the use of $B_0$.

\noindent
\textbf{Lemma:}
If $B_3$ is finite and $I=\prod\left(I(b):b\in B_3\right)$ 
then $I^* = I^*(B)$.

\noindent
\textbf{Proof:}
Let $f \in I^*(B)$. Then $f(b) \in I^*(b)$ for all $b \in B$. Because of (\ref{106}) we need only consider 
\[
b \in B_2 \supseteq B_3\cup B_1(\bS,f)
\]
We can choose $J(b) \supseteq I(b)$, $J(b) \in \T(b)$, so that
\[
\mbX_{B_2}\left(J(b), f(b)\right) \in \bS_1(B_2,\bS,f),
\]
by (\ref{102} to (\ref{105}). Then
\[
\left(\prod\left(J(b):b\in B_2\right),f\right) \in \bS,\;\;\;\;\;
\prod\left(J(b):b\in B_2\right) 
\subseteq \prod\left(I(b):b \in B_3\right)
\]
Hence $f \in I^*$, the set of all $f$ satisfying the latter.  
Conversely, if $f \in I^*$, there is a \[
\prod\left(J(b):b \in B_2\right) \subseteq \prod\left(J(b):b \in B_3\right),\;\;\;\;\;
\left(\prod\left(J(b):b \in B_2\right),f\right) \in \bS.
\]
By (\ref{103}), (\ref{105}), (\ref{106}) in the construction of $\bS$,
\[
f(b) \in J(b)^* \subseteq I(b)^*,\;\;\;(b \in B),\;\;\;f \in I^*(B),
\]
as required.   \nproof

If intervals $I_1, \ldots , I_m$ are disjoint we call the union $F$ of the $I_j^*$ ($1 \leq j \leq m$) an \emph{elementary *-set}. If $E$ is the union of the $I_j$ ($1\leq j\leq m$) then $E^* \supseteq F$, but we need not have equality. 

A family $\mcB$ of sets has the \emph{finite intersection property} if every finite subfamily of $\mcB$ has a non-empty intersection. An elementary set $E$ has the \emph{*-intersection property} if for every family $\mcB$ of elementary *-sets constructed from intervals in $E$ with the finite intersection property, the intersection of all $F \in \mcB$ is non-empty.

\begin{theorem}
\label{Theorem 41}
If $T(b)$ has the *-intersection property for each $b \in B$, then so has $T$.
\end{theorem}
\textbf{Proof:}
Let $\mcB$ be a family of elementary *-sets with the finite intersection property, and let $\mathcal{F}$ be the family of all elementary *-sets not in $\mcB$. By the well-ordering principle we can well-order $\mcB$ and $\mathcal{F}$, putting the sets of $\mathcal{F}$ after those of $\mcB$. Let $\leq$ denote the well-ordering relation. We construct a family $\mathcal{C}$ of elementary *-sets with the following property. An elementary *-set $H$ is in $\mathcal{C}$ if and only if $H$, the sets of $\mcB$, and all $K \in \mathcal{C}$ with $K \leq H$, together have the finite intersection property. Clearly
\begin{equation}
\label{107}
\mcB \subseteq \mathcal{C}
\end{equation}
\begin{equation}
\label{108}
\mbox{Further, if }H=\bigcup_{j=1}^m I_j^* \in \mathcal{C}, \mbox{ then, for some }j,\;\;\;I_j^* \in \mathcal{C}.
\end{equation}
For if $I_j^* \notin \mathcal{C}$, there are $H_1^j, \ldots , H_{p_j}^j$ in $\mathcal{C}$ and before $I_j^*$, such that 
\[
H_1^j \cap H_2^j \cap \cdots \cap H_{p_j}^j \cap I_j^*
\]
is empty. If true for $1\leq j\leq m$, there is a 
\[
K:=\bigcap \left\{H_k^j: 1 \leq j \leq m,\;\;1\leq k \leq p_j\right\}
\]
such that $K \cap I_q^*$ is empty. But of the sets $H$, $H_k^j$ ($1\leq j\leq m,\;\;1\leq k\leq p_j$) of $\mathcal{C}$, there is a last set, say $L$, in the well-ordering, and as $L \in \mathcal{C}$ then $K \cap H$ cannot be empty. Hence (\ref{108}). By (\ref{108}) we can take $\mathcal{E} \subseteq \mathcal{C}$ such that if $H \in \mathcal{E}$ then $H=I^*$ for some interval $I$, and conversely; and $\mathcal{E}$ has the finite intersection property. If $b$ is fixed in $B$ ten by the Lemma the set of $I(b)^*$, for all $I^* \in \mathcal{E}$, has the finite intersection property. As $T(b)$ has the *-intersection property, the intersection of the $I(b)^*$ is not empty (?\emph{``projection''}?), and so contains a point, say $g(b)$. Being true for each $b \in B$, we can define a function $g: B\mapsto T_B$ such that $g$ takes the value $g(b)$ for $b \in B$. For each $I^* \in \mathcal{E}$ there is a finite $B_1 \subseteq B$ with
\[
I = \prod\left( I(b): b \in B_1\right).
\]
As $g(b) \in I(b)^*$ ($b \in B$), the Lemma gives $g \in I^*$ for each $I^* \in \mathcal{E}$. Thus $g \in H$ for each $H \in \mathcal{C}$ and so for each $H \in \mathcal{B}$, by (\ref{107}). Hence the result.   \nproof

We now show that an $\bS$ defined for an elementary set $E$ divides $E$ when certain simple conditions are satisfied. It is sufficient to consider the intervals $J \subseteq E$ with
\[
J=J(B) =\prod\left(J(b): b \in B\right).
\]
\begin{theorem}
\label{Theorem 42}
Let $T(b)$ have the *-intersection property for each $b \in B$.\vspace{-10pt}
\begin{equation}
\label{109}
\end{equation}
Let there be a sequence $B_0=(b_j)_j$ of points $b \in B$, depending only on $\bS$, such that $B_3$ consists of $b_1, \ldots , b_m$, say, and that there is a least intege $p=p(\bS,f)\geq m$, such that $B_1(\bS,f)$ is a subset of $b_1, \ldots ,b_p$, for all $f \in T$.\vspace{-10pt}
\begin{equation}
\label{110}
\end{equation}
Let a sequence $(\D_n(b))_n$ of divisions of $T(b)$ exist with the properties:

\noindent
If $I(b) \in \D_1(b)$ then either $I(b) \subseteq J(b)$ or $I(b) \cap J(b)$ is empty.
\vspace{-10pt}
\begin{equation}
\label{111}
\end{equation}
$\D_{n+1}(b) \leq \D_n(b)$ (i.e.~each $K(b) \in \D_{n+1}(b)$ lies in an $I(b) \in \D_n(b)$ ($n \geq 1$).
\vspace{-10pt}
\begin{equation}
\label{112}
\end{equation}
For each $f \in T$ and each $B_2$ with $b_j \in B_1$, ($1\leq j\leq p(\bS,f)$), there is an integer $N=N(\bS,f,B_2)$ such that, if $n \geq N$, $I(b) \in \D_n(b)$, $I(b) \subseteq J(b)$, $f(b) \in I(b)^*$, then $I(b)$ is an interval $L(b)$ or a finite union of disjoint intervals $L(b)$, ($b \in B_2$), with
\begin{equation}
\label{113}
\mbX_{B_2} \left(L(b), f(b)\right) \in \bS_1(B_2, \bS,f).
\end{equation}
Then $\bS$ divides $J$.
\end{theorem}
\textbf{Proof:}
We suppose the theorem false, and use (114) repeatedly.

\noindent
\emph{\textbf{(114):
If $I \in \T$ is  divided by $bS$, and if $\D$  is a division of $I$, then there is an interval $K \in \D$ such that $K$ is not divided by $\bS$. (Because the union of divisions from $\bS$ of the $K \in \D$ would be a division of $I$.)} }

\noindent
Let $B \supseteq B_4 \supseteq B_3$, $B_4$ finite.
For all $I(b) \in \D_n(b)$ ($b \in B_4$), the $\prod(I(b): b\in B_4)$ form a division $\D_n(B_4)$ of $T$. By (\ref{112}), $\D_{n+1}(B_4) \subseteq \D_n(B_4)$, and by (\ref{111}) each $I \in \D_n(B_4)$ is either in $J$ or is disjoint from $J$. By (114), as $J$ is not divided by $\bS$, there is an $I \in \D_1(B_4)$ not divided by $\bS$, with $I \subseteq J$. If $I \subseteq J$, $I \in \D_n(B_4)$ for some $n \geq 1$ is a $K \in \D_{n+1}(B_4)$ not divided by $\bS$, with $K \subseteq I \subseteq J$. If $E_n(B_4)$ is the union of all $I \in \D_n(B_4)$ not divided by $\bS$, and if 
\[\vspace{-5pt}
E^+_p(B_4) = \bigcap_{n=1}^p E_n(B_4)
\vspace{-3pt}
\]
it follows that each $E^+_p(B_4)$ is a non-empty elementary set. Similarly
\[\vspace{-5pt}
E_{pq} = \bigcap_{j=m}^q E^+_p\left((b_1, \ldots , b_j)\right)
\vspace{-3pt}
\]
is a non-empty elementary set. Hence $F_{pq}$, the union of the $I^*$ for $I \subseteq E_{pq}$, is a non-empty elementary *-set. By construction, any finite number of the $F_{pq}$ have a non-empty intersection. By (\ref{109}) and Theorem \ref{Theorem 41}, there is an $f \in F_{pq}$ for all $p,q$. Let $B_2 = (b_1, \ldots ,b_q)$ for $q = p(\bS,f)$ and let $N$ be the integer of (\ref{113}) corresponding to $\bS$, $f$, $B_2$. Then $f \in F_{Nq}$, and yet by (\ref{113}) the $\prod(I(b):b\in B_2)$ constructed from the $\D_N(b)$, with $f(b) \in I(b)^*$,  is divided by $\bS$. This contradiction proves the theorem.  \nproof

\noindent
\textbf{Note by P.~Muldowney:} \emph{This proof does not work; see citation [107] of MTRV page 512, \cite{MTRV}. A correct method of proof is given in Theorem 4, pages 121--124 of MTRV, \cite{MTRV}.}

\begin{theorem}
\label{Theorem 43}
The $(T,\T,\A)$ of Theorem \ref{Theorem 42} is a decomposable division space.
\end{theorem}
\textbf{Proof:}
We have already seen that it is decomposable by construction. By Theorem \ref{Theorem 42}, if $\bS_2, \bS_3 \in \A$ both divide an elementary set $E$, there is an $\bS_4 \in \A$ that also divides $E$, with $\bS_4 \subseteq \bS_2 \cap \bS_3$. For we need only put
\[
\begin{array}{rll}
B_1(\bS_4,f)& =& B_1(\bS_2,f) \cup B_1(\bS_3,f),\vt
\bS_1(B_2,\bS_4,f)& \subseteq &\bS_1(B_2, \bS_2,f) \cap \bS_1(B_2, \bS_3,f),
\end{array}
\]
and take the $(b_j)$ for $\bS_4$ to be such that $(b_{2j})$
is the sequence for $\bS_2$, and $(b_{2j+1})$ the sequence for $\bS_3$. Hence $\A$ is directed in the sense of divisions, by (\ref{104}), (\ref{105}) (\ref{106}). To show that $\A$ is additive we use a similar proof, since the $\bS_1(B_2,\bS,f)$ divide $T(B_2)$. Each point $f$ might be in several $I_1^*, \ldots , I_k^*$, but as these are only finite in number we need only take the unions of the separate $B_1(\bS,f)$ to obtain a suitable new $B_1(\bS,f)$, and by direction in the sense of divisions, we obtain a suitable $\bS_1(B_2, \bS,f)$. As each $\bS_1(B_2, \bS, f)$ divides $T(B_2)$, (\ref{104}) and Theorem \ref{Theorem 42}) show that $\A$ has the restriction property. Thus $(T,\T,\A)$ has the required properties. \nproof

\begin{theorem}
\label{Theorem 44}
Let the conditions of Theorem \ref{Theorem 42} hold, let $B = H\cup K$ with $H \cap K$ empty, and let $B$ be ordered in such a way that 
\[
U = \mbX_H T(b),\;\;\;\;V=\mbX_K T(b),\;\;\;\; U \times V =T.
\]
If the $\bS_1(B_2,\bS,f)$ are constructed from product division spaces using $\bS_1^b(f)$ for $b \in B_2$, then Fubini's theorem holds.
\end{theorem}
\textbf{Proof:}
Let $f=(u,v)$ with $u \in U, v\in V$, and first fix $u$. For various $v \in V$, there are intervals $I_U \subseteq U$, $I_V \subseteq V$, with $(I_U,u) \times (I_V,v) \in \bS$ where $\bS$ is as in Theorem \ref{Theorem 42}. Here, for some finite $B_2 \supseteq B_1(\bS(u,v))$,
\[
I_U = \prod_H \left(I(b) :b \in B_2 \cap H\right),\;\;
\left(I(b), u(b)\right) \in \bS_1^b(f) = \bS_1^b((u,v)),\;\;(b \in B_2 \cap H),
\]
\vspace{-15pt}
\[
I_V = \prod_K \left(I(b): b \in B_2 \cap K\right),\;\;
\left(I(b), u(b)\right) \in \bS_1^b(f) = \bS_1^b((u,v)),\;\;(b \in B_2 \cap K).
\]
(For each $(u,v)$ we can suppose that $B_1(\bS((u,v)))$ contains certain points of $H$ and of $K$.)
We have written $\prod_H, \prod_K$ for products over $ b\in H$, $b \in K$, respectively. The $(I_V,v)$ form an $\bS_V$ for $V$ that is like $\bS$, except that $T,B$ are replaced by $V,K$, respectively. Hence, by Theorem \ref{Theorem 42}, there is a division $\D(u)$ of $V$ formed by intervals $J_V \subseteq V$ with associated points $y$, say. The corresponding finite union of the finite subsets $B_1(\bS, (u,y))$ gives a finite subset $B_3(\bS,u)$, while there is an $\bS_2^b(u) \in \A(b)$ that divides $T(b)$ and lies in the $\bigcap \bS_1^b(u,y)$ for the finite number of $y$. Hence there is an $\bS_U$ of pairs $(I_U,u)$ with the property that that for some finite set $B_4 \supseteq B_3(\bS,u)$,
\[
I_U = \prod_H \left(I(b): b \in B_4 \cap H\right),\;\;\;\;\;
(I(b), u(b)) \in \bS_2^b(u),\;\;\;(b \in B_4\cap H).
\]
This $\bS_U$ for $U$ is like $\bS$ except that $T,B$ are replaced by $U,H$ respectively, and $\bS_1$ by $\bS_2$.
Hence, by Theorem \ref{Theorem 42}, $\bS_U$ divides $U$, and the first Fubini property is true. Similarly for the second Fubini property; and Theorem \ref{Theorem 38} shows that the present theorem is true.  \nproof

\section{The Special Case}
The origins of integration in function space occur in A.~Einstein, ``\emph{Zur Theorie der Brownschen Bewegung}'', Ann.~d.~Physik 19(1906) 371--381; and M.~Smoluchowski, Ann.~d.~Physik 48(1915) 1103. For example, see S.~Chandrasekhar, ``\emph{Stochastic problems in physics and astronomy}'', Dover, New York (1954).

We consider the event that a free, spherical Brownian particle, starting at time $b=1$ from position $x=0$, will satisfy at time $b=b_j$,
\[
u(b_j) \leq x<v(b_j),\;\;\;\;(1 \leq j \leq n),
\]
with $u(b_j) < v(b_j)$, $a<b_1<b_2 < \cdots <b_n \leq c$, $(a,c] =$ set $B$.

If $B_5 = (b_1, \ldots , b_n)$ the probability of the event is taken to be
\[
P(I) = \int_{u(b_1)}^{v(b_1)} \cdots \int_{u(b_n)}^{v(b_n)} p(x_1, \ldots , x_n;B) dx_1 \ldots dx_n
\]
where 
\[
I = \prod \left([u(b_j), v(b_j)):b \in B_5\right),\] \vspace{-15pt}
\[
p(x_1, \ldots , x_n;B_1) := p(x_1;b_1-a|0)p(x_2;b_2-b_1|x_1) \cdots p(x_n;b_n-b_{n-1}|x_{n-1}).
\]
In the Brownian case
$p(y;t|x)=g(y-x;4Dt)$ where
\[
g(x;t) = \frac 1{\sqrt{\pi t}} e^{-\frac{x^2}t},\;\;\;\;\;
\mbox{i.e.}\;\;\;p(y;t|x) = \frac 1{\sqrt{4\pi Dt}} e^{-\frac{(y-x)^2}{4Dt}}.
\]
Wiener integration uses $D=\frac 12$; and the  
functions $f \in T$ satisfy
\[
A|y-x|^{\frac 12 + \ve} \leq |f(y)-f(x)| \leq B|y-x|^{\frac 12 - \ve},
\]
i.e.~continuous but not differentiable.

Feynman uses complex-valued $g$ with $D= \frac \iota 4 = \frac{\sqrt{-1}}4$. The Feynman measure is complex and not even VBG*. (\textbf{Note by P.~Muldowney:} \emph{See MTRV Section 6.7 pages 282--284, \cite{MTRV}.}) 

Cameron and his pupils define Feynman integration as the limit as $\sigma \rightarrow 0+$ of the Lebesgue-type integration got using $D = \sigma + \iota/4$.

 \vspace{20pt}
 
 \noindent
 ****

 \vspace{20pt}
 
 For functionals $F$ of $f$, we integrate $h(I,f):=F(f)P(I)$ over $T$. The $p(x_1, \ldots ,x_n;B_5)$ is taken to be continuous in $x_1, \ldots ,x_n$ for fixed $B_5$. Then we can replace $P(I)$ by 
 \[
 p(f(b_1), \ldots , f(b_n);B_5) \prod_{j=1}^n \left(v(b_j) - u(b_j)\right).
 \]
\begin{theorem}
\label{Theorem 45}
Let $K$ be the integral of the bounded $F$ in $T$ with respect to a finitely additive 
\[
P(I) = \prod_{j=1}^n P_j\left(I(b_j)\right)\;\;\mbox{ where }\;\; I=\prod \left(I(b): b=b_j, 1 \leq j \leq n\right)
\]
Given $\ve>0$, let $\bS$ in Theorem \ref{Theorem 42} (?) be such that
\begin{equation}
\label{115}
(\D)\sum |K-FP| < \ve.
\end{equation}
Let $C(q) = (b_1, \ldots , b_q)$, $\D_{nq} = \D_n(C(q))$, $=\D_n(b_1) \times \cdots \times \D_n(b_q) \times \mbX_{b\notin C(q)} T(b)$. Let $p(\bS,f)$, $N(\bS, f, C(q))$ be the functions of Theorem \ref{Theorem 42}. Let $Y_q$ be the set of $f\in T$ with $p(\bS, f) \leq q$, and let $Y_{qn}$ be the set of $f\in Y_q$ with $N(\bS,f,C(q)) \leq n$. If ech $Y_{qn}$ is measurable with
\begin{equation}
\label{116}
V(P;\A;T;I^*) \leq V(P;\A;I) < \infty,\;\;\;(I\in \T),
\end{equation}
then
\begin{equation}
\label{117}
\lim_{q\rightarrow \infty} \int_{T(C(q))} F(f_8+f^*) dP = \int_T F(f)dP
\end{equation}
for almost all $f^* \in T(C^*(q))$, $C^*(q) :=B\setminus C(q)$, where
\[
f_*(b) = \left\{
\begin{array}{rl}
f(b),& b\in C(q),\vt
0, & b \in C^*(q),
\end{array}\right.
f^*(b) = f(b) - f_*(b).
\]
\end{theorem}
\textbf{Proof:}
As $Y_q$ is monotone increasing in $q$ to $T$, and $Y_{qn}$ monotone increasing in $n$ to $Y_q$, and as $(T,\T,\A)$ is a decomposable division space, we have
\[
\lim_{q \rightarrow \infty}V(P;\A;T;Y_q) = V(P;\A;T) < \infty,\] \vspace{-10pt}\[
\lim_{n \rightarrow \infty}V(P;\A;T;Y_{qn}) = V(P;\A;T;Y_q).
\]
Thus for $\ve>0$ there are a $q=g(\ve)$ and an $m=m(q,\ve)=m(\ve)$, with
\begin{equation}
\label{118}
\left\{ \begin{array}{rll}
V(P;\A;T;Y_{qn}) &>& V(P;\A;T) - \ve^2,\vt
V(P;\A;T;T\setminus Y_{qn})&<& \ve^2\;\;\;\;\;(n \geq m).
\end{array} \right.
\end{equation}
The latter inequality is by measurability of $Y_{qn}$ and the expression\footnote{\emph{If $h$ is integrable to $H$ and of bounded variation, then $|h|$ and $|H|$ are integrable to $V(h;\A;\cdot)$.}} of the variation as an integral.  If $Z(C(q),f^*)$ is the set of all $f_*$ with $f_*+f^* \in Z$ (some given set), then by Fubini's theorem (using $P(I) = \prod P_j(I(b_j))$) 
\begin{equation}
\label{119}
\begin{array}{r}
V(P;\A(C(q));T(C(q));\setminus Y_{qn}(C(q),f^*))< \ve \vt
\mbox{except in a set } X \subseteq T(C^*(q)) \mbox{ with}
\end{array}
\end{equation}
\vspace{-10pt}
\begin{equation}
\label{129}
V(P;\A(C^*(q));T(C^*(q));X) < \ve.
\end{equation}
For fixed $f^* \in X$ let $\Q$ be the set of all $I \in \D_{mq}$ such that there is an $f_* \in T(C(q))$ with $f=f_*+f^* \in Y_{qm} \cap I^*$. Then $f$ is an extra associated point of $I$. Let $R$ be the union of the $I^*$ with $I \in \Q$. Then
\begin{equation}
\label{121}
R(C(q),f^*) \supseteq Y_{qm}(C(q),f^*).
\end{equation}
As $\D_{mq}$ only uses the $b\in C(q)$,
\[
\begin{array}{rl}
&V(P;\A(C(q));T(C(q)); Y_{qm}(C(q), f^*) \vt
=& V(P;\A;T; Y_{qm}(C(q), f^*)\times T(C^*(q))) \vt
\leq & V(P;\A;T;I^*) \;\; \leq \;\;(\Q)\sum V(P;\A;I). 
\vspace{10pt}
\vt
& (\D_{mq}\setminus \Q) \sum V(P;\A;I) \;\;=\;\;
(\D_{mq})\sum V - (\Q)\sum V \vt
\leq & V(P;\A;T) - V(P;\A(C(q));T(C(q));Y_{qm}(C(q,f^*)) \vt
=& V(P;\A(C(q));T(C(q)); \setminus Y_{qm} (C(q), f^*)) \;\;<\;\;\ve.
\end{array}
\]
If $|F|\leq M$, then, as $P$ is finitely additive, so that $|P|\leq V$, we have
\[
\begin{array}{rllll}
(\D_{mq}\setminus \Q) \sum |K|& \leq & M (\D_{mq}\setminus \Q) \sum V &<& M\ve,\vt
(\D_{mq}\setminus \Q) \sum |FP| &\leq & M (\D_{mq}\setminus \Q) \sum V &<& M\ve.
\end{array}
\]
As $\Q$ is from $\bS$, though $\D_{mq}$ need not be, we have
\[
(\D_{mq})\sum |K-FP| \leq (\Q)\sum  + (\D_{mq} \setminus \Q) \sum < \ve +2M\ve.
\]
By Fubini's theorem the first integral (= l.h.s.~of (\ref{117})), call it $K(C(q), F^*)$, exists for almost all $f^*$. Also, we can choose the $\bS_1$ to satisfy
\[
\bS_1(C(q);\bS;f) \subseteq \bS^*(C(q);\bS;F^*)
\]
where $\bS^*$ is such that, for all divisions $\D$ of $T(C(q))$ from $\bS^*$,
\[
(\D)\sum |K(C(q),f^*)(I) - F(f_*+f^*)P(I)| < \ve.
\]
By a similar proof we have
\[
|K(C(q),f^*) -K(T)| < 2\ve + 4M\ve,
\]
giving the result outside $X$. \nproof

\vspace{20pt}

\noindent
******
\vspace{20pt}

If $P(J)$ ($=\int_j p\Delta x$) exists for some interval $J$ in $T$ we can replace $P(I)$ by
\begin{equation}
\label{122}
p(x_1, \ldots ,x_n;C)\Delta x_1 \cdots \Delta x_n,
\end{equation}
where $x_j = f(b_j),\;\;\;$ $\Delta x_j = v(b_j)-u(b_j)$ and 
\[
I=\prod \left([u(b_j), v(b_j)):b = bj \in C\right),\;\;\;\;\;\;C=(b_1, \ldots ,b_n)
\]
for some $\bS(C) \in \A(C)$ dividing $J(C)$ and for all $(I(C), (x_1, \ldots ,x_n))\in \bS(C)$ with $I \subseteq J$. But in general there is no $g$ with \[(J(C),(g(b_1), \ldots ,g(b_n))) \in \bS(C)\]
so that we cannot always take $I=J$. 

The continuity of $p$ bridges the gap. Fixing $C$, the continuity and positiveness of the Wiener density function
\[
w(x_1, \ldots ,x_n;C) = \frac 1
{(2\pi)^{\frac 12 n} \prod_{j=1}^n(b_j-a_j)^{\frac 12}}
\prod_{j=1}^n e^{
-\frac 12 \frac{x_j-x_{j-1})^2}{b_j-b_{j-1}}}
\]
in $x_1, \ldots ,x_n$ implies that there is in $T(C)$ a nbd $N$ of each $(x_1, \ldots ,x_n)$ in which
\[
w(y_1, \ldots , y_n) > \frac 12 w(x_1, \ldots , x_n) >0\;\;\;\;\;\;\left((y_1, \ldots ,y_n) \in N\right).
\]
By continuity of $p(x_1, \ldots ,x_n;C)$ in $(x_1, \ldots ,x_n)$, given $\ve >0$, there is in $T(C)$ a nbd $N^* \subseteq N$ of $(x_1, \ldots , x_n)$ in which
\[
\left|p(y_1, \ldots ,y_n;C) - p(x_1, \ldots ,x_n;C) \right|
<\frac 12 \ve w(x_1, \ldots , x_n;C)\]
for $(y_1, \ldots ,y_n) \in N^*$.
There are intervals $I$ of $T$ with 
\[
I(C) \subseteq N^*,\;\;\;\;\;I=\prod\left([u^*(b_j), v^*(b_j);C\right),\;\;\;\;\;u^*(b_j) = u_j, v^*(b_j) = v_j,
\]
and with $f(b_j)=x_j, =u^*(b_j)$ or $v^*(b_j)$, $1\leq j\leq n$, with a suitable $(I,f)$ when $f(b_j) = \pm\infty$ for some $j$. Such $(I,f)$ form an $\bS\in \A$, and
\[
\begin{array}{rl}
& \left|P(I) - p(x_1, \ldots ,y_n;C)\Delta x_1 \ldots \Delta x_n\right| \vt
\leq & \int_I \left|p(y_1, \ldots , y_n;C) - p(x_1, \ldots , x_n;C)\right| dy_1 \ldots dy_n \vt
\leq & \frac 12 \ve w(x_1, \ldots ,x_n;C)\Delta x_1 \ldots \Delta x_n \vt
\leq & \ve \inf w(y_1, \ldots ,y_n;C)\Delta x_1 \ldots \Delta x_n \vt
\leq & W(I),
\end{array}
\]
$W(I)$ being Wiener measure $\int w(I)$, the pieces adding up to $\leq 1$.

This being possible for all $C \subseteq B$ and all $\ve >0$, it follows that (\ref{122}) is variationally equivalent to $P(I)$, and there is a theorem to show that (\ref{122}) can replace $P(I)$ in the integration. 

Then Fubini's theorem will go through for all collections $C$ of points of $B$ if $p(x_1, \ldots , x_n;C) = \prod_{j=1}^n p(x_j, b_j)$.

\noindent
\textbf{Note by P.~Muldowney:}\\ \emph{Henstock commented at this point that the Wiener integrator $w=p(x_1, \ldots , x_n;C)$ involves factors $x_j - x_{j-1}$, and does not resolve into a product of terms $\prod_{j=1}^n p(x_j, b_j)$, so Fubini's theorem does not work as he envisaged it. Actually, this version of Fubini's theorem turns out to be quite serviceable in this subject area. 
\\
At various points Henstock also suggested that the Feynman integrator is not VBG*. If this were strictly true it is hard to see how routine mathematical analysis could be applied to it. But the Feynman integrator is essentially based on ``very nice'' functions such as $\sin$ and $\cos$. This issue is addressed in MTRV (\cite{MTRV} pages 282, 283, and 288).
\\ These notes conclude with the following remarks showing that his concept of the limit theorems (``integral of limit equals limit of integrals'')---as he envisaged them at that time---do not work for Feynman integrands.}

To show that the $Y$-sets are useless for Feynman integration, suppose that they are open nbds. Then they contain open circles and so can be taken to be such circles. 

Then $(F_m -F)p\Delta x_1 \ldots \Delta x_n \in Y^j(I)$ implies
\[
\left|(F_m-F) p\Delta x_1 \ldots \Delta x_n\right| < \ve(I)
\]
for a collection of $\ve(I)>0$ with $(\D)\sum \ve(I) \leq \ve$ and
\[
|p| = \pi^{-\frac 12 n} \left((b_1-a)(b_2-b_1) \cdots (b_n-b_{n-1})\right)^{-\frac 12} :=s^{-1}
\]
in the Feynman sense. Hence for $m \geq k(\ve(I), \ve,f)$,
\[
\left|(F_m-F)\Delta x_1 \ldots \Delta x_n\right| < \ve(I)s.
\]
The set $X$ where $|F_m-F|\geq \frac 1r$ for some $m \geq k(\ve(I), \ve,f)$ has
\[
\left|\Delta x_1 \ldots \Delta x_n\right| < r\ve(I) s,\;\;\;\;\;\;\;V(\Delta x_1 \ldots \Delta x_n;\A;T;X) \leq r \ve \sup_\D s.
\]
As the supremum of $b_j-b_{j-1}$ tends to $0$ we will have
\[
V(\Delta x_1 \ldots \Delta x_n;\A;T;X)=0,\;\;\;\;\;\;\;\;\;
V(p\Delta x_1 \ldots \Delta x_n;\A;T;X)=0,
\]
so 
$
V(P;\A;T;X)=0$ and
$F_m=F$ almost everywhere.

\end{document}